\documentclass[final,onefignum,onetabnum,a4paper]{siamart250211}
\usepackage{hyperref}

\usepackage{amssymb,amsfonts}
\usepackage{amsmath}
\usepackage{graphicx,subfigure}
\usepackage{epstopdf}
\usepackage{algorithm,algorithmic}
\usepackage{color}
\usepackage{bm}
\usepackage{booktabs,array}
\usepackage{multirow,multicol}
\usepackage[T1]{fontenc}
\usepackage{geometry}

\ifpdf
  \DeclareGraphicsExtensions{.eps,.pdf,.png,.jpg}
\else
  \DeclareGraphicsExtensions{.eps}
\fi

\newsiamremark{remark}{Remark}
\newsiamremark{hypothesis}{Hypothesis}
\crefname{hypothesis}{Hypothesis}{Hypotheses}
\newsiamremark{assumption}{Assumption}
\crefname{assumption}{Assumption}{Assumption}
\newsiamthm{claim}{Claim}

\headers{A generalized BDF2 for the 2D (modified) Fisher-KPP equation}{L. Ge, Y.-L. Zhao, Q.-Y. Shu}

\title{A generalized BDF2 method for the two-dimensional (modified) Fisher-Kolmogorov-Petrovsky-Piskunov equation\thanks{Received by the editors DATE;
		\funding{This work is supported by the National Natural Science Foundation of China (12401536), and Sichuan Science and Technology Program (2024NSFSC0441).}}}

\author{Lei Ge\thanks{School of Mathematical Sciences, 
		Sichuan Normal University, Chengdu, Sichuan 610066, P.R. China 
		(\email{lgemaill@163.com}).}
		\and  Yong-Liang Zhao\thanks{Corresponding author. School of Mathematical Sciences, 
			Sichuan Normal University, Chengdu, Sichuan 610066, P.R. China 
			(\email{ylzhaofde@sina.com}).}
			\and  Qian-Yu Shu\thanks{School of Mathematical Sciences, 
				Sichuan Normal University, Chengdu, Sichuan 610066, P.R. China 
				(\email{shuqy@sicnu.edu.cn}).}
}

\usepackage{amsopn}

\ifpdf
\hypersetup{
pdftitle={A generalized BDF2 for the 2D (modified) Fisher-KPP equation},
pdfauthor={L. Ge, Y.-L. Zhao, Q.-Y. Shu}
}
\fi

\begin{document}

\maketitle

\begin{abstract}
	The Fisher-Kolmogorov-Petrovsky-Piskunov (Fisher-KPP) equation is a classical reaction-diffusion equation with broad applications in fields such as biology, chemistry, and physics.
	In this paper, an alternative second-order scheme is proposed by employing the generalized BDF2 method to approximate the two-dimensional (modified) Fisher-KPP equation. We  consider both  uniform and nonuniform time steps  in the proposed scheme. 
Stability and convergence of the scheme with  uniform time steps are proved. 
 Several numerical experiments demonstrate that our uniform and nonuniform schemes are robust and accurate.
\end{abstract}



\section{Introduction}
\label{sec1}
Reaction-diffusion equations play a central role in describing spatial propagation phenomena. 
Among them, the Fisher-Kolmogorov-Petrovsky-Piskunov (Fisher-KPP) equation, 
as a typical semilinear parabolic equation, 
has served as a fundamental model in studying population dynamics \cite{fisher1937wave} and chemical reactions  \cite{turing1990chemical} since its introduction in 1937. 
This equation has significant applications in various  scientific and engineering fields,
such as modeling epidermal wound healing \cite{habbal2014assessing,sherratt1990models}, 
describing nonlinear dynamics in nuclear reactors \cite{kastenberg1968stability}, 
and studying the propagation behavior of advantageous genes in multidimensional population genetics \cite{aronson1978multidimensional}.
Its two-dimensional form can be expressed as:
\begin{equation}
	\begin{cases} 
		\dfrac{\partial u(x,y,t)}{\partial t} = \kappa\left(\dfrac{\partial ^2 }{\partial x^2}+\dfrac{\partial^2 }{\partial y^2}\right)u(x,y,t) + f(u(x,y,t)), & (x, y) \in \Omega, \ 0 < t \leq T, \\ 
		u(x, y, 0) = u_0(x, y), & (x, y) \in \Omega, \\ 
		u(x, y, t) = \varphi(x, y, t), & (x, y) \in \partial \Omega, \ 0 \leq t \leq T,
	\end{cases}
	\label {eq1.1}
\end{equation}
where  $u(x,y,t)$ represents a certain density or concentration (e.g.,  population density,  chemical substance concentration), $u_0(x, y)$ and $\varphi(x, y, t)$ are known functions,
 $\kappa>0$ is the diffusion coefficient, 
  $\Omega = [x_L, x_R] \times [y_L,y_R] \subset \mathbb{R}^2$, 
$\partial \Omega$ is the boundary of $\Omega$  
and $f$ is a smooth function that satisfies 
\begin{equation}
	f(u) = 
	\begin{cases}
		K\*u(1-u^p), \\
		K_{pq}\*u^p(1-u)^q,
	\end{cases}
	\label{eq1.2}
\end{equation}
where $K$, $K_{pq}$, $p$ and $q$  are  all positive  constants. 
  The common form of Eq.~\cref{eq1.1} is the first type of  Eq.~\cref{eq1.2} \cite{qin2015two}, while the second type represents a modified nonlinear term \cite{kwak2023unconditionally}.
  For the modified Fisher-KPP equation, the coefficient $K_{pq}$ needs to satisfy the normalization condition:
\begin{equation}
	\int_0^1 K_{pq}u^p(1-u)^q du=1,\\
	K_{pq} = \frac{\Gamma(p + q + 2)}{\Gamma(p + 1)\Gamma(q + 1)} \text{ and } \Gamma(z) = \int_0^\infty t^{z-1} e^{-t} dt. 
	\label{eq1.3}
\end{equation}

Compared to the one-dimensional case of  Eq.~\cref{eq1.1}, solutions to the two-dimensional case  may exhibit more complex spatial patterns, e.g., spiral waves \cite{mikhailov1994complex}, which impose stricter requirements on the precision and stability of numerical approaches.
 Many excellent numerical approaches have been applied to Eq.~\cref{eq1.1}, such as finite difference methods \cite{macias2011bounded,macias2012explicit}, a wavelet collocation method 
 based on Chebyshev wavelets \cite{orucc2020efficient}, pseudospectral  methods \cite{balyan2020stability,chandraker2018implicit}, wavelet methods \cite{hariharan2009haar,mittal2006numerical}, least-squares finite element methods \cite{olmos2006pseudospectral}, and B-spline finite element methods \cite{daug2010numerical}.
 A class of stabilized and structure-preserving finite difference methods for the Fisher-KPP equation is introduced in \cite{deng2025analysis}.
   However,  the stabilizing coefficient $\beta$   needs to satisfy certain conditions, such as $\beta \geq \max \{(p - 1)b,~ b/2\}$. 
   These methods are only first-order accurate in time. 
   To mitigate the limitations of low time accuracy, some temporal second-order methods are chosen.
   In \cite{izadi2021second},  a two-parameter numerical scheme 
   is developed for the one-dimensional  Fisher-KPP equation, which achieves second-order accuracy in time.
  In \cite{kwak2023unconditionally}, the two-dimensional  modified Fisher–KPP equation is solved using an alternating direction implicit  and interpolation method, achieving second-order temporal accuracy.
  For more numerical methods for the (modified) Fisher-KPP equation,  see \cite{macias2019positivity,roquejoffre2018nontrivial,salako2021long,youssef2022solving}.
  
  Research on high-order time schemes  and higher-dimensions of Eq.~\cref{eq1.1} remain quite limited, even for the more complex case,  such as nonuniform time-stepping formulations.
 The second-order backward differentiation formula (BDF2) has attracted widespread attention from researchers due to its second-order accuracy, unconditional stability, and excellent dissipation properties.
   For instance, key progress has been made in the stability analysis of the variable-step BDF2 method \cite{li2022stability}.
    Using discrete orthogonal convolution  kernels and elliptic matrix norms,  the BDF2 method is stable on arbitrary time grids for initial value problems of ordinary differential equations. 
    Based on this,  Liao and Tang \cite{liao2020energy} adopt the nonuniform time-step BDF2 method, combined with the kernel recombination and complementary  technique and adaptive time-stepping strategy, to solve the Allen-Cahn equation.
    To ensure the maximum principle of the numerical solution, strict constraints on the time-step ratio are imposed, requiring, for instance,  that \( r_k < 1+\sqrt{2} \approx 2.414 \).
    In \cite{li2025structure}, they employ the uniform/variable-time-step weighted and shifted BDF2  methods combined with the scalar auxiliary variable  approach and stabilization techniques to solve the anisotropic Cahn-Hilliard model, where the weighting parameter satisfies $\theta \in \left[\frac{1}{2}, 1\right]$.
   For other partial differential equations and other  variable time-step methods, interested readers may refer to \cite{balyan2020stability,chandraker2018implicit,gu2021note,liao2022mesh,liao2021analysis,liu2023positivity,wang2021stability, zhao2021efficient}.
   
   In \cite{huang2024new}, the authors construct a new class of backward difference formula  based on the Taylor expansion at time \(t^{n+\beta}\) (where \(\beta>1\) is a tunable parameter).
   Then, an implicit-explicit (IMEX) scheme is proposed by using such the formula for parabolic equations.
   Inspired by \cite{huang2024new}, we employ  a generalized second-order backward differentiation formula (GBDF2) to derive an IMEX scheme for solving  Eq.~\cref{eq1.1}.
  Solutions of Eq.~\cref{eq1.1} often exhibit multiscale behaviors (e.g., rapid reaction phases followed by slow diffusion).
  Uniform time steps may fail to efficiently resolve such dynamics. 
     However, nonuniform grids adaptively refine temporal resolution, improving computational efficiency.
  Thus, inspired by the nonuniform time steps in \cite{zhang2021error}, we extend our scheme to nonuniform time steps.
   To address the computational bottleneck of our scheme, we propose a fast algorithm  based on the structure of the resulting systems.
 The main contributions of this work can be summarized as follows:
  \begin{enumerate}
  	\item[(a)]	We provide  an alternative second-order scheme for Eq.~\cref{eq1.1} and analyze its stability and convergence.  	
  	\item[(b)]  For higher solution resolution,  we extend the GBDF2 formula to variable time steps.
  	\item[(c)] We design a fast solver based on the  Discrete Sine Transform (DST) \cite{MORONEY2013304} to efficiently solve the resulting linear systems, addressing the computational bottleneck for large-scale problems.
\end{enumerate}

The rest  of this paper is  organized as follows: 
 \Cref{sec2} proposes a GBDF2 based IMEX (denoted as GBDF2-IMEX) scheme with uniform time steps for  Eq.~\cref{eq1.1}. Subsequently, the stability and convergence of the scheme are studied.  
 \Cref{sec3} extends the scheme to nonuniform time steps. 
  In \Cref{sec4}, numerical results are reported to validate the theoretical temporal and spatial convergence orders. 
 \Cref{sec5}  concludes the work and discusses potential future research directions.
\section{Uniform GBDF2-IMEX scheme for the Fisher-KPP equation}
\label{sec2}
 Let $V$ and $H$ be two Hilbert spaces, whose properties are similar to those in \cite{huang2024new}. 
 Let $V\subset H=H^{\prime}\subset V^{\prime}$, where $V$ densely and continuously embedded in $H$,  and $V^{\prime}$ is the dual space of $V$.
 $\mathcal{L}_h\colon V\to V^{\prime}$ is a negative definite, self-adjoint, linear operator.
 The inner product in the space $ H $ is denoted by $ (\cdot, \cdot) $,   and the corresponding norm is denoted by  $| \cdot |$. 
 The norm on $ V $ is represented by $ \| \cdot \| $, i.e., $\|u\|:=|\mathcal{L}_h^{1/2}u|=(\mathcal{L}_hu,u)^{1/2}$.  
 We assume that the nonlinear operator $\mathcal{G}$ is locally Lipschitz continuous within a ball $\mathcal{B}_{u(t)} := \{ v \in V : \| v - u(t) \| \leq 1 \}$ centered at the exact solution $u(t)$.
 The norm on the dual space \( V' \) is given by   
\begin{equation}
	\|v\|_* := \sup_{u \in V \setminus \{0\}} \frac{|\langle v, u \rangle|}{\|u\|}, \quad \forall v \in V'.
	\label{eq2.1}
\end{equation}
  Specifically, for all $v, \tilde{v} \in \mathcal{B}_{u(t)}$ and for all $t \in [0, T]$, the following inequality holds:
 \begin{equation}
 	\| \mathcal{G}(v) - \mathcal{G}(\tilde{v}) \|^{2}_{*} \leq \gamma \| v - \tilde{v} \|^{2} + \mu | v - \tilde{v} |^{2},
 	\label{eq2.2}
 \end{equation}
 where $\gamma$ is a non-negative constant and $\mu$ is an arbitrary constant.
Next, we derive the GBDF2-IMEX scheme with uniform time steps for Eq.~\cref{eq1.1} and analyze its stability and convergence.
\subsection{Discretization of the Fisher-KPP}\label{sec2.1}
For a given positive integer $M$, 
we discretize time interval $[0,T]$ as $0=t^0< t^1 < \cdots<t^M=T$, where $t^n = n \Delta t$  $(n=0,1,\ldots,M)$ with $\Delta t = T/M$. Let $h_x =\dfrac{x_R-x_L}{N_x} , h_y =\dfrac{y_R-y_L}{N_y}$, where $N_x$ and $N_y$ are two positive integers. Then, the spatial domain $\Omega$ is discretized as $\Omega_s =\left\{ {(x_i,y_j) | x_i = x_L + i h_x, y_j = y_L + j h_y, i = 0,1,\ldots,N_x, j=0,1,\ldots,N_y}\right\}$.

 The Taylor expansions of $u(t^{n+1})$, $u(t^n)$, and $u(t^{n-1})$ at time $t^{n+\beta}=t^n+\beta\Delta t$     $(\beta >1 \in \mathbb{R})$ are given in the following: 
\begin{equation*}
u(t^{n-1+i})=u(t^{n+\beta})+(-1+i-\beta)\Delta t\partial_{t}u(t^{n+\beta}) + \dfrac{1}{2}(-1+i-\beta)^2\Delta t^2\partial_{tt}u(t^{n+\beta})+ \mathcal{O}((\Delta t)^{3}),~\mathrm{for}~ 0 \leq i \leq 2 . 
\end{equation*}

From this, we can derive a second-order difference formula for approximating $ \partial_{t}u(t^{n + \beta})$. 
More precisely,
\begin{equation}
\dfrac{1}{\Delta t}	\sum_{q = 0}^{2}a_{q}u(t^{n-1+q}) = \partial_{t}u(t^{n + \beta}) + \mathcal{O}((\Delta t)^{2}),
		\label{eq2.3}
\end{equation}
where the coefficients  $a_{q} ~(q=0,1,2)$ are determined uniquely by solving the following linear system with the Vandermonde matrix 
\begin{equation*}
	\begin{bmatrix}1 & 1 &   1\\1-\beta &-\beta  &  -1-\beta \\(1-\beta )^2 & (-\beta)^2  & (-1-\beta  )^2\end{bmatrix}\begin{bmatrix}a_{2}\\a_{1}\\a_{0}\\\end{bmatrix}=\begin{bmatrix}0\\1\\0\\\end{bmatrix}. 
\end{equation*}
Solving this  system yields:
\[
a_2 = \frac{1 + 2\beta}{2}, \quad
a_1 = -2\beta, \quad
a_0 = \frac{2\beta - 1}{2}.
\]

We use the linear combination of $u(t^{n+1})$ and $u(t^{n})$ to approximate $u(t^{n+\beta})$ in the linear part of Eq.~\cref{eq1.1} as follows:
\begin{equation}
	b_{1}u(t^{n+1})+b_{0}u(t^n) = u(t^{n + \beta}) + \mathcal{O}((\Delta t)^{2}),
	\label{eq2.4}
\end{equation}
where $b_1 = \beta , b_0 = 1 - \beta$ are the unique solution of
\begin{equation*}
	\begin{bmatrix}1 & 1\\ 1-\beta &-\beta  \\ \end{bmatrix}
	\begin{bmatrix}b_{1}\\b_{0}\\ \end{bmatrix}=\begin{bmatrix}1\\0\\\end{bmatrix}.
\end{equation*}


 For obtaining a second-order approximation of $u(t^{n+\beta})$ in the nonlinear part of Eq.~\cref{eq1.1},
  a linear  combination of $u(t^n)$ and $u(t^{n-1})$ is used, i.e.,
\begin{equation}
	c_{1}u(t^{n})+c_{0}u(t^{n-1}) = u(t^{n + \beta}) + \mathcal{O}((\Delta t)^{2}).
		\label{eq2.5}
\end{equation}
Similarly, $ c_1 = 1 + \beta , c_0 = -\beta$ are the solution of
\begin{equation*}
	\begin{bmatrix}1 &    1\\-\beta &-1-\beta \\\end{bmatrix}\begin{bmatrix}c_{1}\\c_{0}\\\end{bmatrix}=\begin{bmatrix}1\\0\\\end{bmatrix}. 
\end{equation*}

We denote $ u_{i,j}^{n}$ as the numerical approximation to the exact solution $u(x_i, y_j,t^n) $ of  Eq.~\cref{eq1.1}.
Then, substituting Eqs.~\cref{eq2.3}--\cref{eq2.5} into Eq.~\cref{eq1.1} at \( u(x_i, y_j, t^{n+\beta}) \) yields the following temporal discretization scheme:
\begin{equation}
	\dfrac{1}{\Delta t}\sum_{q = 0}^{2}a_{q}u^{n -1 + q}=\mathcal{L}_h(b_{1}u^{n+1}+b_{0}u^n)+\mathcal{G}( c_{1}u^{n}+c_{0}u^{n-1}),
	\label{eq2.6}
\end{equation}
where $\mathcal{L}_h=\kappa(\frac{\partial^2}{\partial x^2}+\frac{\partial^2}{\partial y^2})$ and $u^n$ is the approximation of $u(x,y,t^n)$.
Applying a spatial central difference formula to the above equation, our second-order GBDF2-IMEX scheme for Eq.~\cref{eq1.1} is obtained as follows:
\begin{equation}
	\frac{a_2u^{n+1}_{i,j} + a_1 u^n_{i,j}+ a_0 u^{n-1}_{i,j}}{\Delta t} =  \delta_{xy}^2 \left(b_1 u^{n+1}_{i,j} + b_0 u^{n}_{i,j}  \right)
	+\mathcal{G}(c_{1}u_{i,j}^{n}+c_{0}u_{i,j}^{n-1}),
	\label{eq2.7}
\end{equation}
where $\delta_{xy}^2=\kappa(\delta_x^2+\delta_y^2)$  with $\delta_x^2$ and $\delta_y^2$ are the second-order central difference operators along the $x$- and $y$-directions, respectively. 
More precisely,
\begin{equation*}
	\delta_x^2 u_{i,j}^n = \frac{u_{i-1,j}^n - 2u_{i,j}^n + u_{i+1,j}^n}{h_x^2}, \quad
	\delta_y^2 u_{i,j}^n = \frac{u_{i,j-1}^n - 2u_{i,j}^n + u_{i,j+1}^n}{h_y^2}.
\end{equation*}

To convert the spatial discrete scheme into a matrix form, it is first necessary to define the relevant matrix symbols. 
Let $I_x$, $I_y$ and $I_{xy}$ be the identity matrices of sizes $N_x-1$, $N_y-1$ and $(N_x-1)(N_y-1)$, respectively.
By employing the Kronecker product $\otimes$, the spatially discretized matrix can be uniformly written as $L = I_y \otimes A_x + A_y \otimes I_x$, where 
\begin{equation*}
	\begin{aligned}
		A_x &= \frac{1}{h_x^2}
		\begin{bmatrix}
			-2 & 1   &        &           \\
			1  & -2  & 1      &           \\
			& \ddots & \ddots & \ddots   \\
			&        & 1     & -2 & 1   \\
			&        &       & 1  & -2
		\end{bmatrix} 
		\in \mathbb{R}^{(N_x-1) \times (N_x-1)}, \\
		A_y &= \frac{1}{h_y^2}
		\begin{bmatrix}
			-2 & 1   &        &           \\
			1  & -2  & 1      &           \\
			& \ddots & \ddots & \ddots   \\
			&        & 1     & -2 & 1   \\
			&        &       & 1  & -2
		\end{bmatrix}
		\in \mathbb{R}^{(N_y-1) \times (N_y-1)}.
	\end{aligned}
\end{equation*}
Then, the scheme \cref{eq2.7} can be rewritten into the matrix-vector multiplication form: 
\begin{equation}
	\left(\dfrac{1}{\Delta t}a_{2}I_{xy} - \kappa b_{1} L\right) \bm{u}^{n+1} = \left(-\dfrac{1}{\Delta t}a_{1}I_{xy} +\kappa b_{0} L\right) \bm{u}^{n} + \mathcal{G}(c_{0} \bm{u}^{n-1} + c_{1} \bm{u}^{n}) -\dfrac{1}{\Delta t} a_{0}\bm{u}^{n-1} +\kappa \bm{\xi}^n.
	\label{eq2.8}
\end{equation}
Here, $\bm{u}^n = [u_{1,1}^n, u_{2,1}^n, \ldots, u_{N_x-1,1}^n, u_{1,2}^n, u_{2,2}^n, \ldots, u_{N_x-1,2}^n, \ldots, u_{1,N_y-1}^n, u_{2,N_y-1}^n, \ldots, u_{N_x-1,N_y-1}^n]^{\top} ~(n=0,\ldots,M)$ is the approximate solution of Eq.~\cref{eq1.1} at $t^n$,
 and $\bm{\xi}^n
 =  \left[ \xi_{11}^n, \xi_{21}^n, \cdots, \xi_{N_x-1,1}^n, \xi_{12}^n, \cdots, \xi_{N_x-1,N_y-1}^n \right]^{\top}$  is the boundary vector,
where
\[\xi_{ij}^n = 
\begin{cases} 
	\frac{b_1u_{0,j}^{n+1}+b_0u_{0,j}^n}{h_x^2} & (i=1, 1 \leq j \leq N_y - 1) \quad (\text{left boundary contribution}) \\
	\frac{b_1u_{N_x,j}^{n+1}+b_0u_{N_x,j}^n}{h_x^2} & (i=N_x - 1, 1 \leq j \leq N_y - 1) \quad (\text{right boundary contribution}) \\
	\frac{b_1u_{i,0}^{n+1}+b_0u_{i,0}^n}{h_y^2} & (1 \leq i \leq N_x - 1, j = 1) \quad (\text{bottom boundary contribution}) \\
	\frac{b_1u_{i,N_y}^{n+1}+b_0u_{i,N_y}^n} & (1 \leq i \leq N_x - 1, j = N_y - 1) \quad (\text{top boundary contribution}) \\
	0 & (\text{internal nodes not adjacent to boundary, no boundary contribution})
\end{cases}\]

To fully utilize the structure of the coefficient matrix $\dfrac{1}{\Delta t}a_{2}I_{xy} - \kappa b_{1} L$, a fast algorithm is designed.
 We note that the tridiagonal matrices $A_x$ and $A_y$ admit known eigendecompositions in diagonal form:
\[
A_s = Q_s \Sigma_s Q_s^{\top}~ (s=x,y),\]
where $\Sigma_s =\mathrm{diag} (\bm{\lambda}_s)$ with  $ \lambda _{s,i} =\dfrac{-4}{h_s^2} \sin^2\left(\dfrac{\pi i}{2N_s}  \right), (Q_s)_{i j} = \sqrt{\dfrac{2}{N_s}} \sin\left(\dfrac{ij \pi}{N_s}\right)$ $( i,j=1,\ldots,N_s-1)$ is a symmetric and orthogonal matrix. 
Then, $L$ can be rewritten as 
\begin{equation*}
	L = I_y \otimes A_x + A_y \otimes I_x = (Q_x\otimes Q_y)[(I_y\otimes \Sigma_x+\Sigma_y\otimes I_x)](Q_x \otimes Q_y)^{\top} :=Q \Sigma Q^{\top}.
\end{equation*}
Hence, we have
\begin{equation*}
\left(\dfrac{1}{\Delta t}a_{2}I_{xy} - \kappa b_{1} L\right)^{-1}=Q\left[\dfrac{1}{\Delta t} a_2 I_{xy}-\kappa b_1\Sigma \right]Q^{\top}.
\end{equation*} 

Obviously, the core step of solving  \cref{eq2.8} is to compute the matrix-vector multiplication, i.e., $(	\dfrac{1}{\Delta t}a_{2}I_{xy} - \kappa b_{1} L)^{-1} v$ ($v$ is a vector).
 It can be done efficiently via Algorithm \ref{alg2.1}. 
 In this algorithm, \texttt{reshape}, \texttt{dst} and \texttt{idst} are build-in functions in MATLAB.
 In addition,  $\bm{1}_{N_y-1,1}$ and $\bm{1}_{1,N_x-1}$  are column vector and row vector, respectively.
 Their elements are all $1$.

\begin{algorithm}[ht]
	\caption{Compute $Z = \left(\dfrac{1}{\Delta t}a_{2}I_{xy} - \kappa b_{1} L\right)^{-1}v$}
	\label{alg2.1} 
	\begin{algorithmic}[1]
		\STATE Compute 
		\begin{align*}
			\lambda_x &=  \dfrac{-4}{ h_x^2} \left[ \sin^2 \frac{\pi}{2N_x}, \sin^2 \frac{2 \pi}{2N_x}, \ldots, \sin^2 \frac{(N_x - 1) \pi}{2N_x} \right], \\
		\lambda_y &=\dfrac{-4}{  h_y^2 } \left[ \sin^2 \frac{\pi}{2N_y}, \sin^2 \frac{2 \pi}{2N_y}, \ldots, \sin^2 \frac{(N_y - 1) \pi}{2N_y} \right],\\
			\lambda_L &=   \bm{1}_{N_y-1,1}  	\lambda_x + 	\lambda_y^{\top} \bm{1}_{1,N_x-1} ;
		\end{align*}
		
		\STATE { $\Lambda=\dfrac{1}{\Delta t}a_2I_{xy}-\kappa b_1 \Lambda_L$};
	   
	   \STATE {$Z$ =\texttt {reshape}($v,Nx - 1,Ny - 1$)};
      \STATE {$Z$ = \texttt{dst}($Z$);  $Z$ = \texttt{dst}($Z^\top$)};
   \STATE {$Z$ = $Z./ \Lambda$};                 
\STATE	{$Z$ = \texttt{idst}($Z$); ~  $Z$ = \texttt{idst}($Z^\top$)};           

\STATE {	$Z$ = \texttt{reshape}($Z$,($Nx-1)(Ny-1)$,1);} 
	\end{algorithmic}
\end{algorithm}

\subsection{Stability  and convergence of the uniform GBDF2-IMEX scheme}\label{sec2.2}
 To simplify the presentations, we denote
 \begin{equation*}
 	A(u^{n+1}) = \sum_{q=0}^2a_{q}u^{n-1+q}, \quad B(u^{n+1}) = b_{1}u^{n+1}+b_{0}u^n, \quad C(u^{n}) = c_{1}u^n+c_{0}u^{n-1}.
 	\label{2-6}
 \end{equation*}
 
 To facilitate the subsequent proof, we introduce a splitting of  $B(u^{n+1})$ 
 \begin{equation}
 	B(u^{n+1})=\eta C(u^{n+1})+D(u^{n+1}),\quad \eta=\dfrac{1-\beta}{\beta} <0 \quad~\mathrm{and}\quad D(u^{n+1})=(2-2\beta)u^n+(2\beta-\dfrac{1}{\beta})u^{n+1}.
 	\label{eq2.9}
 \end{equation}
 We first consider the  convergence  for  Eq.~\cref{eq2.6}, which can be rewritten as
 \begin{equation}
 	\frac{A(u^{n+1})}{\Delta t}=\mathcal{L}_hB(u^{n+1})+\mathcal{G}C(u^{n}).
 	\label{eq2.10}
 \end{equation}
 	
 Using these notations, we define the following truncation errors:
 \begin{align*}
 	E^{n+1} &:= \Delta t u_t(t^{n+\beta}) - A(u(t^{n+1})),
 	\\
 	R^{n+1} &:= u(t^{n+\beta}) - B(u(t^{n+1})), 
 	\\
 	P^n &:= u(t^{n+\beta}) - C(u(t^n)).
 \end{align*}
 
 Based on the Taylor expansion and the integral remainder term, and under the regularity requirements for $u$, 
 we derive 
 \begin{equation*}
 	\begin{aligned}
 		E^{n+1} 
 		&= \dfrac{1}{2}a_2\int_{t^{n+1}}^{t^{n+\beta}}(t^{n+1}-s)^2u^{(3)}(s) ds 
 		+ \dfrac{1}{2}a_1\int_{t^{n}}^{t^{n+\beta}}(t^{n}-s)^2u^{(3)}(s) ds  + \dfrac{1}{2}a_0\int_{t^{n-1}}^{t^{n+\beta}}(t^{n-1}-s)^2u^{(3)}(s) ds \\
 		&= \dfrac{1}{2}\sum_{q=0}^{2}a_q\int_{t^{n-1+q}}^{t^{n+\beta}}(t^{n-1+q}-s)^2u^{(3)}(s)ds
 		\leq  \dfrac{\hat{C_1}}{2}\sum_{q=0}^{2}a_q\int_{t^{n-1+q}}^{t^{n+\beta}}(t^{n-1+q}-s)^2ds,
 	\end{aligned}
 \end{equation*}
where $\hat{C_1}$ is a constant such that $|u^{(3)}(s)| \leq \hat{C_1}$ for all $s$ in the relevant time interval.
 Analogously,  we obtain
 \begin{equation}
 	\vert E^{n + 1}\vert^2 \leq \hat{C}(\Delta t)^{6}, \quad \Vert R^{n + 1}\Vert^2 \leq \hat{C}(\Delta t)^{4}, \quad \Vert P^n\Vert^2 \leq \hat{C}(\Delta t)^{4}, \quad \forall n + 1 \leq \frac{T}{\Delta t},
 	\label{eq2.11}
 \end{equation}
where  $\hat{C}$ is a constant independent of $\Delta t$.  
Unless otherwise specified, we denote by $\hat{C}$ a constant whose value varies as needed in practice.

We denote  $e^m =  {u}^m - u(t^m)$, which  measures the numerical error at time $t^m(m=0,1,\ldots,M)$, where
$u(t^m)$ being the exact solution of Eq.~\cref{eq1.1} at $t=t^m$.
Obviously, $e^0=0$.
 At  $t^{n+\beta}$,  Eq.~\cref{eq1.1} becomes
\begin{equation}
	u_{t}(t^{n+\beta})=\mathcal{L}_hu(t^{n+\beta})+\mathcal{G}[u(t^{n+\beta})].
	\label{eq2.12}
\end{equation}

This will be used later when employing the multiplier.
The proof follows a similar structure to \cite[Theorem 3]{huang2024new}. 
For completeness, a concise version is provided below. 

\begin{theorem}
	\label{theorem2.1}
	Assume Eq.~\cref{eq2.2} and the solution of Eq.~\cref{eq1.1}  is sufficiently smooth such that Eq.~\cref{eq2.11} is true, and the following stability condition  	is satisfied
	\begin{equation}
		-\eta-\sqrt{\gamma}\geq\rho>0,
		\label{eq2.13}
	\end{equation}
 where $\rho$ is a constant. Given $u^{0}=u(0)\in V$, we assume   $u^{1}$ is computed with a proper initialization procedure such that
	\begin{equation}
		|u^{1}-u(t^{1})|^{2},~\|u^{1}-u(t^{1})\|^{2}\leq \hat{C}(\Delta t)^{4}~\mathrm{and}~ C(u^{1})\in\mathcal{B}_{u(t^{1+\beta})}.
		\label{eq2.14}
	\end{equation}
	Then, for $\Delta t$ sufficiently small and $\beta>1$, we conclude that
	\begin{equation}
		C(u^{n+1})\in\mathcal{B}_{u(t^{n+1+\beta})},\quad\forall n+1\leq\frac{T}{\Delta t},
		\label{eq2.15}
	\end{equation}
	and
\begin{equation}
	g|e^{m+1}|^{2}+\dfrac{\rho}{2} \Delta t\sum_{n=1}^{m}\|C(e^{n+1})\|^{2}\leq \hat{C} \exp(1-\hat{C}(\Delta t)^{-1}T)(\Delta t)^4,\quad\forall\,n\leq m.
	\label{eq2.16}
\end{equation}
	where $g$ is a positive constant and $\hat{C}$ is independent of $\Delta t$.
\end{theorem}

\begin{proof}
Suppose that Eq.~\eqref{eq2.16} holds for all $n \leq m-1$, and assume we already have
\begin{equation}
	C(u^{n}) \in \mathcal{B}_{u(t^{n+\beta})} \quad \text{for all } n \leq m.
	\label{eq2.15}
\end{equation}
It then remains to prove that Eq.~\eqref{eq2.16} also holds for all $n \leq m$, and that
\begin{equation}
	C(u^{m+1}) \in \mathcal{B}_{u(t^{m+1+\beta})}.
	\label{eq2.18}
\end{equation}
	
	Subtracting  Eq.~\cref{eq2.12} from Eq.~\cref{eq2.10} and multiplying by $\Delta t$, we arrive at
	\begin{equation*}
		A(u^{n+1})-\Delta t u_t(t^{n+\beta}) = \Delta t (\mathcal{L}_hB(u^{n+1})-\mathcal{L}_hu(t^{n+\beta})) + \Delta t \left(\mathcal{G}[C(u^n)]-\mathcal{G}[u(t^{n+\beta})]\right).
	\end{equation*}
	
	To introduce the error, we derive
\begin{align*}
	&A(u^{n+1}) - A(u(t^{n+1})) - \left(\Delta t u_t(t^{n+\beta}) - A(u(t^{n+1}))\right) \\
	&\quad = \Delta t \left(\mathcal{L}_hB(u^{n+1}) - \mathcal{L}_hB(u(t^{n+1}))\right) \quad - \Delta t \left(\mathcal{L}_hu(t^{n+\beta}) - \mathcal{L}_hB(u(t^{n+1}))\right) \\
	&\quad + \Delta t \left(\mathcal{G}[C(u^n)] - \mathcal{G}[u(t^{n+\beta})]\right).
\end{align*} 

Substituting the error expressions and rearranging yields the error equation: 
	\begin{equation}
		A(e^{n+1})-\Delta t\mathcal{L}_hB(e^{n+1})=\Delta t\big{(}\mathcal{G}[C(u^{n})]-\mathcal{G}[u(t^{n+\beta})]\big{)}+E^{n+1}-\Delta t\mathcal{L}_hR^{n+1},
		\label{eq2.19}
	\end{equation}

 We split $\mathcal{G}[C(u^{n})]-\mathcal{G}[u(t^{n+\beta})]$ as
	\begin{align*}
		\mathcal{G}[C(u^{n})]-\mathcal{G}[u(t^{n+\beta})]&=\big{(}\mathcal{G}[C(u^{n})]-\mathcal{G}[C(u(t^{n}))]\big{)}+\big{(}\mathcal{G}[C(u(t^{n}))]-\mathcal{G}[u(t^{n+\beta})]\big{)}\nonumber\\
		&=:T^{n}_{1}+T^{n}_{2}.
	\end{align*}
	
	We denote 
	\[	(\mathcal{L}_hC(e^{n+1}),C(e^{n+1}))=\|C(e^{n+1})\|^2. 
	\]
	
	Now, we take the inner product of  Eq.~\eqref{eq2.19} with $C(e^{n+1})$.
	 On the left-hand side, we apply the splitting  $B(e^{n+1})$ for Eq.~\eqref{eq2.9} which gives:
	\begin{align}
		&\big{(}A(e^{n+1}),C(e^{n+1})\big{)}-\Delta t\eta\|C(e^{n+1})\|^{2}-\Delta t\big{(}\mathcal{L}_hD(e^{n+1}),C(e^{n+1})\big{)}\nonumber\\
		&=\Delta t\big{(}T^{n}_{1},C(e^{n+1})\big{)}+\Delta t\big{(}T^{n}_{2},C(e^{n+1})\big{)}+\big{(}E^{n+1},C(e^{n+1})\big{)}-\Delta t\big{(}\mathcal{L}_hR^{n+1},C(e^{n+1})\big{)}.
		\label{eq2.20}
	\end{align}
	
    With $\Delta t$ sufficiently small, we obtain $C(u(t^{n})) \in \mathcal{B}_{u(t^{n+\beta})}$ from Eq.~\eqref{eq2.5}. 
    Subsequently, applying Eq.~\eqref{eq2.2} and the Young inequality with an arbitrary $\varepsilon > 0$ to bound the right-hand side of Eq.~\eqref{eq2.20},  it follows that
	\begin{align}
		\big{|}\big{(}T^{n}_{1},C(e^{n+1})\big{)}\big{|}&\leq\|T^{n}_{1}\|_{*}\|C(e^{n+1})\|\leq\frac{\varepsilon}{2}(\gamma\|C(e^{n})\|^{2}+\mu|C(e^{n})|^{2})+\frac{1}{2\varepsilon}\|C(e^{n+1})\|^{2},
		\label{eq2.21}
	\end{align}

Similarly, we bound the  term:
	\begin{align}
		\big{|}\big{(}T^{n}_{2},C(e^{n+1})\big{)}\big{|}&\leq\|T^{n}_{2}\|_{*}\|C(e^{n+1})\|\leq\frac{1}{\rho}(\gamma\|P^{n}\|^{2}+\mu|P^{n}|^{2})+\frac{\rho}{4}\|C(e^{n+1})\|^{2},\nonumber\\
		&\leq \hat{C}(\Delta t)^{4}+\frac{\rho}{4}\|C(e^{n+1})\|^{2}.
		\label{eq2.22}
	\end{align}
	Here, $\hat{C}$ is  independent of $\Delta t$.
	
	For the truncation error terms, we proceed as follows:
	\begin{equation}
		\big{(}E^{n+1},C(e^{n+1})\big{)}\leq\frac{1}{2\Delta t}|E^{n+1}|^{2}+\frac{\Delta t}{2}|C(e^{n+1})|^{2}\leq \hat{C}(\Delta t)^{5}+\frac{\Delta t}{2}|C(e^{n+1})|^{2},
		\label{eq2.23}
	\end{equation}
	and
	\begin{equation}
		\big{(}\mathcal{L}_hR^{n+1},C(e^{n+1})\big{)}\leq\frac{1}{\rho}\|R^{n+1}\|^{2}+\frac{\rho}{4}\|C(e^{n+1})\|^{2}\leq \hat{C}(\Delta t)^{4}+\frac{\rho}{4}\|C(e^{n+1})\|^{2}.
		\label{eq2.25}
	\end{equation}
	
	Denote \(\phi (e)^{n+1}:=(e^{n-1},e^n,e^{n+1})^{T}\).
	For the left-hand side terms of Eq.~\cref{eq2.20}, we make use of a symmetric positive definite matrix $ G =(g_{ij})\in \mathbb{R}^{2 \times 2} $ (see \cite{huang2024new} for details) and a positive constant $h$ such that
\begin{align}
	\left( A( e^{n+1} ), C( e^{n+1} )\right) &\geq \sum_{i,j=1}^{2} g_{ij} \left( e^{n+i-1}, e^{n+j-1} \right) - \sum_{i,j=1}^{2} g_{ij} \left( e^{n+i-2}, e^{n+j-2} \right) \notag \\
	&=: \left| \phi (e)^{n+1} \right|_{G}^{2} - \left| \phi (e)^{n} \right|_{G}^{2}. \\
	-\left( \mathcal{L}_hD( e^{n+1} ), C( e^{n+1} )\right) &\geq h \left( \mathcal{L}_h e^{n+1}, e^{n+1} \right) - h \left( \mathcal{L}_h e^{n}, e^{n} \right) \notag \\
	&=: \left| \phi (e)^{n+1} \right|_{H}^{2} - \left| \phi (e)^{n} \right|_{H}^{2}.
	\label{eq2.26}
\end{align}

Combining Eqs.~\cref{eq2.20}-\cref{eq2.26}, we obtain
\begin{equation*}
	\begin{aligned}
		&\left| \phi (e)^{n+1} \right|_{G}^{2} - \left| \phi (e)^{n} \right|_{G}^{2}
		- \Delta t \eta \|C(e^{n+1})\|^2
		+ \Delta t \left( \left| \phi (e)^{n+1} \right|_{H}^{2} - \left| \phi (e)^{n} \right|_{H}^{2} \right) \\
		&\leq (\dfrac{1}{2\varepsilon}+\dfrac{\rho}{2}) \Delta t \|C(e^{n+1})\|^2
		+ \dfrac{\varepsilon \gamma }{2} \Delta t \|C(e^n)\|^2  + \dfrac{\varepsilon \mu}{2} \Delta t |C(e^n)|^2
		+ \dfrac{\Delta t}{2} |C(e^{n+1})|^2
		+ \hat{C}(\Delta t)^4.
		\label{2.29}
	\end{aligned}
\end{equation*}

 Summing $n$ from $1$ to $m$  and doing some cancellations,  we arrive at
\begin{equation*}
	\begin{aligned}
		&\left| \phi (e)^{m+1} \right|_{G}^{2} 
		-(\eta +\dfrac{1}{2\varepsilon}+\dfrac{\rho}{2}+\dfrac{\varepsilon \gamma}{2})\Delta t \sum_{n=1}^{m}\|C(e^{n+1})\|^2
		+ \Delta t \left| \phi (e)^{m+1} \right|_{H}^{2} \\
		&\leq \left| \phi (e)^{1} \right|_{G}^{2} + \Delta t\left| \phi (e)^{1} \right|_{H}^{2} 
		+ (\dfrac{\varepsilon \mu}{2}+\dfrac{1}{2})\Delta t\sum_{n=1}^{m} |C(e^{n+1})|^2 \\
		&\quad +\dfrac{\varepsilon \gamma}{2} \Delta t\|C(e^1)\|^2 + \dfrac{\varepsilon \mu}{2}\Delta t |C(e^{1})|^2
		+ \hat{C}(\Delta t)^4.
	\end{aligned}
\end{equation*}
	
Taking  $\varepsilon=\dfrac{1}{\sqrt{\gamma}}$ and  $-\eta-\sqrt{\gamma}\geq\rho>0$ based on the condition Eq.~\cref{eq2.13},  it yields
		\begin{equation}
			\begin{aligned}
				&\left| \phi (e)^{m+1} \right|_{G}^{2} +
				\dfrac{\rho}{2}\Delta t \sum_{n=1}^{m}\|C(e^{n+1})\|^2
				+ \Delta t \left| \phi (e)^{m+1} \right|_{H}^{2} \\
				&\leq \left| \phi (e)^{1} \right|_{G}^{2} + \Delta t\left| \phi (e)^{1} \right|_{H}^{2} 
				+ (\dfrac{\varepsilon \mu}{2}+\dfrac{1}{2})\Delta t\sum_{n=1}^{m} |C(e^{n+1})|^2 \\
				&\quad + \dfrac{\varepsilon \mu}{2}\Delta t |C(e^{1})|^2+\dfrac{\varepsilon \gamma}{2} \Delta t\|C(e^1)\|^2
				+ \hat{C}(\Delta t)^4.
				\label{eq2.27}
			\end{aligned}
		\end{equation}
		
	Let $ g>0 $ be the smallest eigenvalue of the symmetric positive definite matrix  $ G \in \mathbb{R}^{2 \times 2} $.
	This yields the lower bound:
		\begin{equation}
			|\phi (e)^{m+1}|_{G}^{2}\geq g|e^{m+1}|^{2}.
			\label{eq2.28}
		\end{equation}
		
		Furthermore, based on the initialization condition in Eq.~\cref{eq2.14}, this gives  $\hat{C}$  such that the following bounds:
		\begin{equation}
			|\phi (e)^{1}|_{G}^{2}\leq \hat{C}\sum_{i=0}^{1}|e^{i}|^{2} =\hat{C}|e^1|^2\leq \hat{C} (\Delta t)^4 \qquad \mathrm{and}
			\label{eq2.29}
		\end{equation}
		\begin{equation}
			\Delta t\|\phi (e)^{1}\|_{H}^{2}\leq \hat{C}\Delta t\sum_{i=0}^{1}\|e^{i}\|^{2}=\hat{C} \Delta t \|e^1\|^2\leq \hat{C} (\Delta t)^4.
			\label{eq2.30}
		\end{equation}
		
	Given the linear form of \( C(e^{n+1}) \), its square can be bounded by the sum of squares of \( e^n \) and \( e^{n+1} \) via the inequality \( (c_1 e^{n+1} + c_0 e^n)^2 \leq 2(c_1^2 |e^{n+1}|^2 + c_0^2 |e^n|^2) \). This leads to
	\begin{equation}
			\left(\frac{\mu}{2\sqrt{\gamma}}+\frac{1}{2}\right)\Delta t\sum_{n=1}^{m}|C(e^{n+1})|^{2} \leq \hat{C}\Delta t\sum_{n=1}^{m}(|e^{n}|^{2}+|e^{n+1}|^{2}) \leq \hat{C}\Delta t\sum_{q=0}^{m+1}|e^{q}|^{2}.
		\label{eq2.31}
	\end{equation}
	
		We notice that $C(e^1) = c_1e^1$ and Eq.~\cref{eq2.14}, the term
		$\dfrac{\varepsilon \mu}{2}\Delta t |C(e^{1})|^2+\dfrac{\varepsilon \gamma}{2} \Delta t\|C(e^1)\|^2$ 
		in Eq.~\cref{eq2.27} can be bounded by $\hat{C}(\Delta t)^4$ (i.e., 
		$\dfrac{\varepsilon \mu}{2}\Delta t |C(e^{1})|^2+\dfrac{\varepsilon \gamma}{2} \Delta t\|C(e^1)\|^2 \leq \hat{C}(\Delta t)^4$).
		 
		Combining Eqs.~\cref{eq2.27}-\cref{eq2.31}, we arrive at
		\begin{equation}
			g|e^{m+1}|^{2}+\dfrac{\rho}{2} \Delta t\sum_{n=1}^{m}\|C(e^{n+1})\|^{2}\leq \hat{C}\Delta t\sum_{q=0}^{m+1}|e^{q}|^{2}+\hat{C}(\Delta t)^{4},\quad\forall\,n\leq m.
			\label{eq2.32}
		\end{equation}

Therefore, an application of the discrete Gronwall inequality \cite[Lemma 2]{huang2024new} leads to
\begin{equation}
	g|e^{m+1}|^{2}+\frac{\rho}{2}\Delta t\sum_{q=1}^{m+1}|C(e^{q})|^{2}\leq \hat{C}\exp\big{(}(1-\hat{C}\Delta t)^{-1}T\big{)}(\Delta t)^{4},\quad\forall m+1\leq\frac{T}{\Delta t},
	\label{eq2.33}
\end{equation}
where $\hat{C}$ is  independent of $\Delta t$. This establishes Eq.~\eqref{eq2.16}.

Combining Eq.~\eqref{eq2.33} with \eqref{eq2.5} yields
\begin{align}
	\|C(u^{m+1})-u(t^{m+1+\beta})\|^{2}&\leq 2\|C(u^{m+1})-C(u(t^{m+1}))\|^{2}+2\|C(u(t^{m+1}))-u(t^{m+1+\beta})\|^{2}\nonumber\\
	&\leq 2\|C(e^{m+1})\|^{2}+\mathcal{O}((\Delta t)^{4})\nonumber\\
	&\leq\hat{C}(\Delta t)^{3}.
	\label{eq2.34}
\end{align}
According to Eq.\cref{eq2.34}, Eq.~\eqref{eq2.18} is proved for sufficiently small $\Delta t$.
	\end{proof}
	\begin{remark}
		From Eq.~\eqref{eq2.33},
		we find that $g|e^{m+1}|^2 \leq \hat{C} (\Delta t)^4$, i.e., $|e^{m+1}| \leq \hat{C}(\Delta t)^2$.
		 For the spatial discretization by the central difference scheme employed in this paper, the local truncation error in space is of order $\mathcal{O}(h_x^2 + h_y^2)$. 
		 Integrating the aforementioned temporal error estimate with the spatial discretization accuracy, the overall error of the fully discrete solution $|u_{i,j}^{n} - u(x_i, y_j, t^n)|$ can be bounded via the triangle inequality as follows
		\[
		|u_{i,j}^{n} - u(x_i, y_j, t^n)| \leq |u_{i,j}^{n} - u^{n}|+ |u^{n} - u(x_i, y_j, t^n)|  \leq \hat{C} \left( (\Delta t)^2 + h_x^2 + h_y^2 \right),
		\]
	where $u^{n}$ is the solution of  Eq.~\cref{eq2.6}  and $\hat{C}$ is independent of $\Delta t$. 
	 This establishes the second-order convergence in both time and space for the fully discrete Eq.~\cref{eq2.7}.
	\end{remark}

	On this basis, we next prove the stability of Eq.~\cref{eq2.7}.
	We define the error $\tilde{e}^n_{ij}=u^n_{ij}-U^n_{ij}(n=0,\ldots,M,i=1,\ldots,N_x,j=1,\ldots,N_y)$, 
	where $U^n_{ij}$ is the approximation of $u^{n}_{ij}$ given in Eq.~\cref{eq2.7}. 
	Obviously, $\tilde{e}^0_{ij} = 0$.
	Then, we denote $\tilde{\bm{e}}^n=\bm{u}^n - U^n$, where $U^n = [U_{1,1}^n, \ldots, U_{N_x,1}^n, \ldots, U_{N_x,N_y}^n]^{\top }$.
 \begin{theorem}
 Under the assumptions of  Theorem~\ref{theorem2.1} for $\beta>1$, the scheme \cref{eq2.10} is stable in the sense that
 	\[
 	g|\tilde{\bm{e}}^{m+1}|^{2} -  \dfrac{\eta \Delta t}{2}\|C(\tilde{\bm{e}}^{m+1})\|^2 \leq \hat{C}(|\tilde{\bm{e}}^{1}|^{2}+\Delta t \|\tilde{\bm{e}}^{1}\|^2)- \dfrac{\hat{C} \Delta t}{2\eta}\sum_{n=1}^{m}\left( \|C(\tilde{\bm{e}}^{n})\|^2+\|C(\tilde{\bm{e}}^{n})\|\right).
 	\]
 The $g$  and $\eta$ are defined respectively in Eqs.~\cref{eq2.28} and  \cref{eq2.9}.
\end{theorem}

\begin{proof}
	 From Eq.~\eqref{eq2.7}, we arrive at
	\begin{equation}
		\dfrac{1}{\Delta t}A(\tilde{\bm{e}}^{n+1})=\delta_{xy}^2B(\tilde{\bm{e}}^{n+1})+\mathcal{G}(C(u^{n}))-\mathcal{G}(C(U^{n})).
		\label{eq2.35}
	\end{equation}
	
We consider the inner product of Eq.~\cref{eq2.35} with $\Delta t C(\widetilde{\bm{e}}^{n+1})$ and split $B(\widetilde{\bm{e}}^{n+1})$, we  conclude
\begin{equation}
	\begin{aligned}
		\big(A(\tilde{\bm{e}}^{n+1}), C(\tilde{\bm{e}}^{n+1})\big)
		= &\ \Delta t \eta \|C(\tilde{\bm{e}}^{n+1})\|^2 
		+ \Delta t \big(\delta_{xy}^2 D(\tilde{\bm{e}}^{n+1}), C(\tilde{\bm{e}}^{n+1})\big) \\
		&+ \Delta t \big(\mathcal{G}(C(u^{n})) - \mathcal{G}(C(U^n)), C(\tilde{\bm{e}}^{n+1})\big).
	\end{aligned}
	\label{eq2.36}
\end{equation}

Let \( \phi (\tilde{\bm{e}})^{n+1}:= (\tilde{\bm{e}}^{n-1}, \tilde{\bm{e}}^n, \tilde{\bm{e}}^{n+1})^T \). 
For  Eq.~\cref{eq2.36},  similar to Theorem~\ref{theorem2.1}, we make use of the same symmetric positive definite matrix \( G\) and the positive constant  \( h  \) such that
	\begin{equation}
\begin{aligned}
	\left(A(\tilde{\bm{e}}^{n+1}),C(\tilde{\bm{e}}^{n+1})\right) &\geq \sum_{i,j=1}^{2} g_{ij} \left( \bm{e}^{n+i-1}, \bm{e}^{n+j-1} \right) - \sum_{i,j=1}^{2} g_{ij} \left( \bm{e}^{n+i-2}, \bm{e}^{n+j-2} \right) \\
	&=: \left| \phi (\tilde{\bm{e}})^{n+1} \right|_{G}^{2} - \left| \phi (\tilde{\bm{e}})^{n} \right|_{G}^{2}.
\end{aligned}
\label{eq2.37}
	\end{equation}
		\begin{equation}
\begin{aligned}
	-\left(\delta_{xy}D(\tilde{\bm{e}}^{n+1}),C(\tilde{\bm{e}}^{n+1})\right) &\geq h(\delta_{xy}^2(\tilde{\bm{e}}^{n+1}),\tilde{\bm{e}}^{n+1})-h(\delta_{xy}^2(\tilde{\bm{e}}^n,\tilde{\bm{e}}^n)) \\
	&=: \| \phi (\tilde{\bm{e}})^{n+1} \|_{H}^{2} - \| \phi (\tilde{\bm{e}})^{n} \|_{H}^{2}.
\end{aligned}
\label{eq2.38}
	\end{equation}
	
	The nonlinear term on the right-hand side of Eq.~\cref{eq2.36} is bounded as follows
\begin{align}
	\bigl(\mathcal{G}(C(u^{n}))-\mathcal{G}(C(U^n)),C(\tilde{\bm{e}}^{n+1})\bigr)
	&\leq \|\mathcal{G}(C(u^{n}))-\mathcal{G}(C(U^n))\|_{*}\,\|C(\tilde{\bm{e}}^{n+1})\| \nonumber \\
	&\leq \frac{1}{2\varepsilon}\|\mathcal{G}(C(u^{n}))-\mathcal{G}(C(U^n))\|_{*}^2 
	+ \frac{\varepsilon}{2}\|C(\tilde{\bm{e}}^{n+1})\|^2 \nonumber \\
	&\leq \frac{1}{2\varepsilon}\Bigl(\gamma\|C(\widetilde{\bm{e}}^{n})\|^2
	+ \mu|C(\tilde{\bm{e}}^{n})|^2\Bigr)
	+ \frac{\varepsilon}{2}\|C(\tilde{\bm{e}}^{n+1})\|^2,
	\label{eq2.39}
\end{align}
where  the condition \eqref{eq2.2} is applied.

 From Eqs.~\cref{eq2.36}-\cref{eq2.39},  let $\varepsilon=-\eta>0$, we get
\begin{equation}
	\begin{aligned}
		&|\phi (\tilde{\bm{e}})^{n+1}|_G^2 - |\phi (\tilde{\bm{e}})^n|_G^2 + \Delta t \left(\|\phi (\tilde{\bm{e}})^{n+1}\|_H^2 - \|\phi (\tilde{\bm{e}})^n\|_H^2\right) \\
		&\quad \leq \frac{\Delta t \eta}{2} \|C(\tilde{\bm{e}}^{n+1})\|^2 - \frac{\Delta t}{2\eta} \left(\gamma \|C(\tilde{\bm{e}}^{n})\|^2 + \mu |C(\tilde{\bm{e}}^{n})|^2\right).
	\end{aligned}
	\label{eq2.40}
\end{equation}

Summing $n$ from $1$ to $m$  and doing some cancellations for Eq.~\cref{eq2.40}, it follows that
\begin{multline}
	|\phi (\tilde{\bm{e}})^{m+1}|_G^2 + \Delta t \|\phi (\tilde{\bm{e}})^{m+1}\|_H^2
	\leq |\phi (\tilde{\bm{e}})^1|_G^2 + \Delta t \|\phi (\tilde{\bm{e}})^{1}\|_H^2 \\
	+ \dfrac{\Delta t \eta}{2} \sum_{n=1}^m \|C(\tilde{\bm{e}}^{n+1})\|^2
	- \dfrac{\Delta t}{2 \eta} \sum_{n=1}^{m} \left(\gamma \|C(\widetilde{\bm{e}}^{n})\|^2 + \mu |C(\tilde{\bm{e}}^{n})|^2\right).
	\label{eq2.41}
\end{multline}
 
 Further, it can be known from Eq.~\cref{eq2.28} that
\begin{equation}
	|\phi (\tilde{\bm{e}})^{m+1}|_{G}^{2}\geq g|\tilde{\bm{e}}^{m+1}|^{2},
	\label{eq2.42}
\end{equation}
 where $g$ is given in Eq.~\cref{eq2.28}.

 We can choose large enough  \(\hat{C}\) such that
\begin{equation}
	|\phi (\tilde{\bm{e}})^{1}|_{G}^{2}\leq \hat{C}\sum_{i=0}^{1}|\tilde{\bm{e}}^{i}|^{2}=\hat{C}|\tilde{\bm{e}}^{1}|^{2}\qquad \mathrm{and}  \qquad \Delta t \|\phi (\tilde{\bm{e}})^{1}\|_H^2 \leq \hat{C} \Delta t \sum_{i=0}^{1}\|\bm{e}^i\|^2=\hat{C}\|\tilde{\bm{e}}^{1}\|^{2}.
	\label{eq2.43}
\end{equation}

Combining Eqs.~\cref{eq2.41}-\cref{eq2.43}, we arrive at
	\[
g|\tilde{\bm{e}}^{m+1}|^{2} -  \dfrac{\eta \Delta t}{2}\|C(\tilde{\bm{e}}^{m+1})\|^2 \leq \hat{C}(|\tilde{\bm{e}}^{1}|^{2}+\Delta t \|\tilde{\bm{e}}^{1}\|^2)- \dfrac{\hat{C} \Delta t}{2\eta}\sum_{n=1}^{m}\left( \|C(\tilde{\bm{e}}^{n})\|^2+\|C(\tilde{\bm{e}}^{n})\|\right).
\]
\end{proof}

The uniform time-step scheme \cref{eq2.7} possesses second-order accuracy and favorable stability. However, for solving Eq.~\cref{eq1.1} with its multiscale transient behaviors—such as alternating rapid reaction and slow diffusion phases—uniform time steps may not be optimal. Nonuniform time steps provide greater flexibility: they can be refined during rapid changes for accuracy and relaxed during gradual phases for efficiency. Based on this advantage, we will extend the GBDF2-IMEX scheme  \cref{eq2.7} to nonuniform time steps in the next section.

\section{Nonuniform GBDF2-IMEX discretization scheme}
\label{sec3}
In this section, we extend the  GBDF2-IMEX scheme to nonuniform time steps for  Eq.~\cref{eq1.1} and give some related analysis.
Let $\beta=\bar{\beta}+\tilde{\beta} $,
where $\bar{\beta}$ is the integer part of $\beta$  and $\tilde{\beta}=\beta - \bar{\beta}$  is the fractional part of $\beta$, satisfying $0<\tilde{\beta}<1$.
 We reused the notations of time steps in \autoref{sec2} and denote $\Delta t_k = t^{k+1} - t^k ~(k=0,1,...,M-1)$.
We first construct a scheme for Eq.~\cref{eq1.1}  based on the Taylor expansion at time $t^{n+\beta}=t^{n+\bar{\beta}}+\tilde{\beta} \Delta t_{n+\bar{\beta}}~(\beta >1)$.
 According to the Taylor expansion at time $t^{n+\beta}~(\beta>1)$, we get
\begin{equation*}
	u\left(t^{n+1-i}\right)=u(t^{n+\beta})+(t^{n+1-i}-t^{n+\beta})\partial_tu(t^{n+\beta})+\dfrac{1}{2}(t^{n+1-i}-t^{n+\beta})^2\partial_{tt}u(t^{n+\beta})+\mathcal{O}\left(\tau_1^{3}\right), ~\mathrm{for}~  0 \le i \leq 2 ,
\end{equation*}
where $\tau_1=\max \{t^{n+1}-t^{n+\beta},t^{n}-t^{n+\beta},t^{n-1}-t^{n+\beta}\}$.

Consequently, we can derive difference formulae for approaching $\partial_t u(t^{n+\beta})$ and $u(t^{n+\beta})$, respectively
\begin{align*}
	\sum_{q = 0}^{2}a_{n,q}u(t^{n + q-1}) &= \partial_{t}u(t^{n + \beta}) + \mathcal{O}(\tau_1^{2}), \\
	b_{n,1}u(t^{n+1})+b_{n,0}u(t^n) &= u(t^{n + \beta}) + \mathcal{O}(\tau_2^{2}).
\end{align*}
Here, 
\begin{equation*}
	\begin{aligned}
		a_{n,0} &= \frac{2t_{n+\beta} - t_n - t_{n+1}}{(t_{n-1} - t_n)(t_{n-1} - t_{n+1})}, \quad
		a_{n,1} = \frac{2t_{n+\beta} - t_{n-1} - t_{n+1}}{(t_n - t_{n-1})(t_n - t_{n+1})}, \quad
		a_{n,2} = \frac{2t_{n+\beta} - t_{n-1} - t_n}{(t_{n+1} - t_{n-1})(t_{n+1} - t_n)},
	\end{aligned}
\end{equation*}
\begin{equation*}
	\begin{aligned}
		b_{n,0} = \frac{t_{n+1} - t_{n+\beta}}{t_{n+1} - t_n}, \quad
		b_{n,1} = \frac{t_{n+\beta} - t_n}{t_{n+1} - t_n},\quad
		\tau_2=\max\{{t^{n+1}-t^{n+\beta}, t^{n}-t^{n+\beta}}\}.
	\end{aligned}
\end{equation*}

To handle the nonlinear term in Eq.~\cref{eq1.1}, we also need the following explicit difference formula for approximation $u(t^{n+\beta})$
\begin{equation*}
	c_{n,1}u(t^{n})+c_{n,0}u(t^{n-1}) = u(t^{n + \beta}) + \mathcal{O}(\tau_3^{2}),
\end{equation*}
where 
\begin{equation*}
		c_{n,1} = \frac{t_{n+\beta} - t_{n-1}}{t_n - t_{n-1}},  \quad c_ {n,0} = \frac{t_n - t_{n+\beta}}{t_n - t_{n-1}},\quad 	\tau_3=\max\{{t^{n}-t^{n+\beta}, t^{n-1}-t^{n+\beta}}\}.
\end{equation*}

Obviously, if $\Delta t_k=\Delta t~(k=0,1,\ldots,M-1)$, the coefficients $a_{n,q}~(q=0,1,2),~b_{n,q}~(q=0,1)$ and $c_{n,q}~(q=0,1)$ become the uniform one.

With these at hand, a semi-discrete scheme with nonuniform time steps  for Eq.~\cref{eq1.1} is
\begin{equation*}
	\sum_{q = 0}^{2}a_{n,q}u^{n -1 +q}=\mathcal{L}_h(b_{n,1}u^{n+1}+b_{n,0}u^n)+\mathcal{G}( c_{n,1}u^{n}+c_{n,0}u^{n-1}).
\end{equation*}

Combining the spatial discrete method and reusing the notations in \Cref{sec2}, we have the  GBDF2-IMEX scheme with the nonuniform time steps in the form of matrix-vector multiplication  as follows
\begin{equation}
	(a_{n,2}I_{xy} - \kappa b_{n,1} L) \bm{u}^{n+1} = (-a_{n,1}I_{xy} + \kappa b_{n,0} L) \bm{u}^{n} + \mathcal{G} (c_{n,0} \bm{u}^{n-1} + c_{n,1} \bm{u}^{n}) - a_{n,0} \bm{u}^{n-1} +\kappa \bm{\xi}^n.
	\label{eq3.1}
\end{equation}
We conjecture that the convergence order of scheme \cref{eq3.1} is \(\mathcal{O}(\tau^2+h_x^2+h_y^2)\), where $\tau = \max \limits_{0 \leq k \leq M-1} \{\Delta t_k\}$.
Note that the  structure of the coefficient matrix in  \cref{eq3.1} is the same as that of \cref{eq2.8}. 
Therefore, Algorithm \ref{alg2.1} can still be used to solve \cref{eq3.1} efficiently.

 Recent studies demonstrate that tanh-based nonuniform grids significantly enhance error balance for time-fractional equations \cite{zhang2021error}. 
Inspired by this, we extend such grids to Eq.~\cref{eq1.1}.
The time interval $[0, T]$ is split into a series of nonuniform subintervals $[t_n, t_{n+1}]$, for a given positive integer $M$, with
\begin{equation*}
	t^n := T - T\varphi(\beta(n), y(n)),   \qquad  \qquad \mathrm{for }~ n = 0, 1, \ldots, M,
\end{equation*}
where
\begin{equation*}
	\varphi(\mu_1, \mu_2) := \frac{\tanh(\mu_1\mu_2)}{\tanh(\mu_1)},  \qquad \qquad \mathrm{for}~ \forall \mu_1, \mu_2 \in \mathbb{R}, \quad 0 < \mu_1, 0 \leq \mu_2 \leq 1,
\end{equation*}
\begin{equation*}
	\beta(n) := \frac{\gamma}{2}N \frac{\ln(N + 2) - \ln(n + 1)}{N - n + 1}, \qquad \mathrm{for }~ \forall \gamma \in \mathbb{R}^+, \ \forall n \in \left\{0,1,\ldots,M\right\},
\end{equation*}
\begin{equation*}
	y(n) := \frac{N - n}{N} ,  \qquad  \qquad \mathrm{for}~ \forall n \in \left\{0,1,\ldots,M\right\}.
\end{equation*}

The $\tanh$ function clusters grid points near interval boundaries (e.g., $t\to 0 ~\mathrm{or}~ t\to T$), capturing rapid transients or boundary layers. The parameter $\gamma$ controls the density of the grid.
\begin{remark}
	\label{remark1}
	For nonuniform time steps, the coefficients $a_{n,q}, b_{n,q}, c_{n,q}$ can be derived using Lagrange interpolation polynomials. Let $u(t)$ be the Lagrange interpolating polynomial over the nodes $t_{n-1}, t_n$ and $t_{n+1}$:
	\[
	u(t) = u(t^{n-1})L_0(t) + u(t^{n})L_1(t) + u(t^{n+1})L_2(t),
	\]
	where the Lagrange basis functions are:
	\begin{align*}
		L_{0}(t) &= \frac{(t - t_{n+1})(t - t_n)}{(t_{n-1} - t_{n+1})(t_{n-1} - t_n)}, \\
		L_{1}(t) &= \frac{(t - t_{n-1})(t - t_{n+1})}{(t_{n} - t_{n-1})(t_{n} - t_{n+1})}, \\
		L_{2}(t) &= \frac{(t - t_{n})(t - t_{n-1})}{(t_{n+1} - t_{n})(t_{n+1} - t_{n-1})}.
	\end{align*}
	
	The derivative of $u(t)$ at $t_{n+\beta}$ is:
	\begin{equation*}
			u'(t_{n+\beta}) = a_{n,0}u(t^{n-1})  +a_{n,1} u(t^{n}) +a_{n,2} u(t^{n+1}).
	\end{equation*}
	
		For $u(t_{n+\beta})$ in the  linear and nonlinear terms, this yields respectively
			\[u(t_{n+\beta})=b_{n,0} u(t^{n}) + b_{n,1} u(t^{n+1})\qquad \mathrm{and}\qquad u(t_{n+\beta})=c_{n,0} u(t^{n-1})^{n-1} + c_{n,1} u(t^{n}).
		\]
		
		In our future work, we will research the stability and convergence of the nonuniform scheme \cref{eq3.1}.
\end{remark}

\section{Numerical experiments}
\label{sec4}
In this section, several numerical examples are presented to validate the convergence orders of our schemes in both time and space.
It is easy to know that Eqs.~\cref{eq2.8} and \cref{eq3.1} are not self-starting, as they require solutions at two previous time levels. 
This necessitates the use of a separate, high-accuracy method to compute the solution at time $t_1$. 
To minimize the impact of this initialization on the overall convergence, we employ the high-order explicit Dormand-Prince 5(4) Runge-Kutta method \cite{seen2014gpu}, which provides sufficient accuracy to ensure that the initialization error is negligible compared to the errors introduced by the main GBDF2-IMEX scheme.
	
To validate the effectiveness of the proposed numerical method, we compute the errors at $t=T$ in $L_\infty$-norm and $L_2$-norm, respectively. The definitions of these two types of errors are given as follows:
	\begin{equation*}
		Err_\infty (M,N_x,N_y)= \max_{1 \leq i \leq N_x,1 \leq j \leq N_y} |u(x_i, y_j, T) - u_{i,j}^T|,
	\end{equation*}
\begin{equation*}
    Err_2 (M,N_x,N_y)= \left(h_x h_y \sum_{i=1}^{N_x} \sum_{j=1}^{N_y} |u(x_i, y_j, T) - u_{i,j}^M|^2 \right)^{1/2},
\end{equation*}

	The convergence orders are  computed using the following formulae:
\begin{equation*}
	Order_{t,\infty} = \frac{\log \left( \frac{Err_\infty (M, N_x,N_y)}{Err_\infty (2M, N_x,N_y)} \right)}{\log 2}, \quad 
	Order_{t,2} = \frac{\log \left( \frac{Err_2 (M, N_x,N_y)}{Err_2 (M, 2N_x,2N_y)} \right)}{\log 2}.
\end{equation*}
\begin{equation*}
	Order_{h,\infty}  = \frac{\log \left( \frac{Err_\infty  (M, N_x,N_y)}{Err_\infty  (M, 2N_x,2N_y)} \right)}{\log 2}, \quad 
	Order_{h,2}  = \frac{\log \left( \frac{Err_2 (M, N_x,N_y)}{Err_2 (M, 2N_x,2N_y)} \right)}{\log 2}.
\end{equation*}

In all  figures, the dashed line labeled $\mathrm{Slope}=2$ represents the theoretically expected second-order convergence rate.
 Moreover, for convenience, we set $N_x=N_y=N $ in all experiments.
  All experiments are conducted using MATLAB2024b on a Intel(R) Core(TM) i9-14900KF 3.20GHz of processor and 128GB of RAM.

{\noindent \textbf{Example 1.}}
We consider the following two-dimensional Fisher-KPP equation with the zero boundary conditions:
\begin{equation}
	\frac{\partial u(x,y,t)}{\partial t} = \left(\frac{\partial^2 }{\partial x^2} + \frac{\partial^2 }{\partial y^2}\right)u(x,y,t) + u(x,y,t)(1 - u(x,y,t)) + f(x,y,t), \quad (x, y, t) \in [0,\pi]^2 \times [0, 1].
	 \label{eq4.1}
\end{equation}
The exact solution  is  assumed as \( u(x,y,t) = \sin(x) \sin(y) \sin(t) \). 
Then, the source term is 
\begin{equation*}
	f(x,y,t) = \sin(x) \sin(y) (\cos(t)+\sin(t)+\sin(x)\sin(y)\sin^2(t)) .
\end{equation*}

\begin{figure}[htbp]
	\setlength{\tabcolsep}{0.2pt}
	\centering
	\begin{tabular}{m{0.4cm}<{\centering} m{5.18cm}<{\centering} m{5.18cm}<{\centering} m{5.18cm}<{\centering}}
		& $\beta=\sqrt2$ & $\beta=2$ & $\beta=\pi$ \\
		\rotatebox{90}{{\tt $Err_\infty$}} &
		\includegraphics[width=2.0in,height=1.9in]{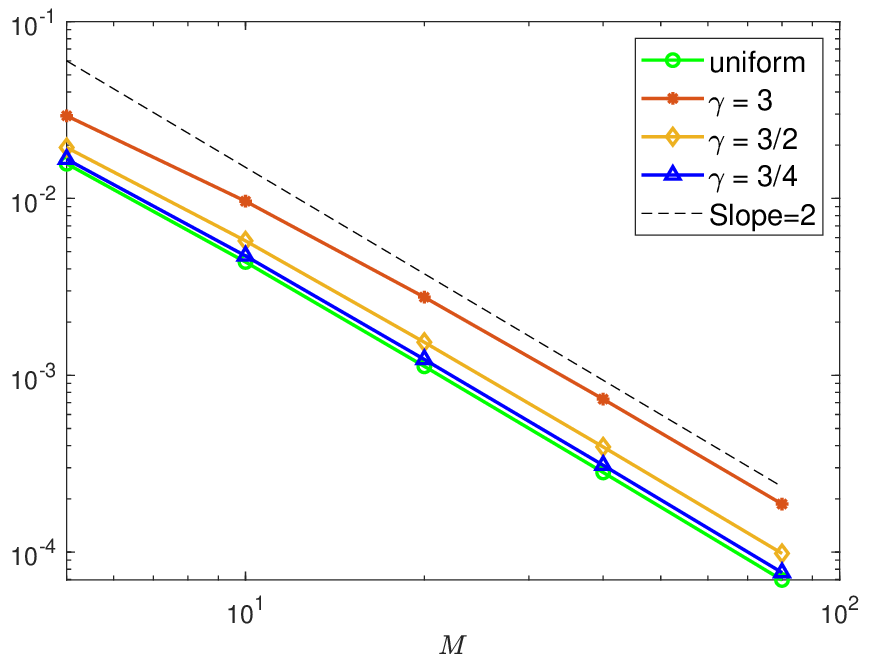} &
		\includegraphics[width=2.0in,height=1.9in]{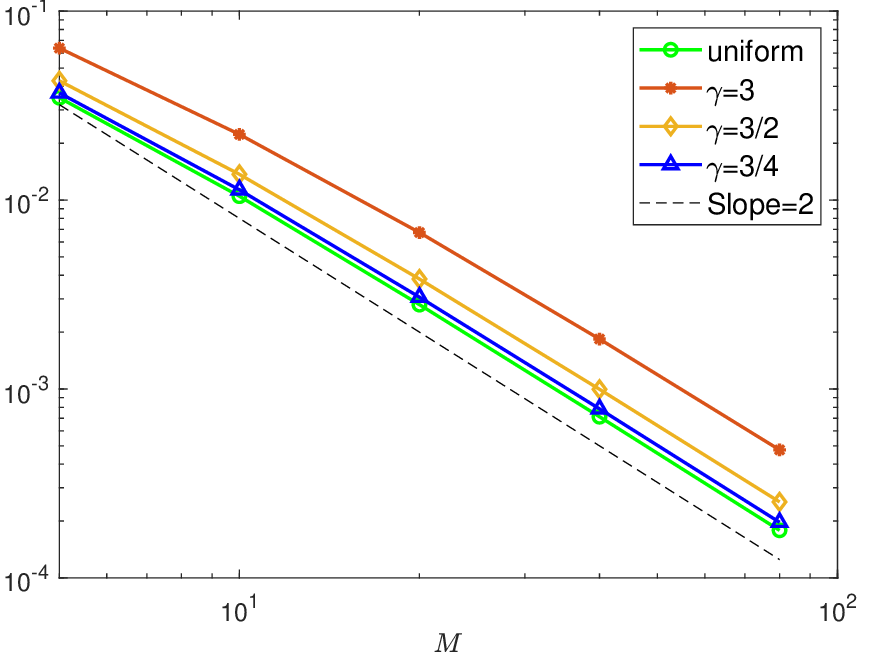} &
		\includegraphics[width=2.0in,height=1.9in]{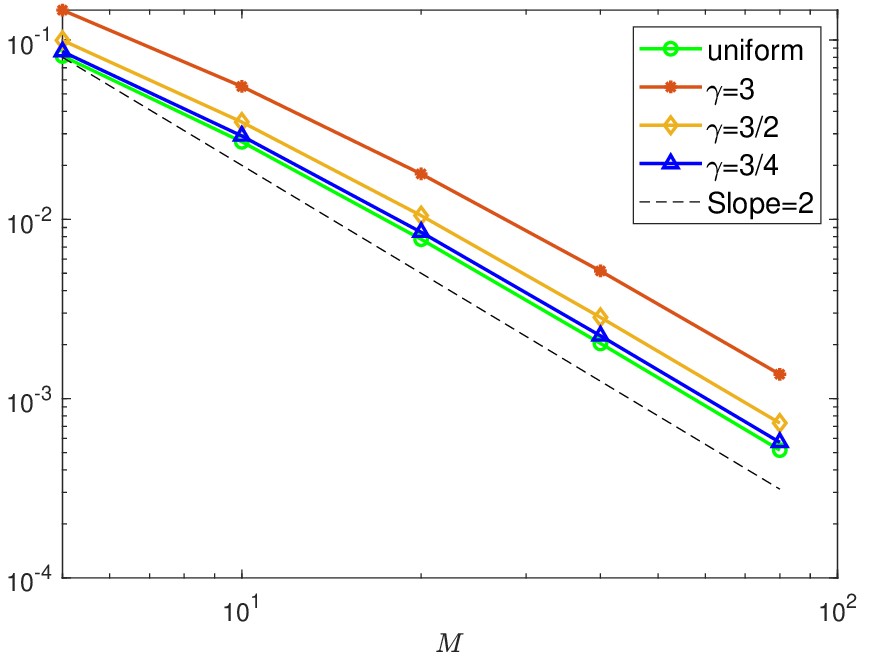} \\
		\rotatebox{90}{{\tt $Err_2$}} &
		\includegraphics[width=2.0in,height=1.9in]{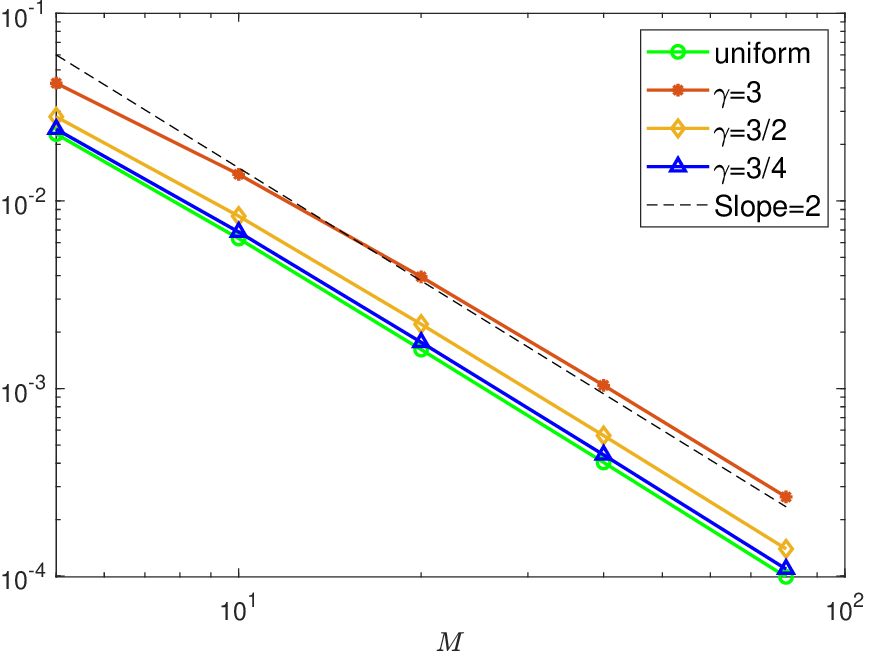} &
		\includegraphics[width=2.0in,height=1.9in]{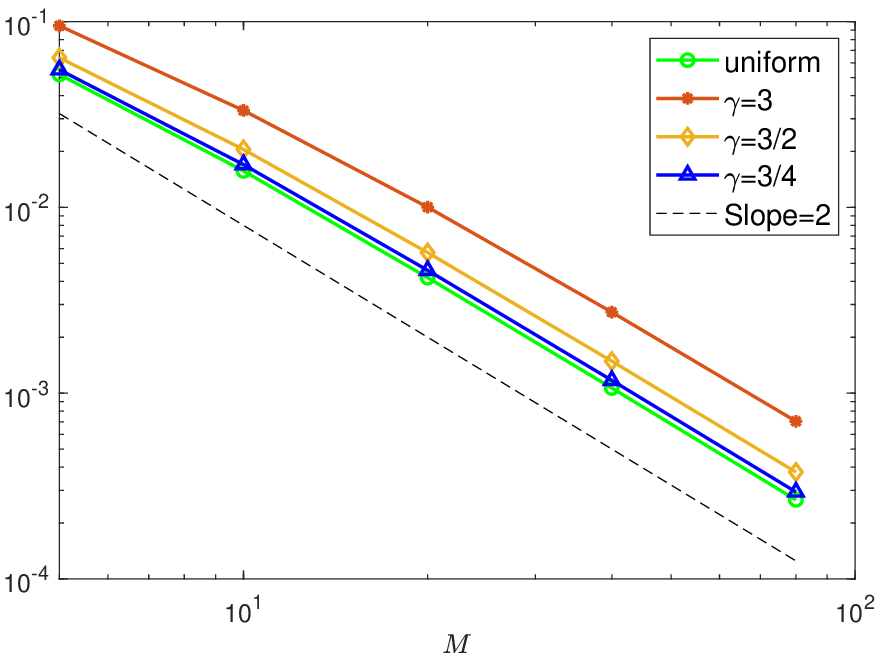} &
		\includegraphics[width=2.0in,height=1.9in]{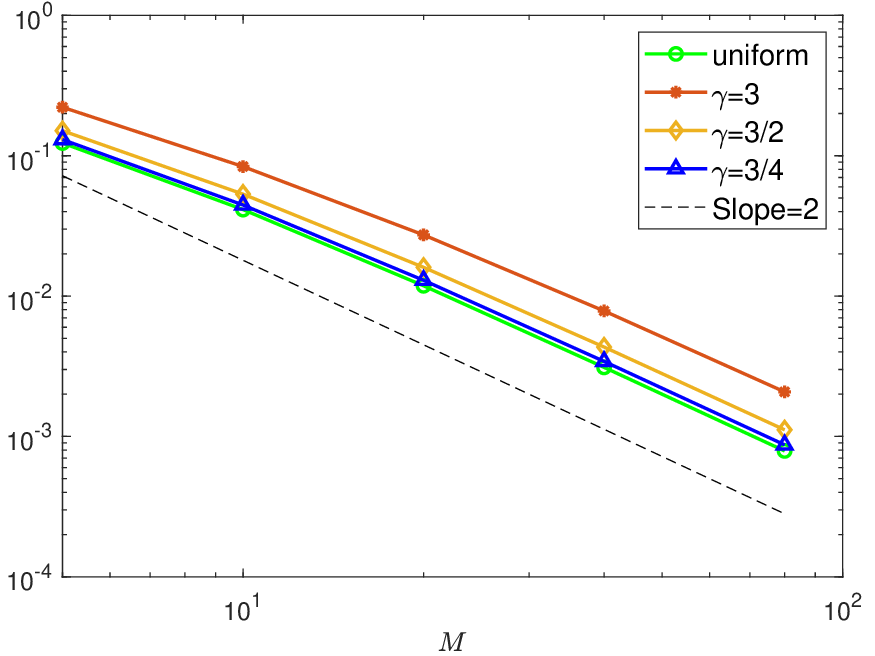} \\
	\end{tabular}
	\caption{Observed temporal errors for Example 1 for different values of $\beta$ and $M$, where $N=500$.}
	\label{fig1}
\end{figure}
\cref{fig1} demonstrates the error convergence behavior in the temporal direction under different parameters $\beta$ for Eq.~\cref{eq4.1}. 
 Both $Err_\infty$ and $Err_2$  curves follow the reference line with $\mathrm{Slope}=2$, consistent with the theoretical second-order accuracy in the temporal direction.
 The uniform mesh and the nonuniform mesh with $\gamma=3/4$ show comparable error performance.
 
\begin{figure}[htbp]
	\setlength{\tabcolsep}{0.2pt}
	\centering
	\begin{tabular}{m{0.4cm}<{\centering} m{5.18cm}<{\centering} m{5.18cm}<{\centering} m{5.18cm}<{\centering}}
		& $\beta=\sqrt2$ & $\beta=2$ & $\beta=\pi$ \\
		\rotatebox{90}{{\tt $Err_\infty$}} &
		\includegraphics[width=2.0in,height=1.9in]{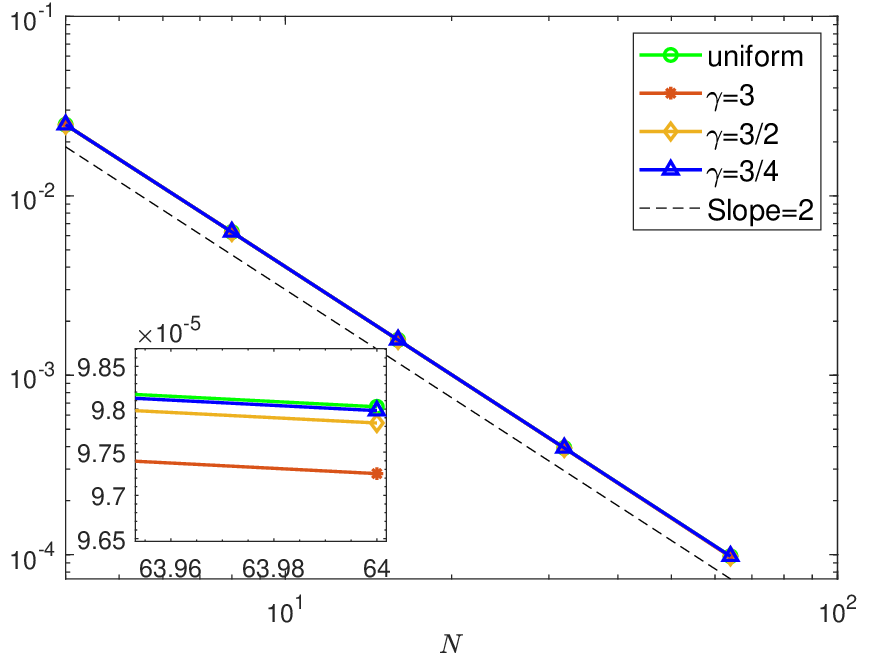} &
		\includegraphics[width=2.0in,height=1.9in]{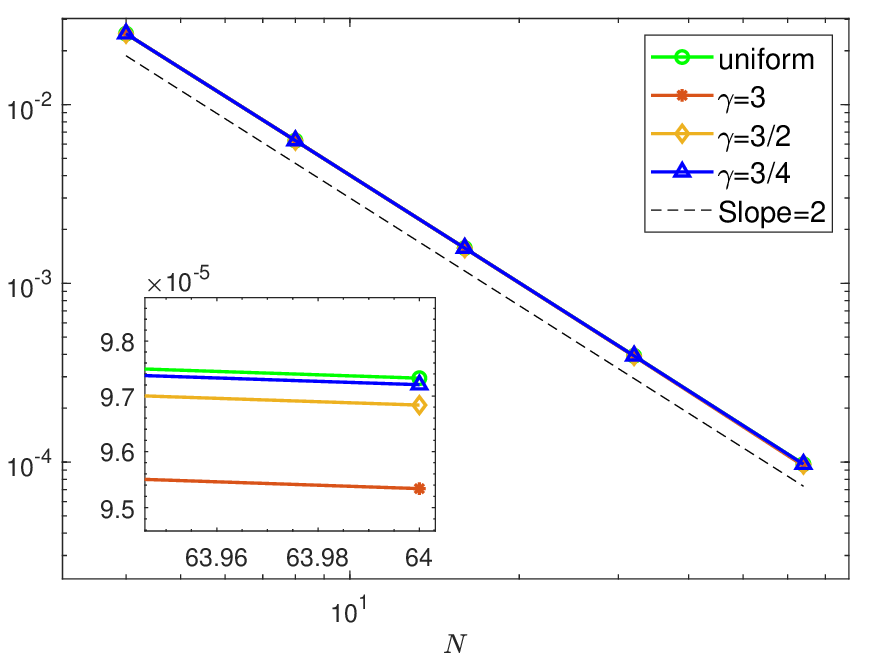} &
		\includegraphics[width=2.0in,height=1.9in]{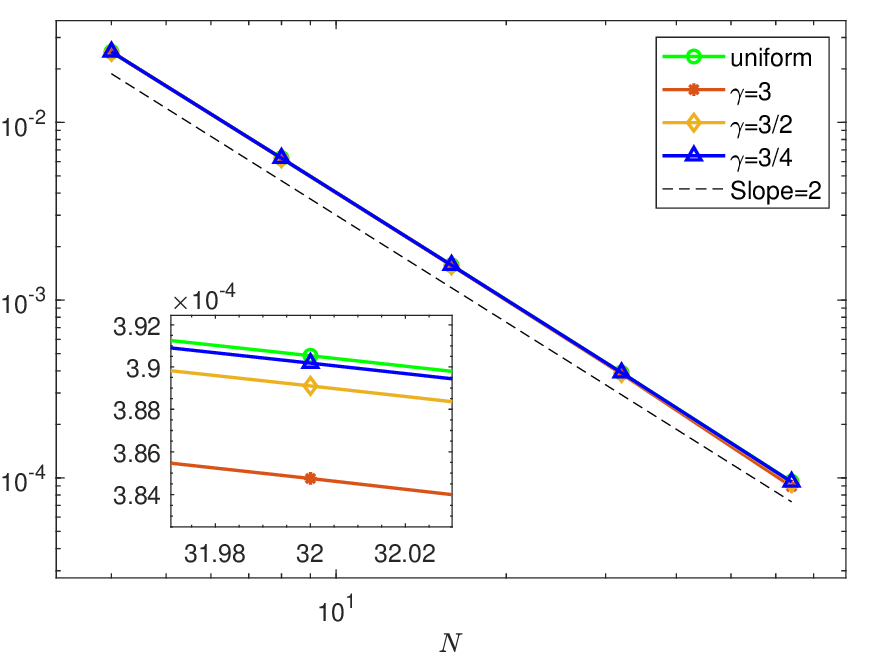} \\
		\rotatebox{90}{{\tt $Err_2$}} &
		\includegraphics[width=2.0in,height=1.9in]{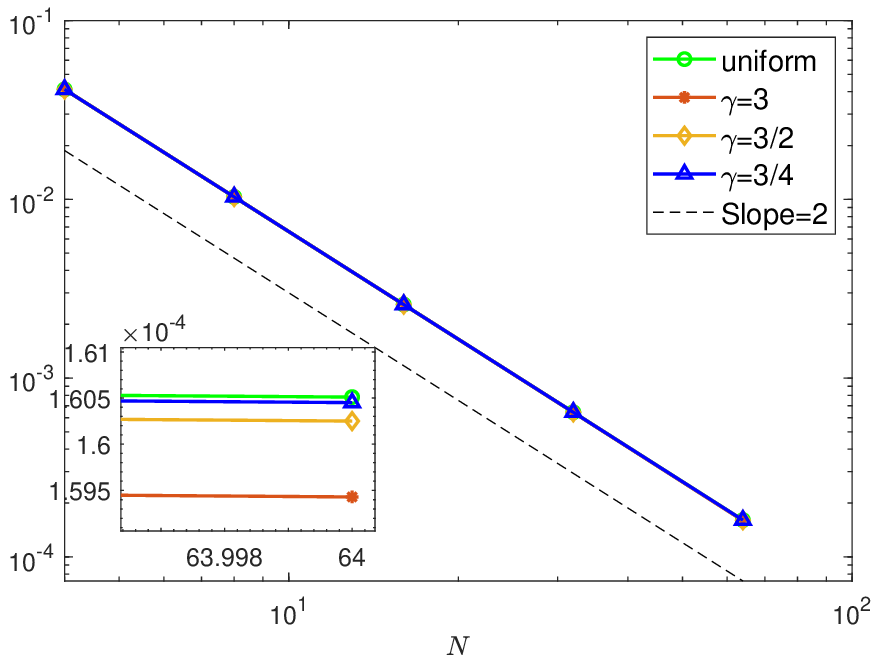} &
		\includegraphics[width=2.0in,height=1.9in]{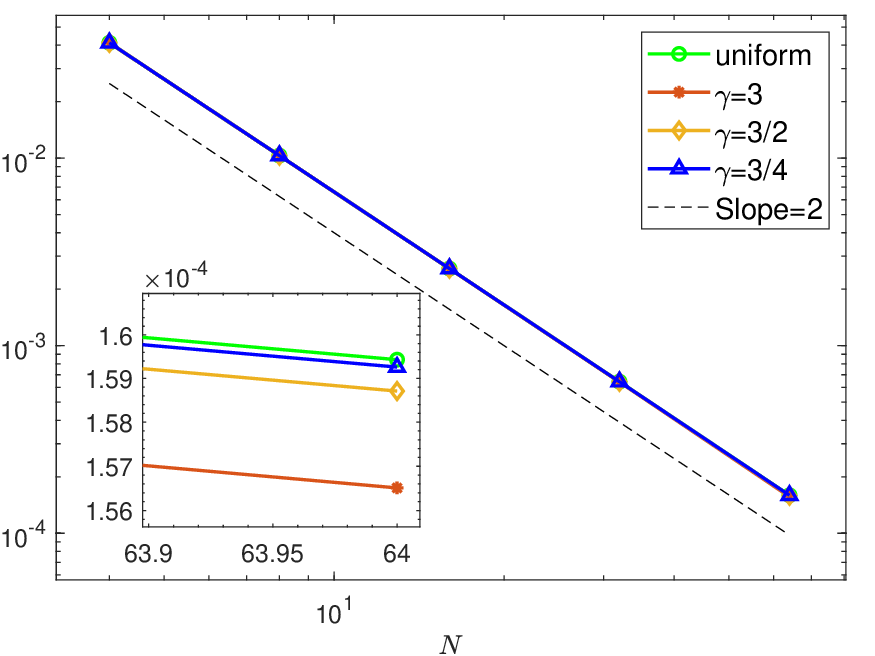} &
		\includegraphics[width=2.0in,height=1.9in]{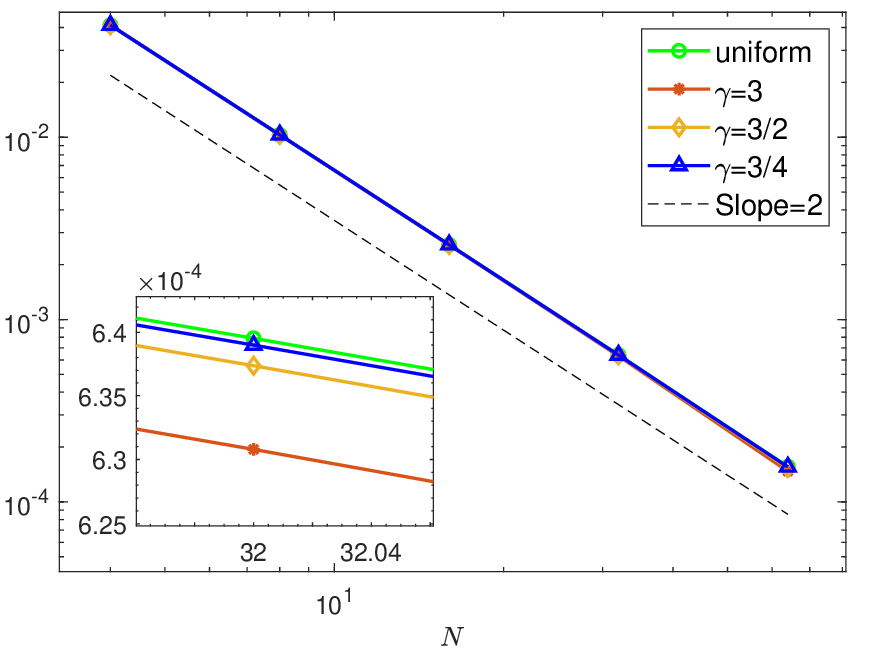} \\
	\end{tabular}
	\caption{Observed space errors for Example 1 for different values of $\beta$ and $N$, where $M=1000$.}
	\label{fig2}
\end{figure}
\cref{fig2} demonstrates the error convergence behavior in the spatial direction for Eq.~\cref{eq4.1}.
Both the $Err_\infty$ and $Err_2$ error curves (including those with different $\beta$ and $\gamma$ parameters) follow the reference line of slope 2, which confirms that the spatial discretization achieves second-order convergence.

{\hypertarget{example2}\noindent \textbf{Example 2.}(\cite{orucc2020efficient})}
In this example, we consider Eq.~\eqref{eq1.1} with the nonlinear term $f(u) = u(1 - u)$ (i.e., the first form in Eq.~\eqref{eq1.2} with $K=p=1$). The other parameters are $\kappa = 1$, $\Omega = [-1, 1]^2$, and $T=1$. The exact solution is given by
\begin{equation*}
	u(x, y, t) = \left[ 1 + \exp\left( \frac{(x - y)/ \sqrt{2} - (5 / \sqrt{6}) t}{\sqrt{6}} \right) \right]^{-2}.
\end{equation*}

\begin{figure}[htbp]
	\setlength{\tabcolsep}{0.2pt}
	\centering
	\begin{tabular}{m{0.4cm}<{\centering} m{5.18cm}<{\centering} m{5.18cm}<{\centering} m{5.18cm}<{\centering}}
		& $\beta=\sqrt2$ & $\beta=2$ & $\beta=\pi$ \\
		\rotatebox{90}{{\tt $Err_\infty$}} &
		\includegraphics[width=2.0in,height=1.9in]{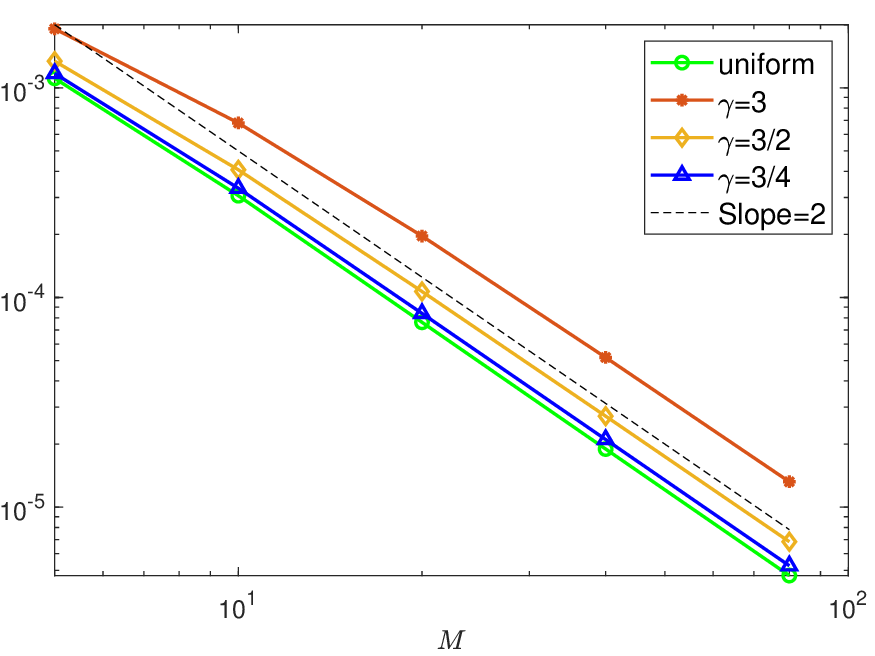} &
		\includegraphics[width=2.0in,height=1.9in]{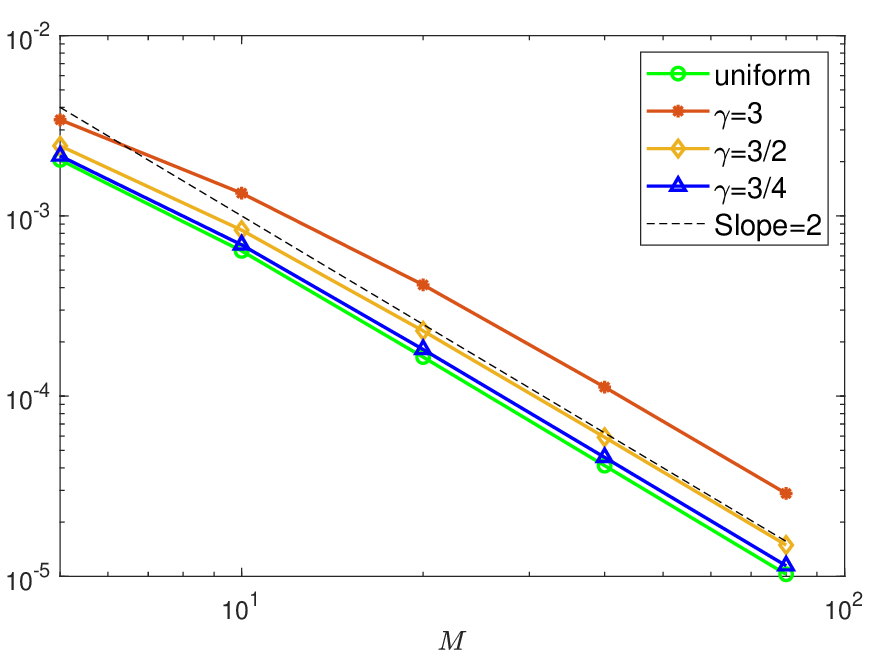} &
		\includegraphics[width=2.0in,height=1.9in]{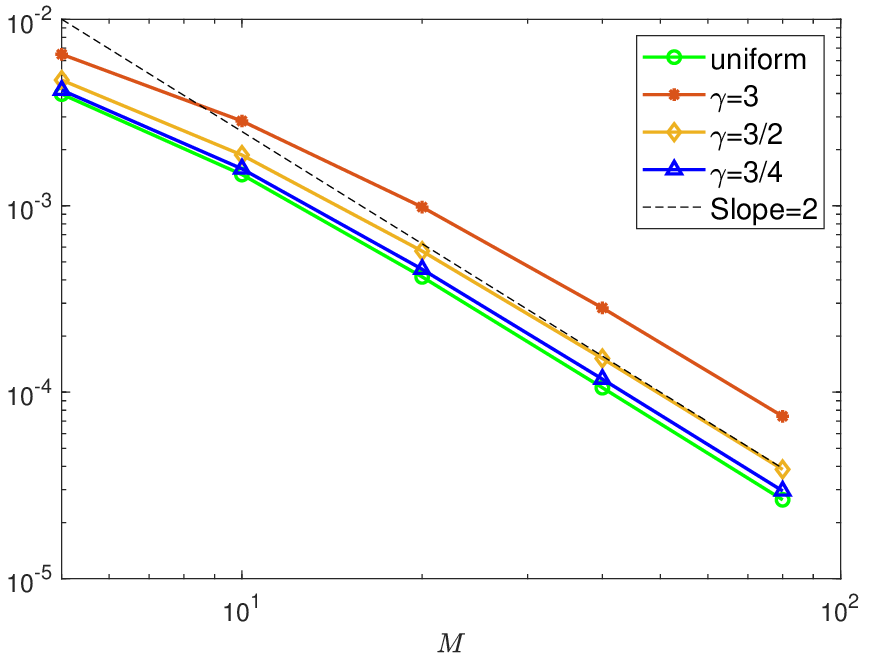} \\
		\rotatebox{90}{{\tt $Err_2$}} &
		\includegraphics[width=2.0in,height=1.9in]{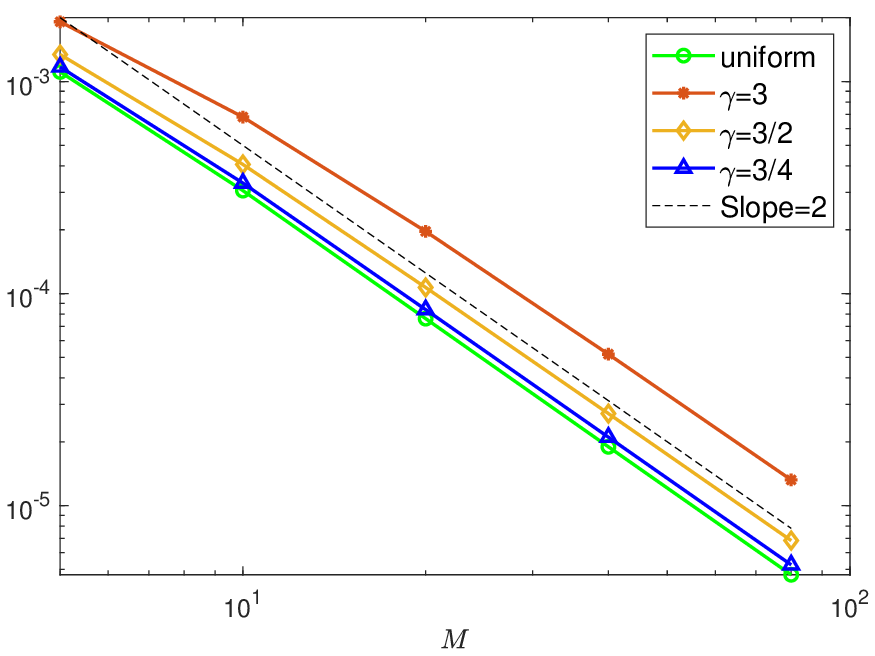} &
		\includegraphics[width=2.0in,height=1.9in]{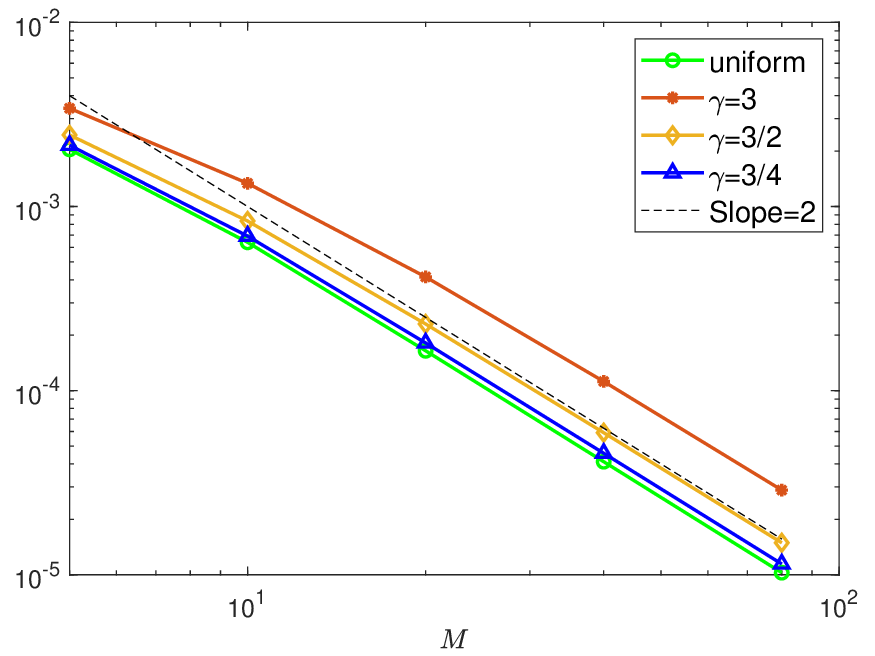} &
		\includegraphics[width=2.0in,height=1.9in]{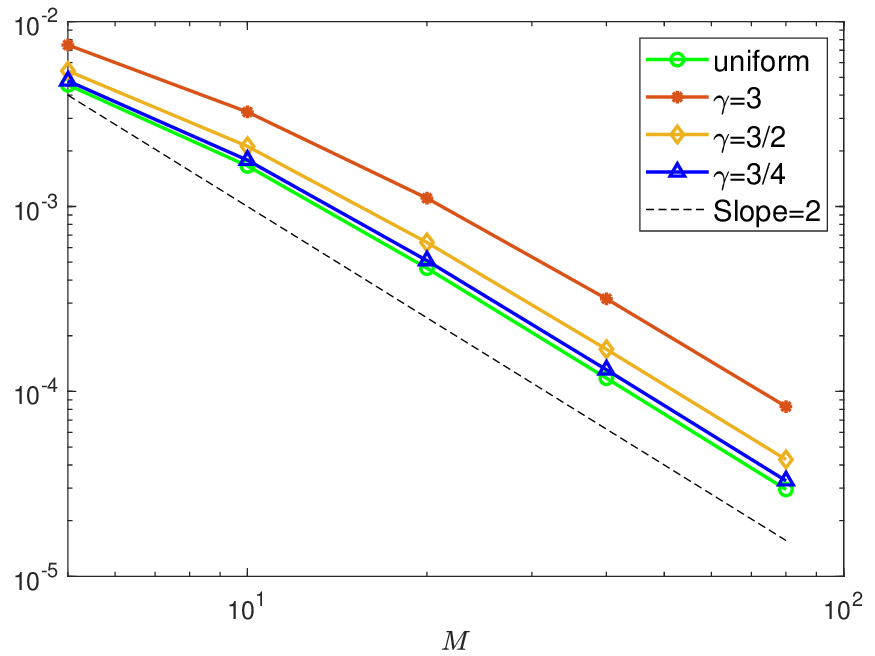} \\
	\end{tabular}
	\caption{Observed temporal errors for Example 2 for different values of $\beta$ and $M$, where $N=500$.}
	\label{fig3}
\end{figure}
 \cref{fig3}  validates the second-order temporal convergence in \hyperlink{example2}{Example 2}.
  Both $Err_\infty$ and $Err_2$ error curves  follow the reference line with $\mathrm{Slope}=2$,  confirming the second-order accuracy in the temporal direction. 
 
\begin{figure}[htbp]
	\setlength{\tabcolsep}{0.2pt}
	\centering
	\begin{tabular}{m{0.4cm}<{\centering} m{5.18cm}<{\centering} m{5.18cm}<{\centering} m{5.18cm}<{\centering}}
		& $\beta=\sqrt2$ & $\beta=2$ & $\beta=\pi$ \\
		\rotatebox{90}{{\tt $Err_\infty$}} &
		\includegraphics[width=2.0in,height=1.9in]{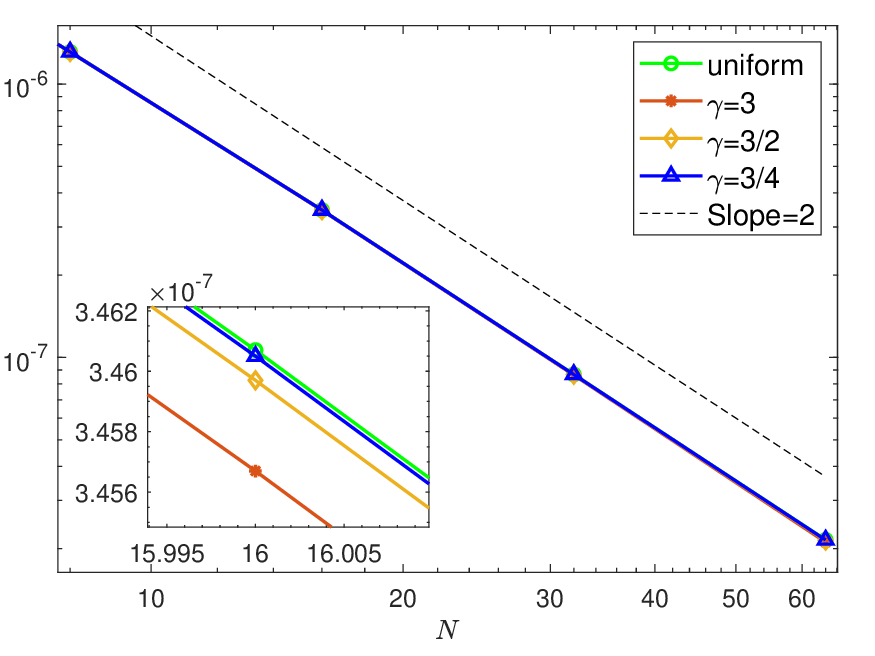} &
		\includegraphics[width=2.0in,height=1.9in]{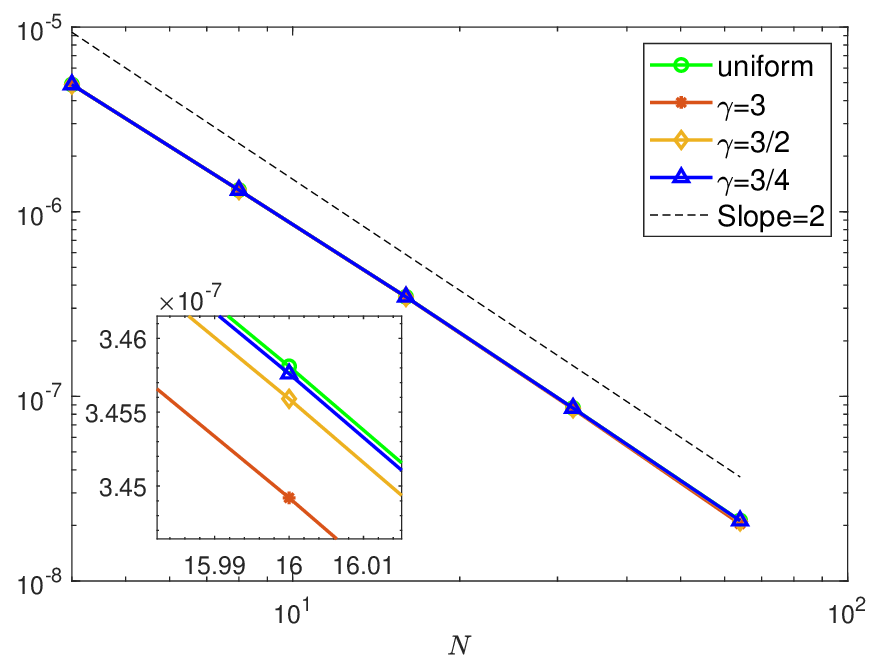} &
		\includegraphics[width=2.0in,height=1.9in]{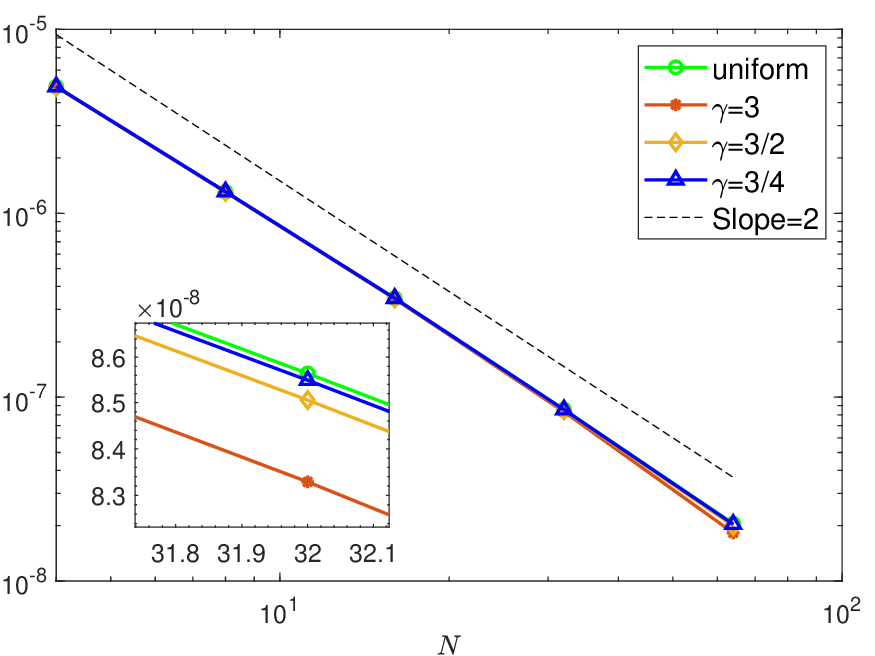} \\
		\rotatebox{90}{{\tt $Err_2$}} &
		\includegraphics[width=2.0in,height=1.9in]{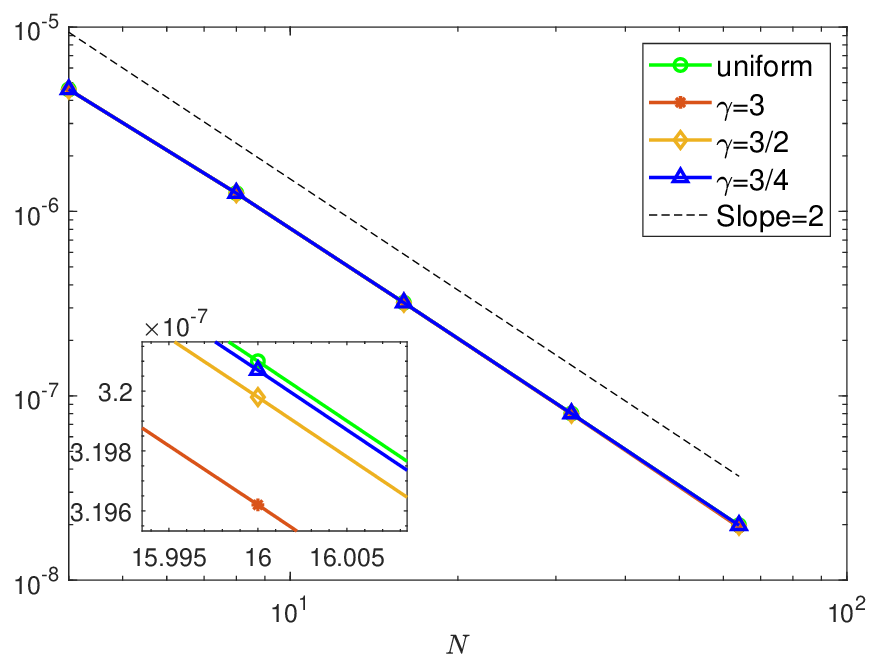} &
		\includegraphics[width=2.0in,height=1.9in]{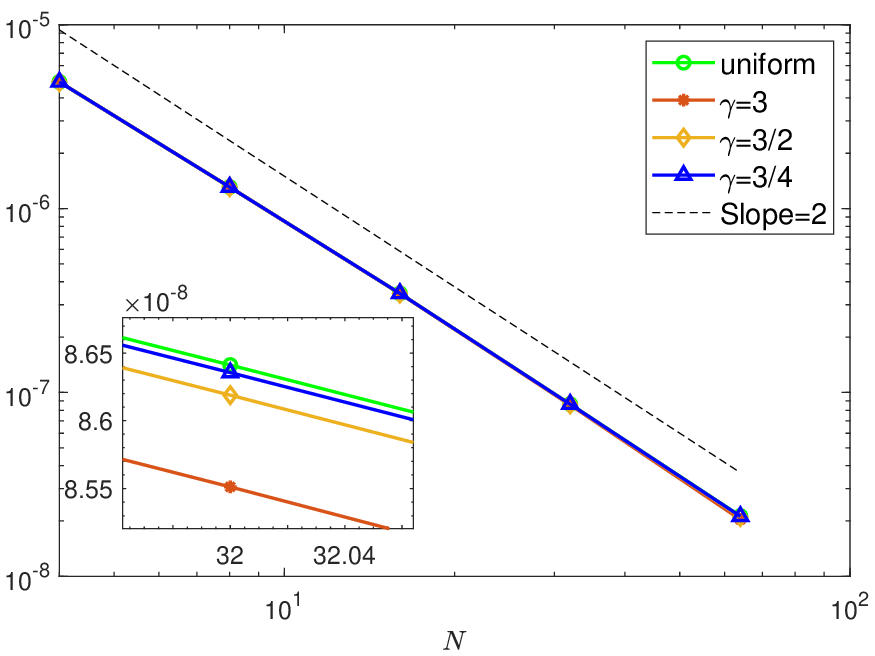} &
		\includegraphics[width=2.0in,height=1.9in]{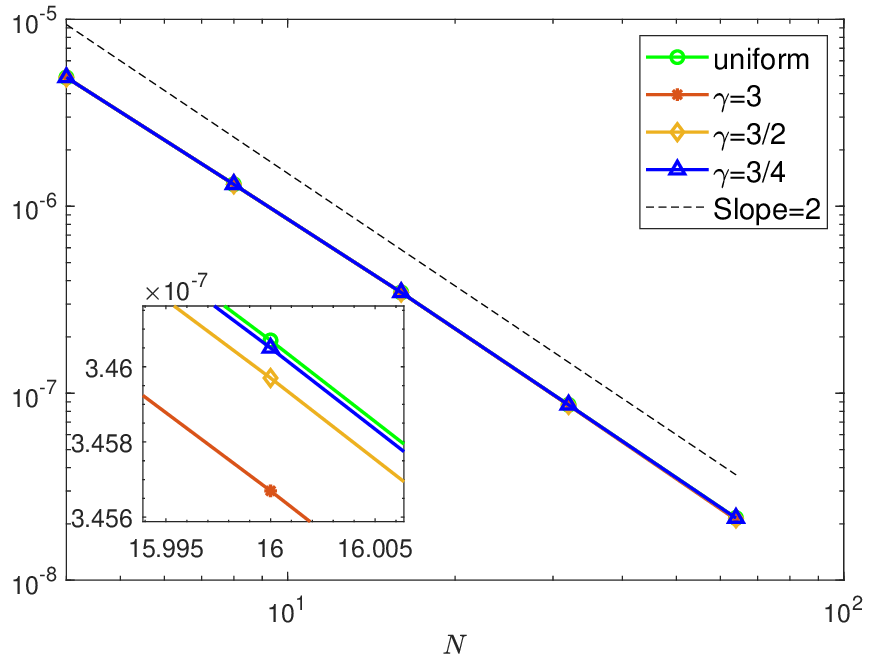} \\
	\end{tabular}
	\caption{Observed space errors for Example 2 for different values of $\beta$ and $N$, where $M=10000$.}
	\label{fig4}
\end{figure}
 \cref{fig4} validates the second-order spatial convergence in \hyperlink{example2}{Example 2}. 
 Both $Err_\infty$ and $Err_2$ error curves  nearly coincide, demonstrating a spatial convergence order of 2.
 
  {\hypertarget{example3}\noindent \textbf{Example 3.}(\cite{tyson2000traveling})
  We consider Eq.~\eqref{eq1.1} with $\kappa = 1$, $f(u) = u(1 - u^2)$ (i.e., the first form in Eq.~\eqref{eq1.2} with $K=1, p=2$).
  The  exact solution  is given
  	\begin{equation*}
  		u(x,y,t)=\left\{\dfrac{1}{2}\tanh[\dfrac{n}{2\sqrt{2n+4}}(x \sin \psi +y \cos \psi)+\dfrac{n(n+4)}{4(n+2)}t+\ell]+\dfrac{1}{2}\right\}^{\dfrac{2}{n}},\quad (x,y,t) \in [-5,5]^2 \times [0,1],
  	\end{equation*}
  where $\psi$ and $\ell$ are arbitrary real constants, and $n$ is a positive integer.
   In this example,  we fix  $\psi =\dfrac{3\pi}{4},~ \ell=0$ and $n=2$.
  	\begin{figure}[htbp]
  		\setlength{\tabcolsep}{0.2pt}
  		\centering
  		\begin{tabular}{m{0.4cm}<{\centering} m{5.18cm}<{\centering} m{5.18cm}<{\centering} m{5.18cm}<{\centering}}
  			& $\beta=\sqrt2$ & $\beta=2$ & $\beta=\pi$ \\
  			\rotatebox{90}{{\tt $Err_\infty$}} &
  			\includegraphics[width=2.0in,height=1.9in]{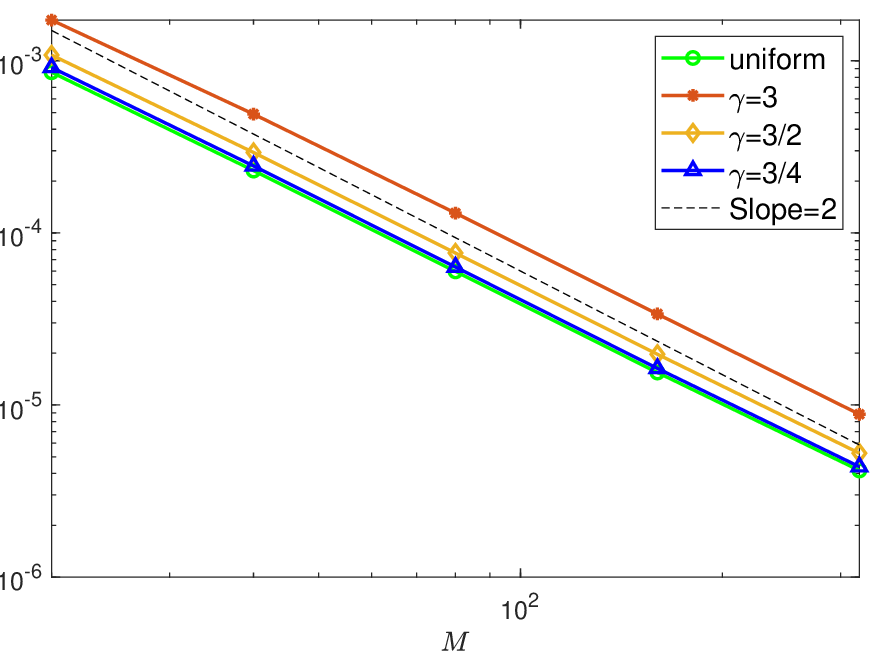} &
  			\includegraphics[width=2.0in,height=1.9in]{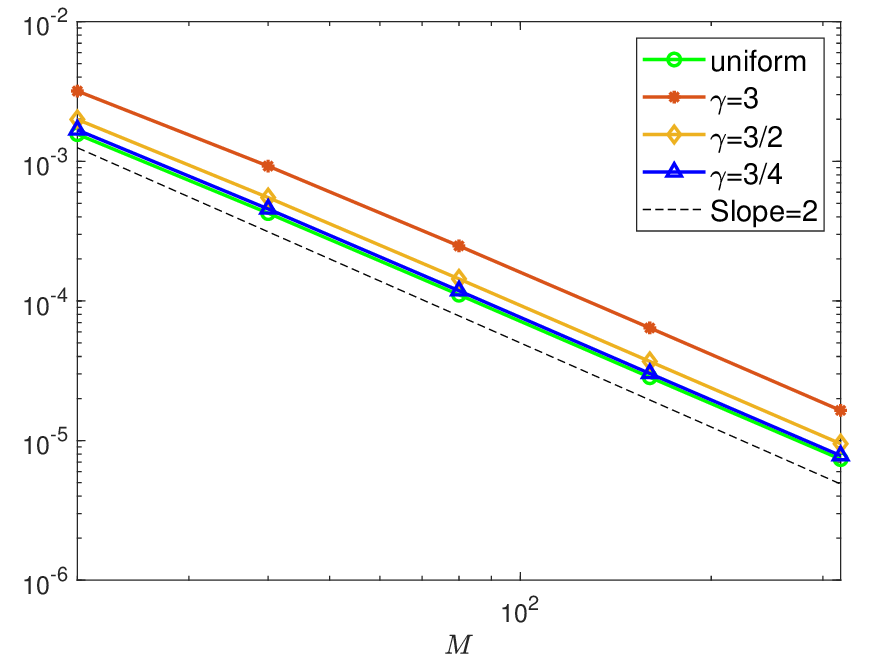} &
  			\includegraphics[width=2.0in,height=1.9in]{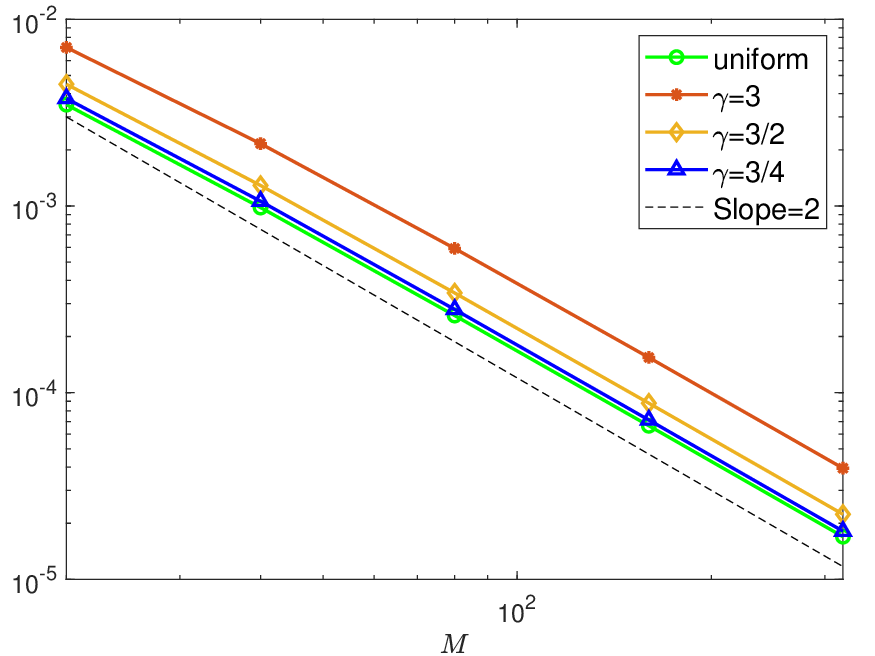} \\
  			\rotatebox{90}{{\tt $Err_2$}} &
  			\includegraphics[width=2.0in,height=1.9in]{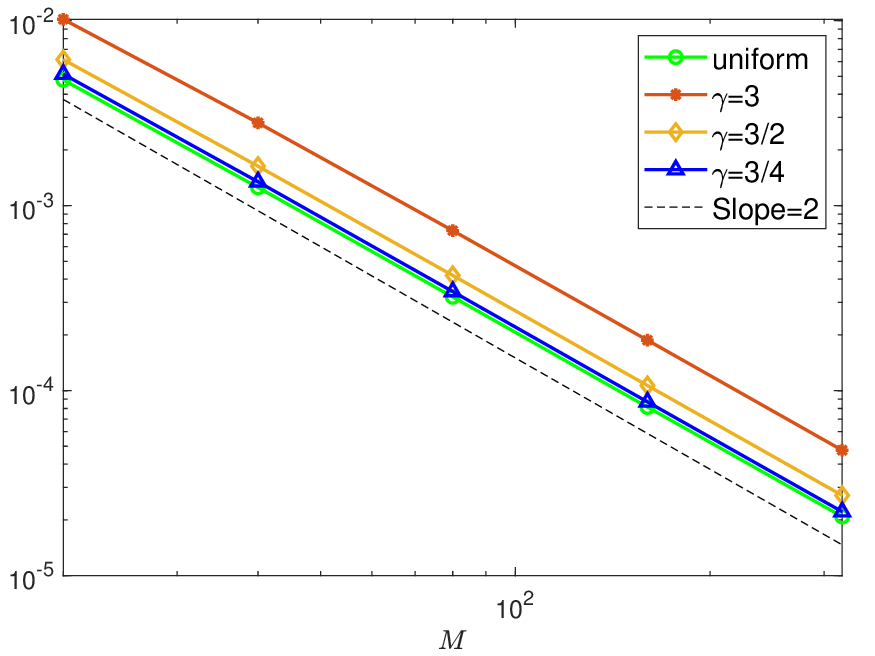} &
  			\includegraphics[width=2.0in,height=1.9in]{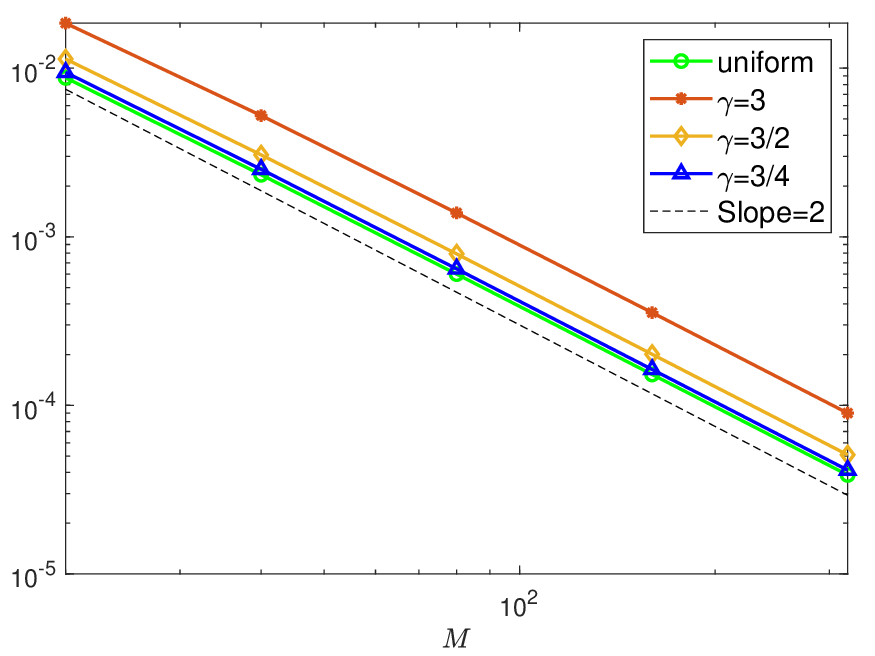} &
  			\includegraphics[width=2.0in,height=1.9in]{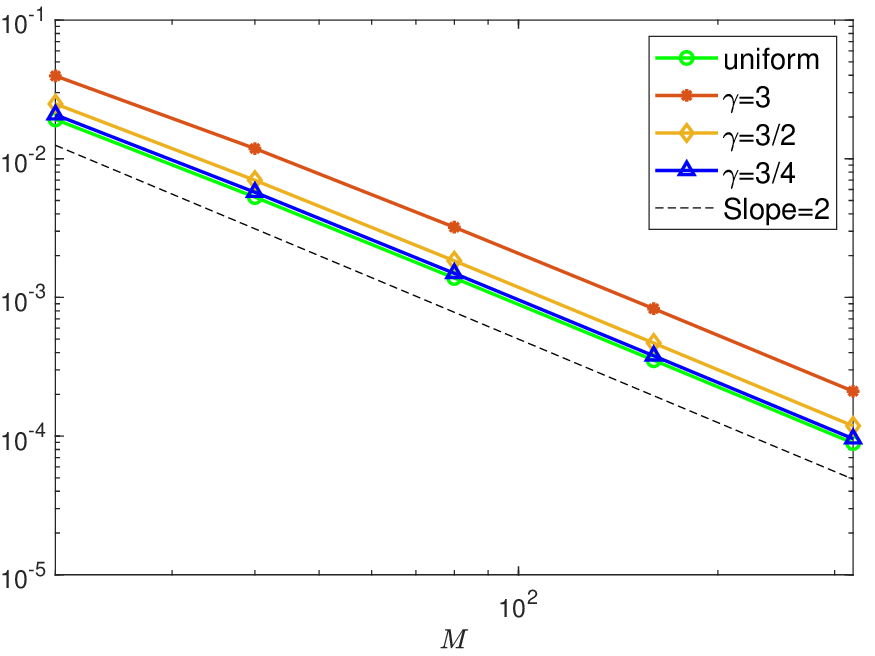} \\
  		\end{tabular}
  		\caption{Observed temporal errors for Example 3 for different values of $\beta$ and $M$, where $N=500$.}
  		\label{fig5}
  	\end{figure}
  
  	\cref{fig5}   validates the second-order temporal convergence of the proposed  GBDF2-IMEX scheme for \hyperlink{example3}{Example 3}.
  	\cref{fig6}  validates the second-order spatial convergence of the proposed numerical scheme for \hyperlink{example3}{Example 3}.
  	\begin{figure}[htbp]
  		\setlength{\tabcolsep}{0.2pt}
  		\centering
  		\begin{tabular}{m{0.4cm}<{\centering} m{5.18cm}<{\centering} m{5.18cm}<{\centering} m{5.18cm}<{\centering}}
  			& $\beta=\sqrt2$ & $\beta=2$ & $\beta=\pi$ \\
  			\rotatebox{90}{{\tt $Err_\infty$}} &
  			\includegraphics[width=2.0in,height=1.9in]{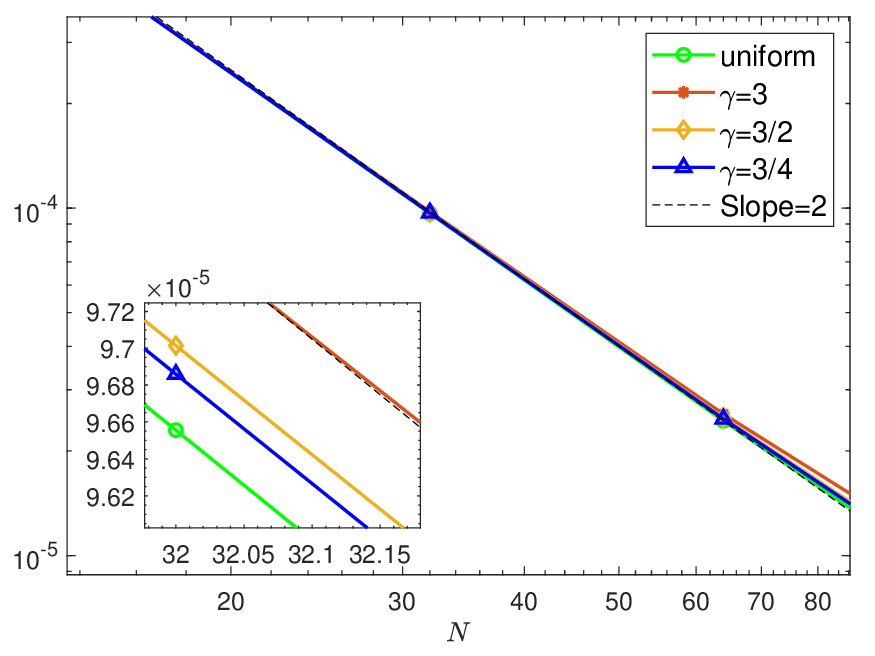} &
  			\includegraphics[width=2.0in,height=1.9in]{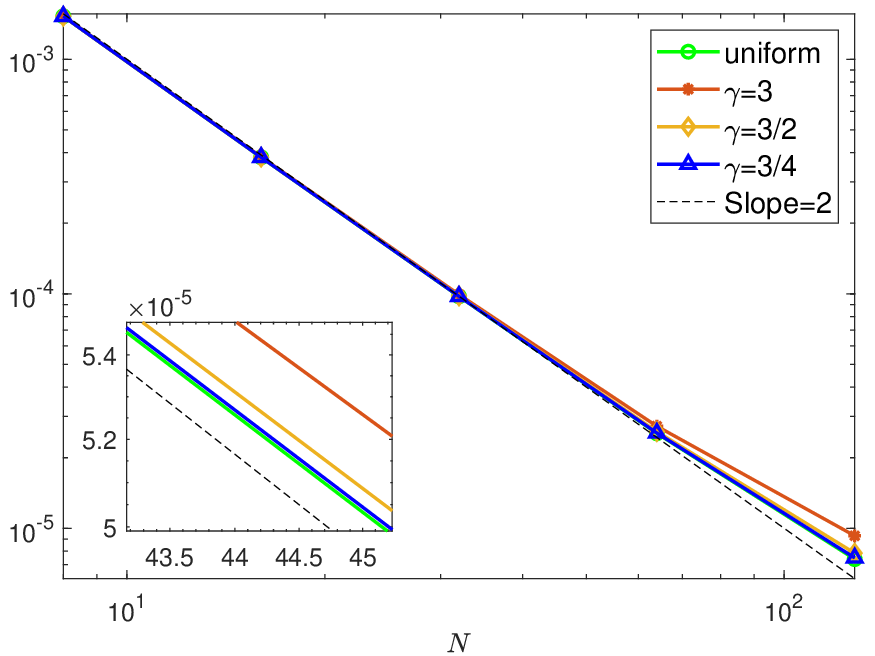} &
  			\includegraphics[width=2.0in,height=1.9in]{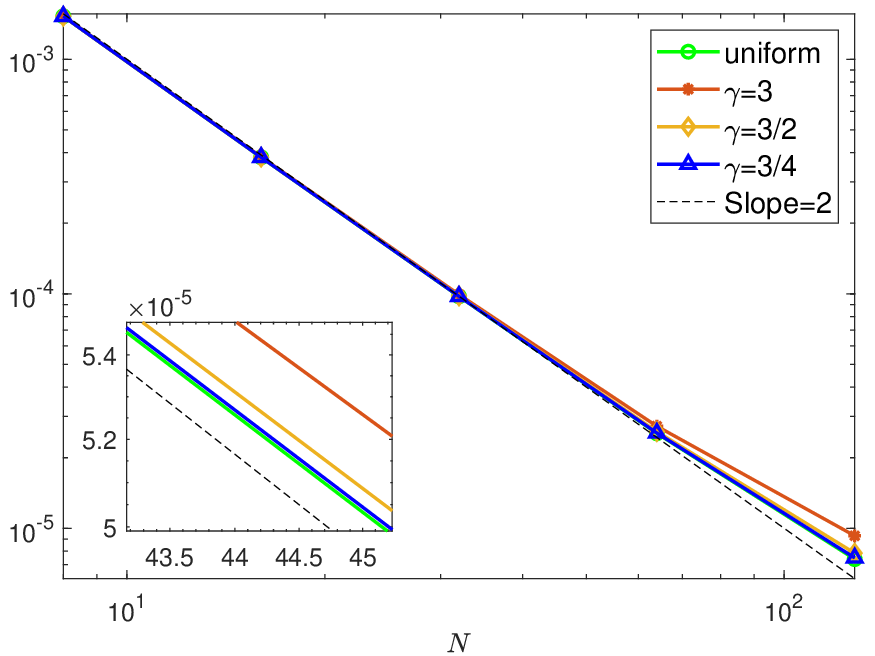} \\
  			\rotatebox{90}{{\tt $Err_2$}} &
  			\includegraphics[width=2.0in,height=1.9in]{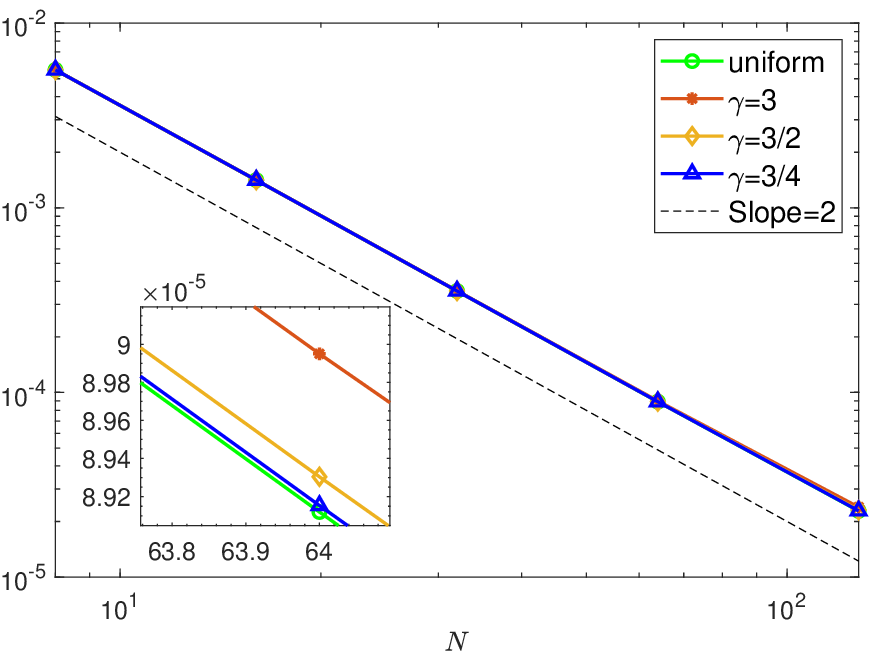} &
  			\includegraphics[width=2.0in,height=1.9in]{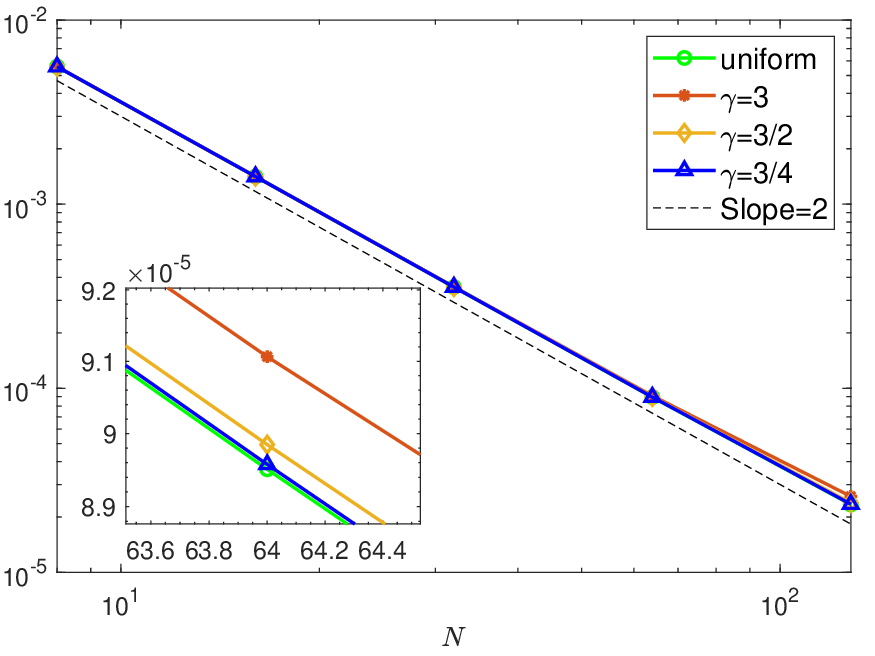} &
  			\includegraphics[width=2.0in,height=1.9in]{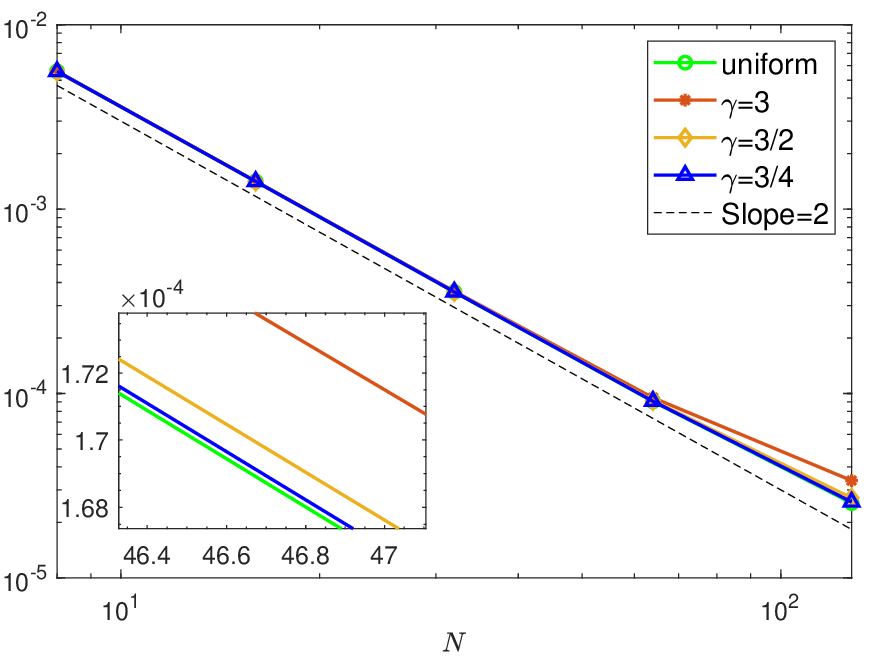} \\
  		\end{tabular}
  		\caption{Observed space errors for Example 3 for different values of $\beta$ and $N$, where $M=1000$.}
  		\label{fig6}
  	\end{figure}
  	
  	{\hypertarget{example4}\noindent \textbf{Example 4.}(\cite{kwak2023unconditionally})} We consider Eq.~\eqref{eq1.1} with  the periodic  boundary condition, the nonlinear term $f(u) = K_{pq}u^p(1-u)^q$ (selecting the second form in Eq.~\eqref{eq1.2}), $\kappa = 0.001$, $\Omega = [0, 1]^2$, and $T = 0.1$.
  	 The coefficient $K_{pq}$ is determined  in Eq.~\eqref{eq1.3}. 
  	 We assume the initial  condition is $u(x,y,0) = 0.2 \cos(2\pi x) \cos(2\pi y) + 0.75$.
 Under this periodic  boundary condition, the resulting system can be derived by doing some modifications in Eq.~\eqref{eq2.8}  or Eq.~\eqref{eq3.1}. 
 More precisely,
 $ A_x, A_y$, $\bm{u}^n$ and $\bm{\xi}^n$ are replaced by 
  \begin{align*}
  	A_x = \frac{1}{h_x^2}
  	\begin{bmatrix}
  		-2 & 1 &        &   & 1 \\
  		1  & -2 & 1     &   &   \\
  		& \ddots & \ddots & \ddots & \\
  		&        & 1 & -2 & 1 \\
  		1  &        &   & 1  & -2
  	\end{bmatrix} \in \mathbb{R}^{N_x \times N_x}, 
  	&&
  	A_y = \frac{1}{h_y^2}
  	\begin{bmatrix}
  		-2 & 1 &        &   & 1 \\
  		1  & -2 & 1     &   &   \\
  		& \ddots & \ddots & \ddots & \\
  		&        & 1 & -2 & 1 \\
  		1  &        &   & 1  & -2
  	\end{bmatrix} \in \mathbb{R}^{N_y \times N_y}.
  \end{align*}
 $\bm{u}^n = [u_{1,1}^n, \ldots, u_{N_x+1,1}^n,\ldots, u_{1,N_y+1}^n, \ldots, u_{N_x+1,N_y+1}^n]^{\top} $ and $\bm{\xi}^n=0~(n=0,\ldots,M)$.
  Since $A_x$ and $A_y$ are circulant matrices, Algorithm \ref{alg2.1} 
  needs the following changes:
  \begin{enumerate}
  	\item [(a)]
  	\begin{align*}
  		\lambda_x &= \dfrac{-4}{\ h_x^2} \left[0, \sin^2 \frac{\pi}{2N_x},\ \sin^2 \frac{2\pi}{2N_x},\ \ldots,\ \sin^2 \frac{(N_x - 1)\pi}{2N_x} \right], \\
  		\lambda_y &= \dfrac{-4}{\ h_y^2} \left[0, \sin^2 \frac{\pi}{2N_y},\ \sin^2 \frac{2\pi}{2N_y},\ \ldots,\ \sin^2 \frac{(N_y - 1)\pi}{2N_y} \right];
  	\end{align*}
  	\item [(b)]	Steps 4 and 6 are replaced by the two-dimensional Fast Fourier Transform, 
  	that is, \texttt{fft2} and \texttt{ifft2} in MATLAB.
  \end{enumerate}
  	  	
     The modified two-dimensional Fisher-KPP equation currently has no known analytical solution. 
     Therefore, the numerical solution with uniform time steps of $M = N =1024$ is adopted as the reference solution.
     In this example, we fix $p=2$.
  	  
  	  	\begin{figure}[htbp]
  	  		\setlength{\tabcolsep}{0.2pt}
  	  		\centering
  	  		\begin{tabular}{m{0.4cm}<{\centering} m{5.18cm}<{\centering} m{5.18cm}<{\centering} m{5.18cm}<{\centering}}
  	  			& $\beta=\sqrt2$ & $\beta=2$ & $\beta=\pi$ \\
  	  			\rotatebox{90}{{\tt $Err_\infty$}} &
  	  			\includegraphics[width=2.0in,height=1.9in]{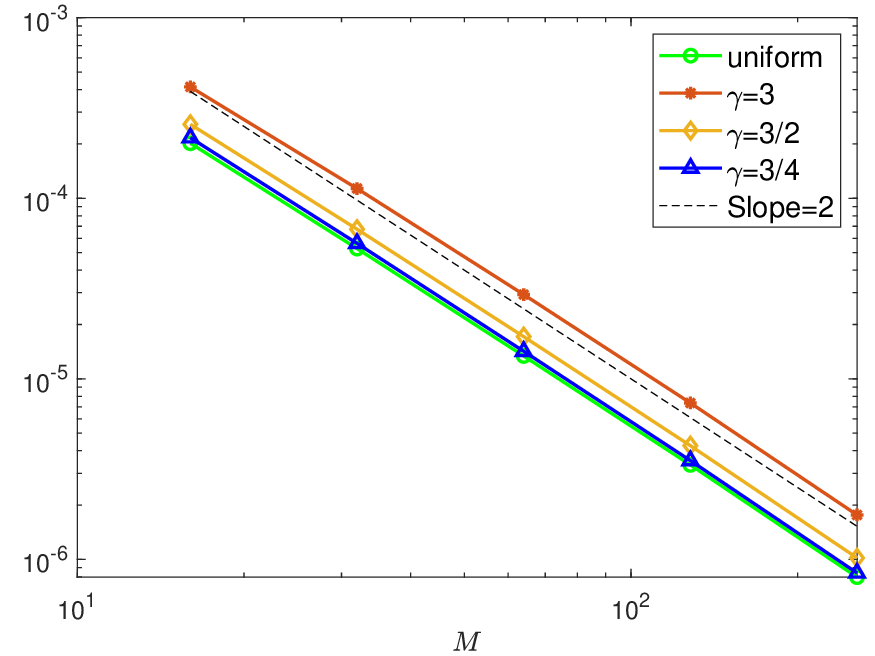} &
  	  			\includegraphics[width=2.0in,height=1.9in]{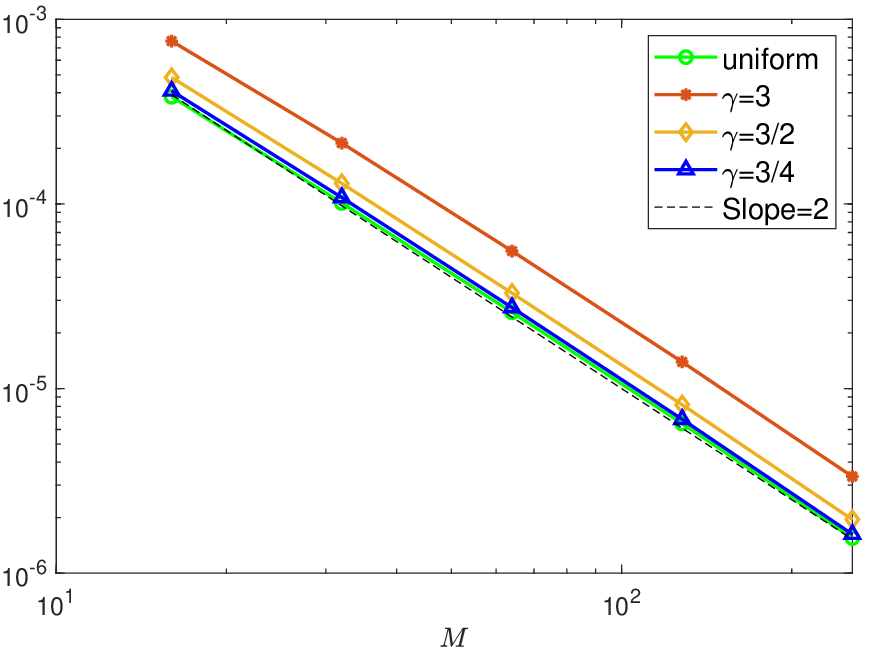} &
  	  			\includegraphics[width=2.0in,height=1.9in]{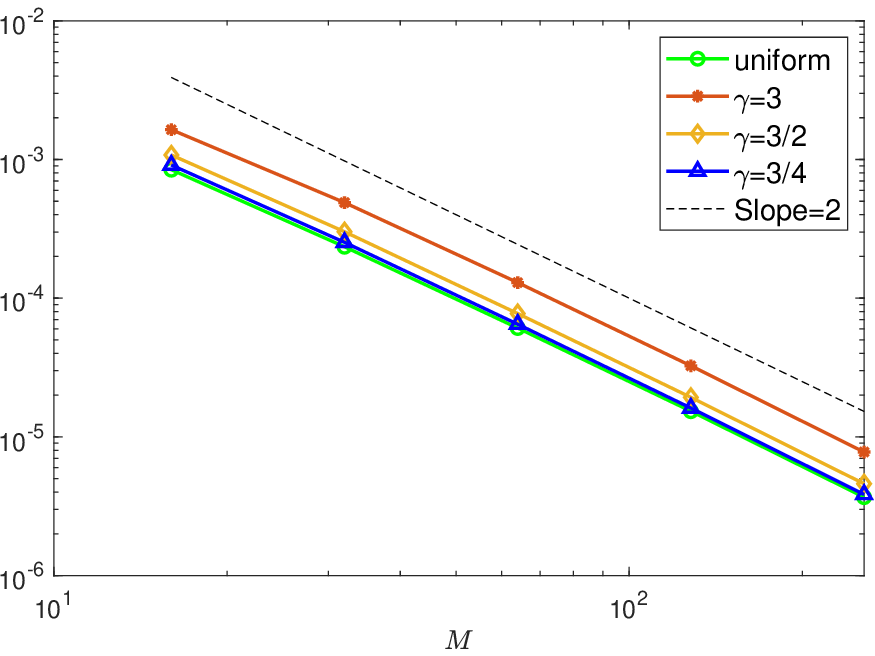} \\
  	  			\rotatebox{90}{{\tt $Err_2$}} &
  	  			\includegraphics[width=2.0in,height=1.9in]{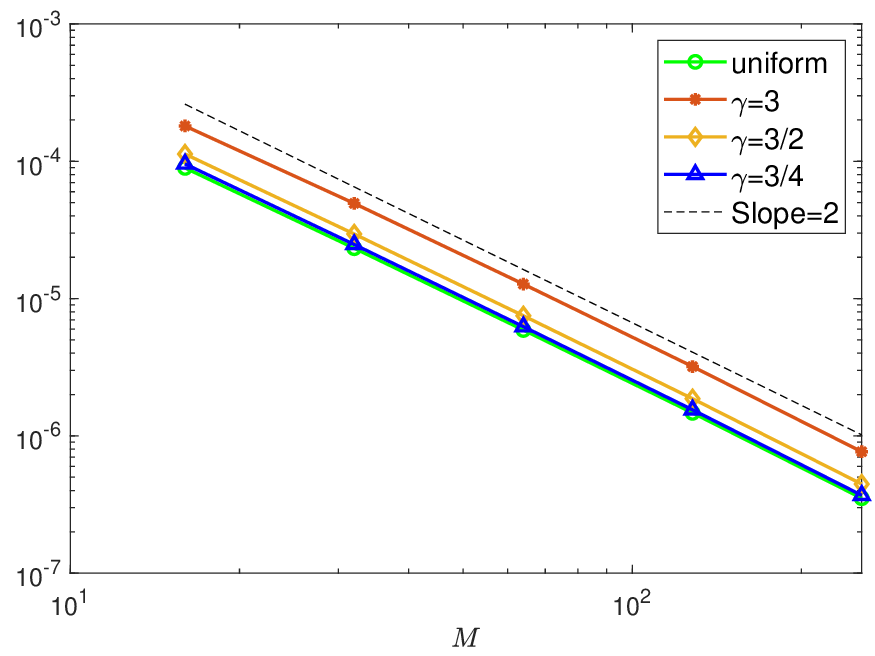} &
  	  			\includegraphics[width=2.0in,height=1.9in]{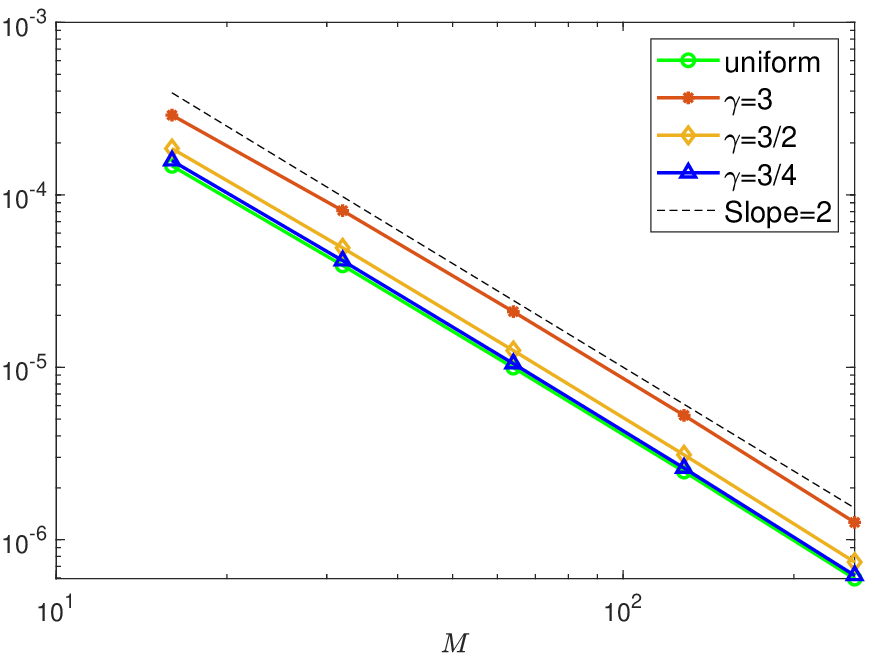} &
  	  			\includegraphics[width=2.0in,height=1.9in]{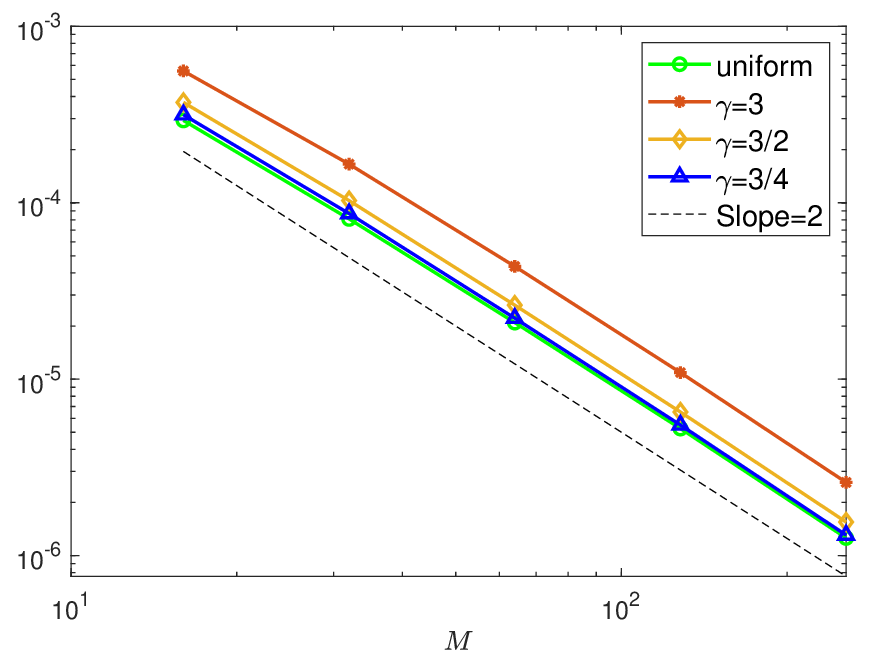} \\
  	  		\end{tabular}
  	  		\caption{Observed temporal errors for Example 4 for different values of $\beta$ and $M$, where $N=1024$ and  $q=2$.}
  	  		\label{fig7}
  	  	\end{figure}
  	  	 	\cref{fig7} demonstrates that for \hyperlink{example4 }{Example 4}  $ (q=2)$, the proposed numerical method achieves the theoretical second-order convergence in time. 
  	  	 	Under different parameters $\beta$, the error curves closely follow the reference line, confirming the effectiveness and robustness of the temporal discretization.
  	  	\begin{figure}[htbp]
  	  	\setlength{\tabcolsep}{0.2pt}
  	  	\centering
  	  	\begin{tabular}{m{0.4cm}<{\centering} m{5.18cm}<{\centering} m{5.18cm}<{\centering} m{5.18cm}<{\centering}}
  	  		& $\beta=\sqrt2$ & $\beta=2$ & $\beta=\pi$ \\
  	  		\rotatebox{90}{{\tt $Err_\infty$}} &
  	  		\includegraphics[width=2.0in,height=1.9in]{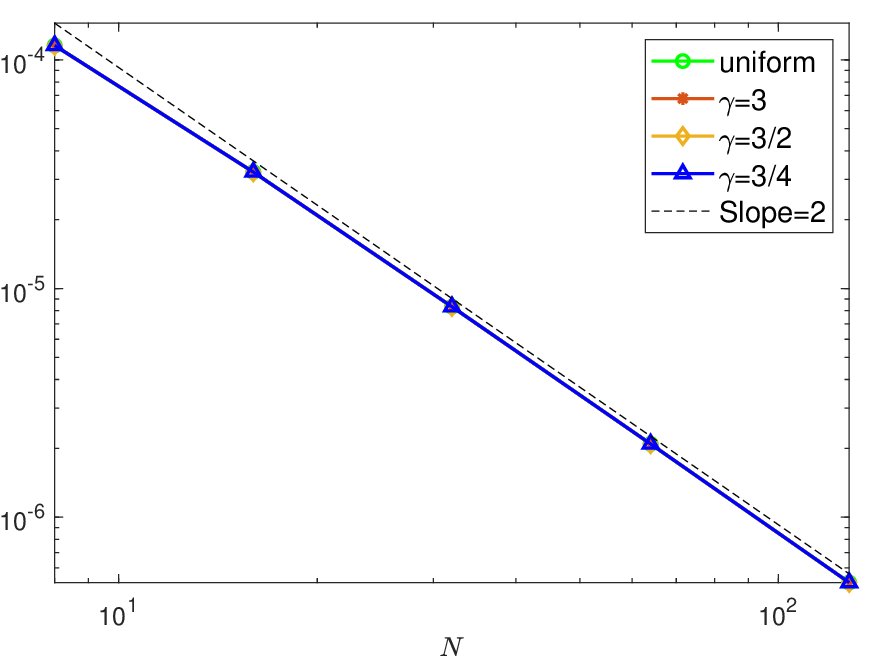} &
  	  		\includegraphics[width=2.0in,height=1.9in]{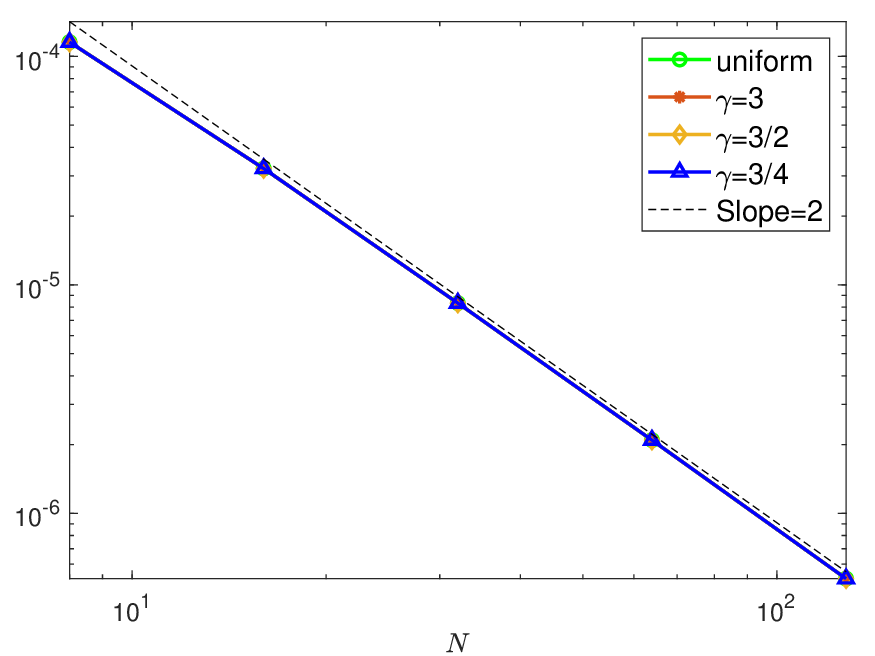} &
  	  		\includegraphics[width=2.0in,height=1.9in]{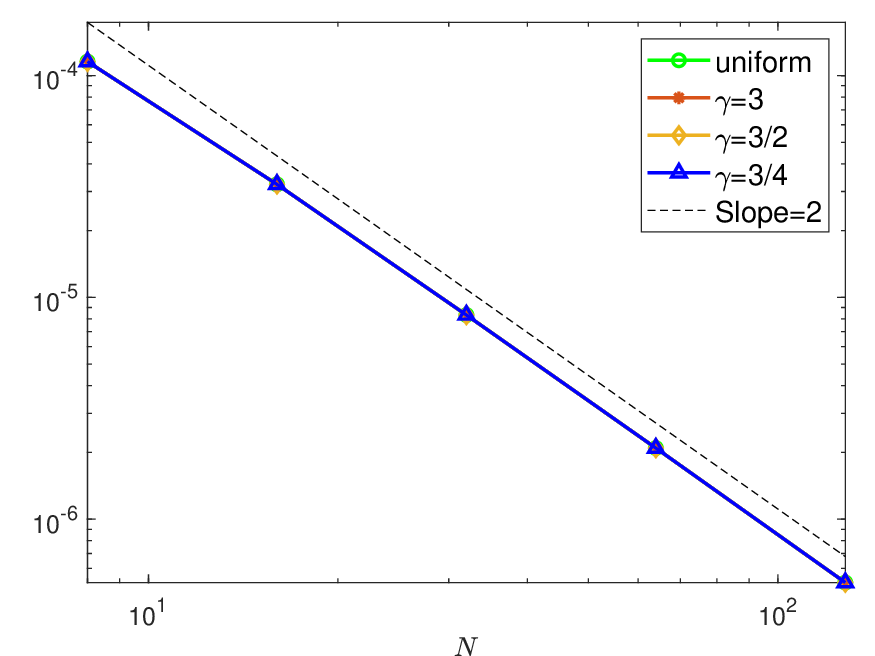} \\
  	  		\rotatebox{90}{{\tt $Err_2$}} &
  	  		\includegraphics[width=2.0in,height=1.9in]{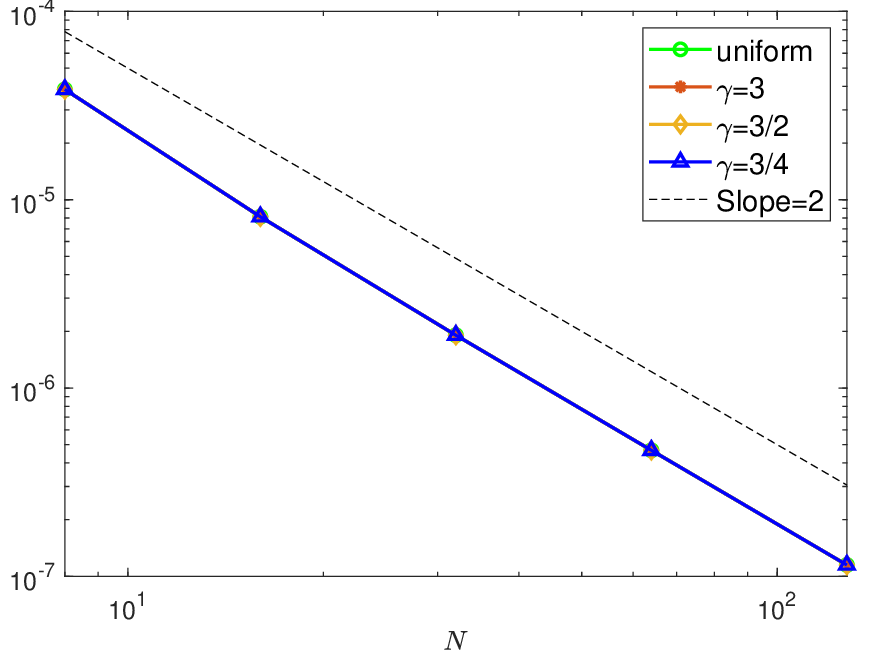} &
  	  		\includegraphics[width=2.0in,height=1.9in]{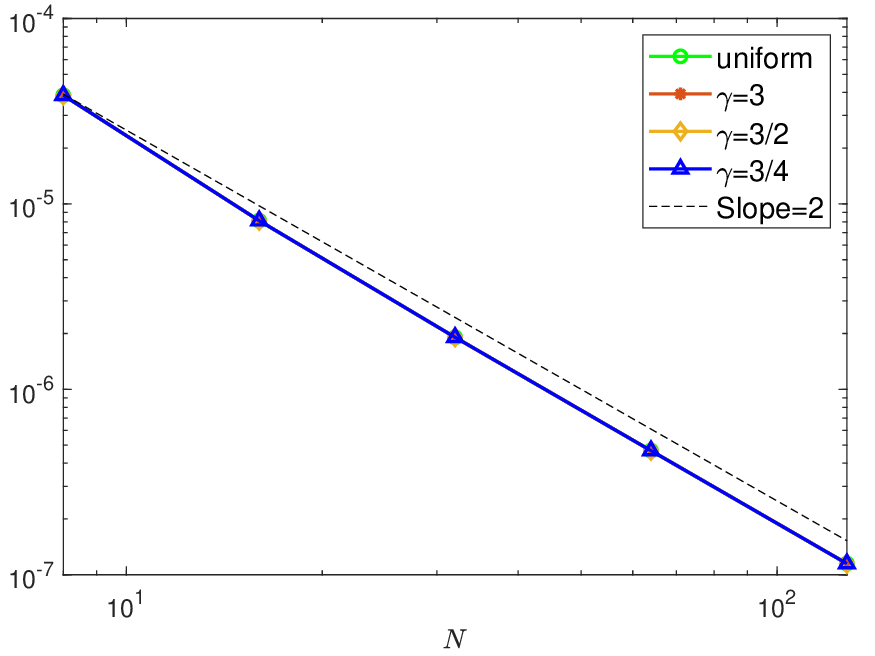} &
  	  		\includegraphics[width=2.0in,height=1.9in]{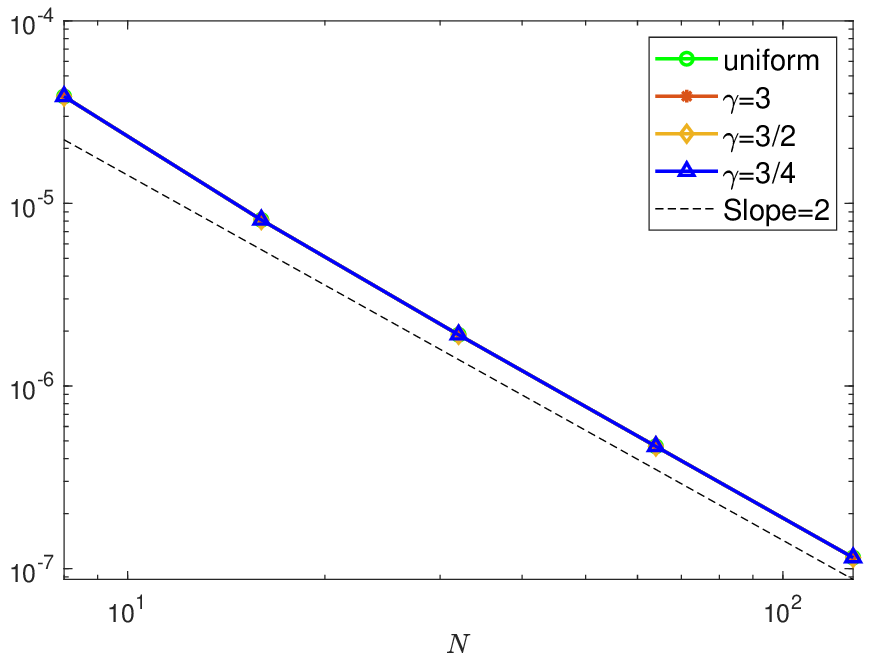} \\
  	  	\end{tabular}
  	  	\caption{Observed space errors for Example 4 for different values of $\beta$ and $N$, where $M=1024$  and $q=2$.}
  	  	\label{fig8}
  	  \end{figure}
  	  
  		\cref{fig8}  shows that for \hyperlink{example4 }{Example 4}  $ (q=2)$, the scheme maintains second-order convergence in space.
  		All error curves align parallel to the reference line, verifying the correctness of the central difference scheme, with convergence independent of the temporal parameter $\beta$.
  	  
  	  	\begin{figure}[htbp]
  	  	\setlength{\tabcolsep}{0.2pt}
  	  	\centering
  	  	\begin{tabular}{m{0.4cm}<{\centering} m{5.18cm}<{\centering} m{5.18cm}<{\centering} m{5.18cm}<{\centering}}
  	  		& $\beta=\sqrt2$ & $\beta=2$ & $\beta=\pi$ \\
  	  		\rotatebox{90}{{\tt $Err_\infty$}} &
  	  		\includegraphics[width=2.0in,height=1.9in]{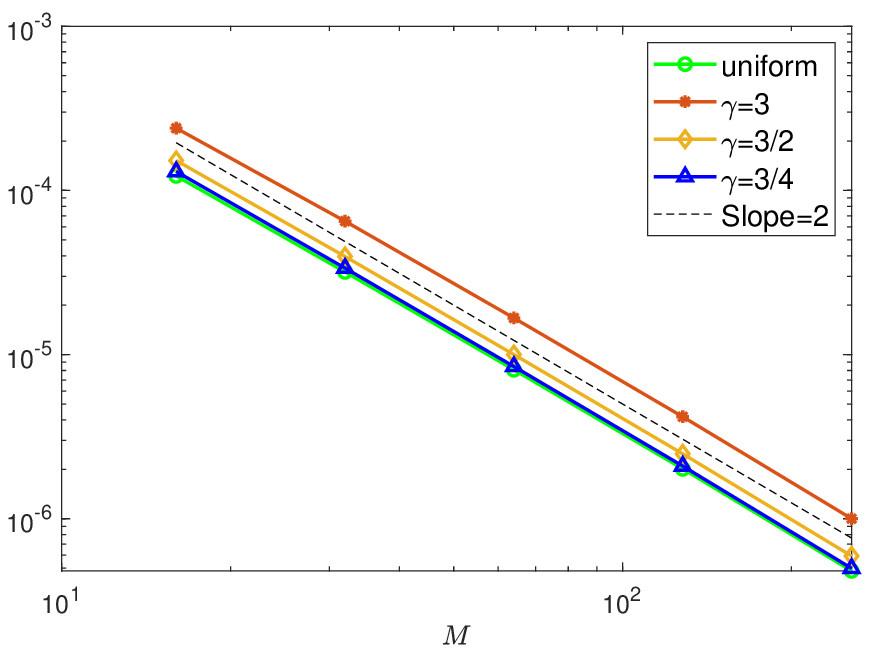} &
  	  		\includegraphics[width=2.0in,height=1.9in]{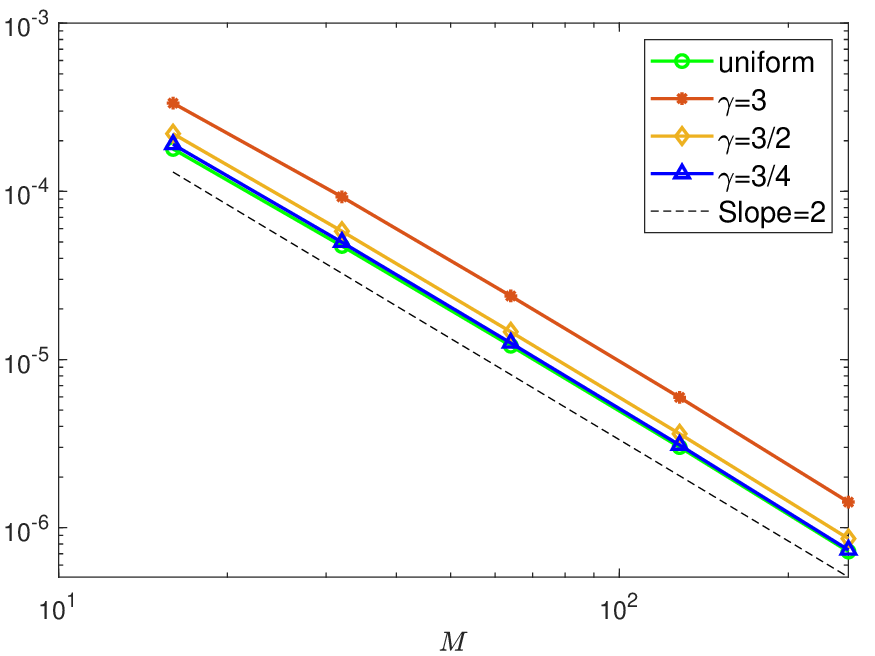} &
  	  		\includegraphics[width=2.0in,height=1.9in]{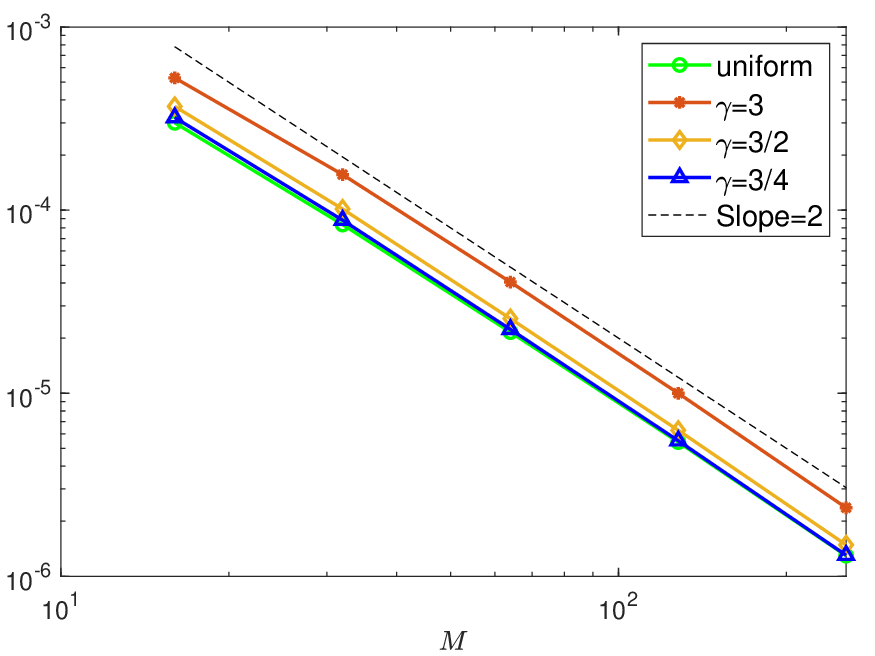} \\
  	  		\rotatebox{90}{{\tt $Err_2$}} &
  	  		\includegraphics[width=2.0in,height=1.9in]{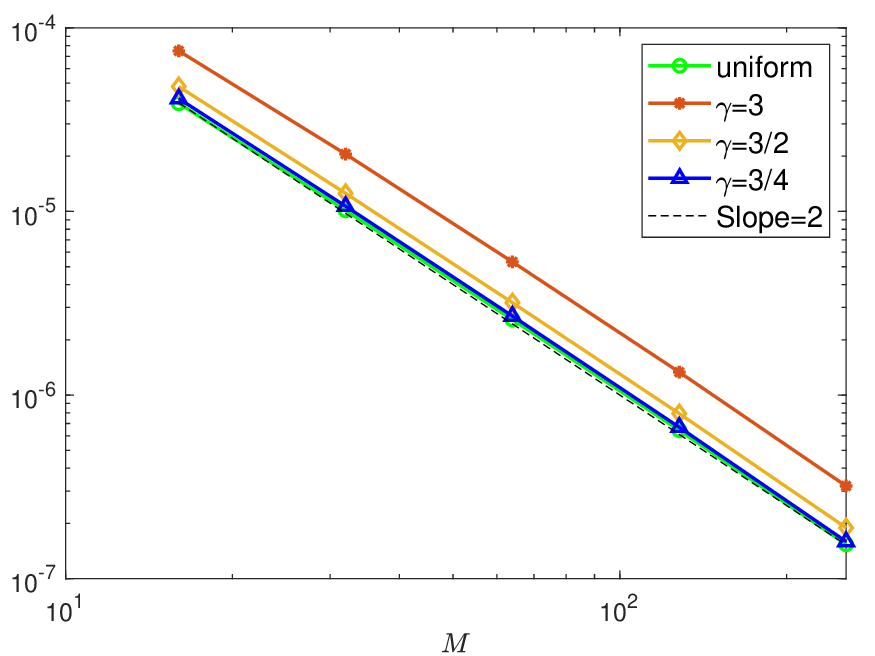} &
  	  		\includegraphics[width=2.0in,height=1.9in]{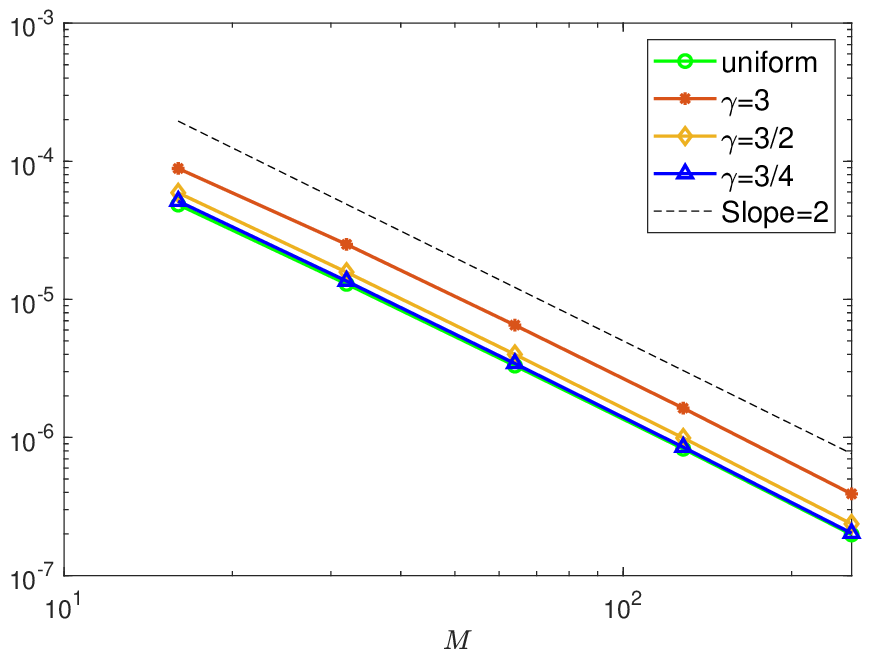} &
  	  		\includegraphics[width=2.0in,height=1.9in]{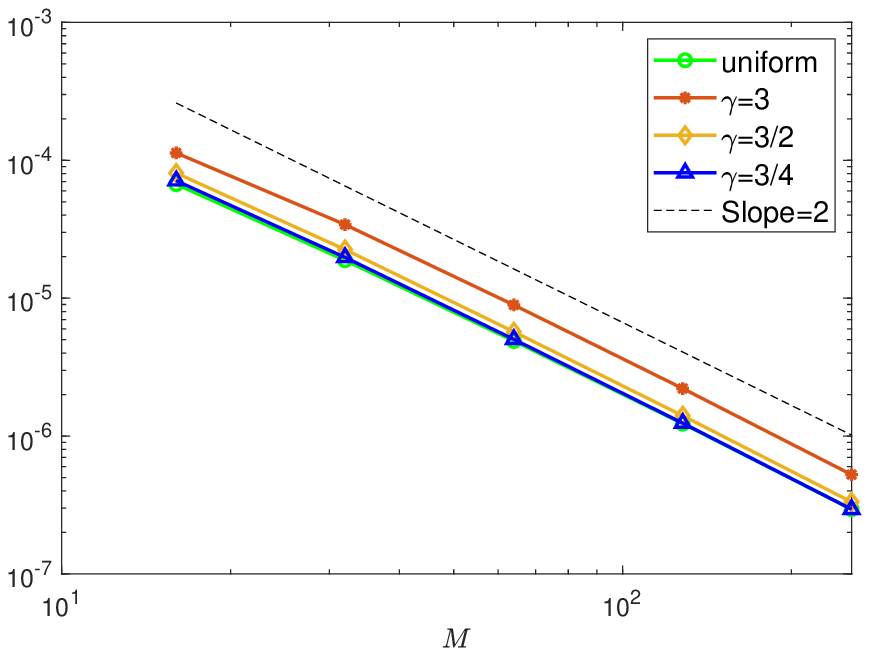} \\
  	  	\end{tabular}
  	  	\caption{Observed temporal errors for Example 4 for different values of $\beta$ and $M$, where $N=1024$  and $q=3$.}
  	  	\label{fig9}
  	  \end{figure}

  	  \begin{figure}[htbp]
  	  	\setlength{\tabcolsep}{0.2pt}
  	  	\centering
  	  	\begin{tabular}{m{0.4cm}<{\centering} m{5.18cm}<{\centering} m{5.18cm}<{\centering} m{5.18cm}<{\centering}}
  	  		& $\beta=\sqrt2$ & $\beta=2$ & $\beta=\pi$ \\
  	  		\rotatebox{90}{{\tt $Err_\infty$}} &
  	  		\includegraphics[width=2.0in,height=1.9in]{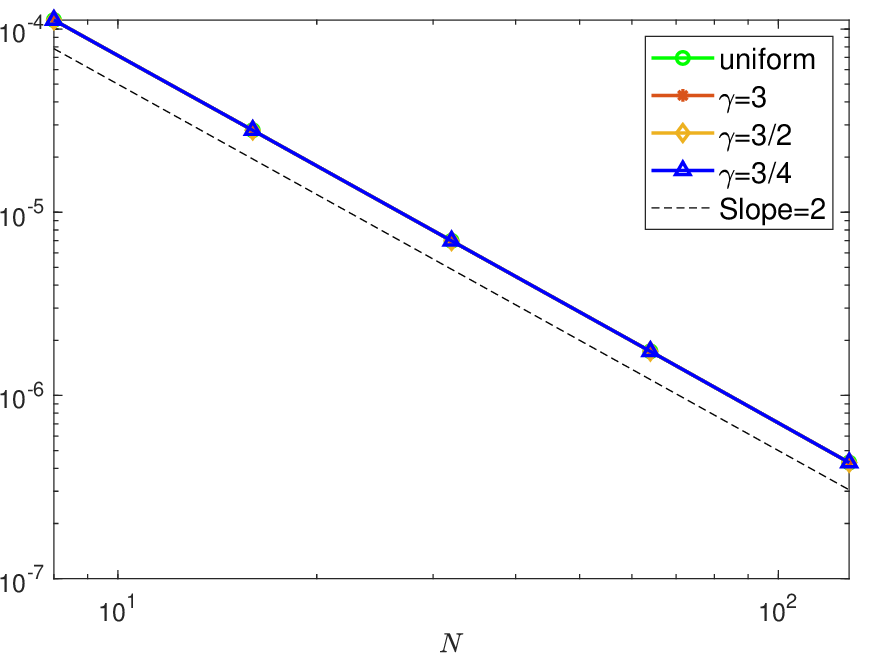} &
  	  		\includegraphics[width=2.0in,height=1.9in]{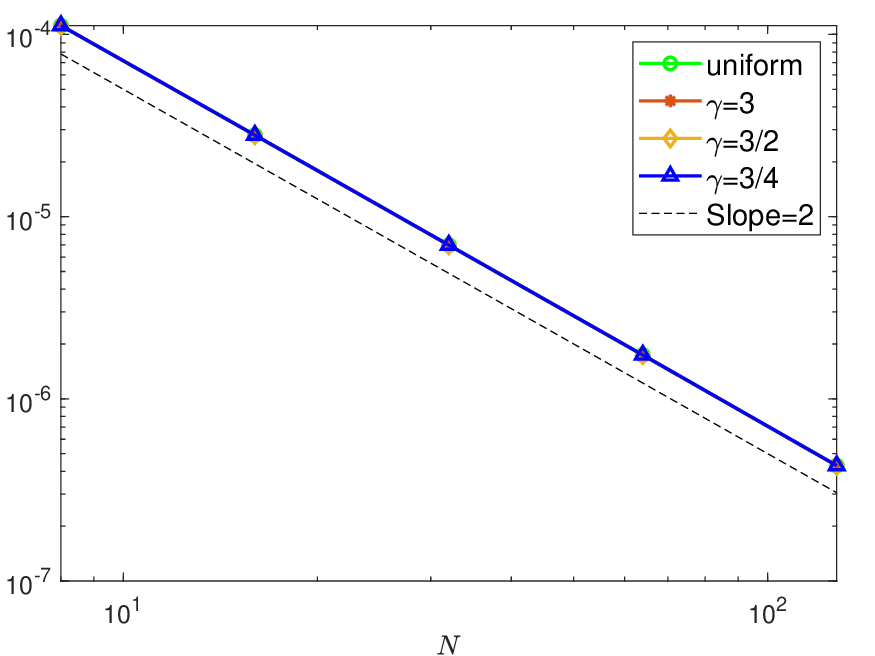} &
  	  		\includegraphics[width=2.0in,height=1.9in]{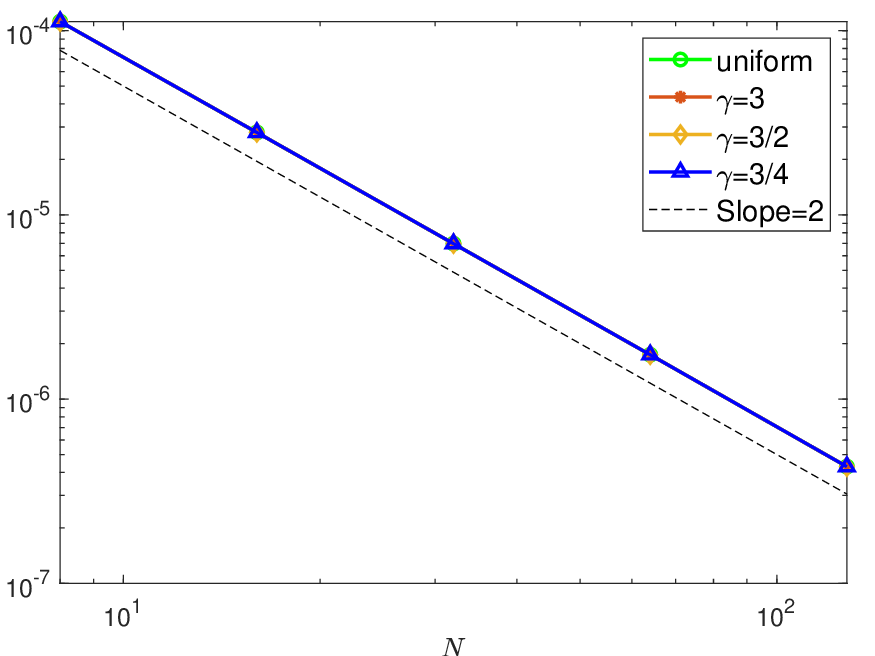} \\
  	  		\rotatebox{90}{{\tt $Err_2$}} &
  	  		\includegraphics[width=2.0in,height=1.9in]{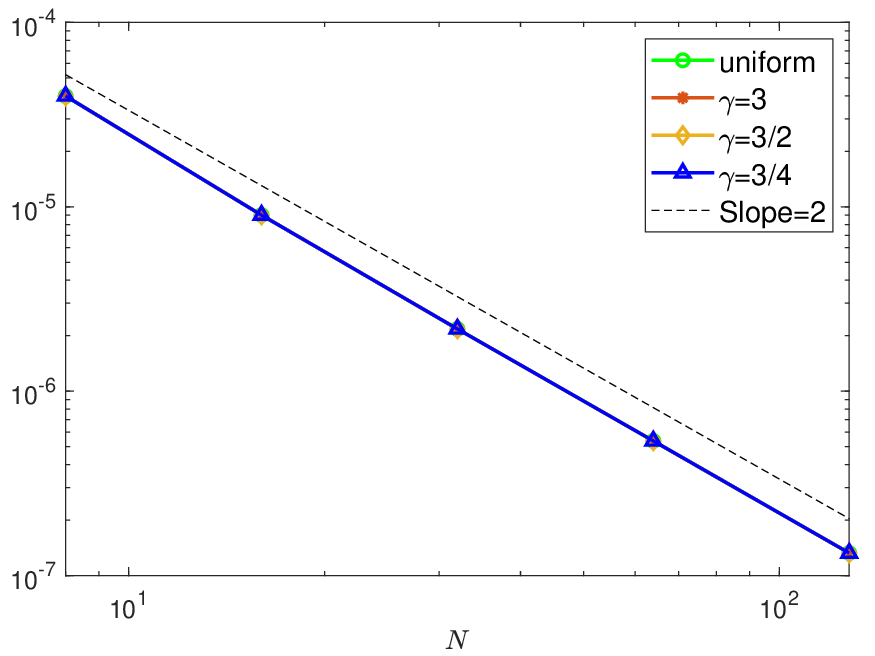} &
  	  		\includegraphics[width=2.0in,height=1.9in]{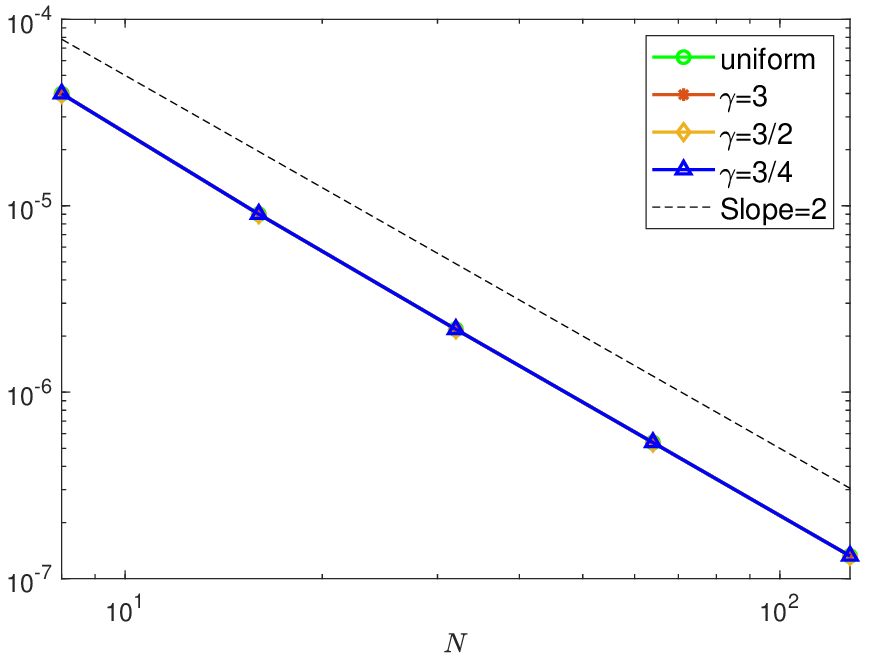} &
  	  		\includegraphics[width=2.0in,height=1.9in]{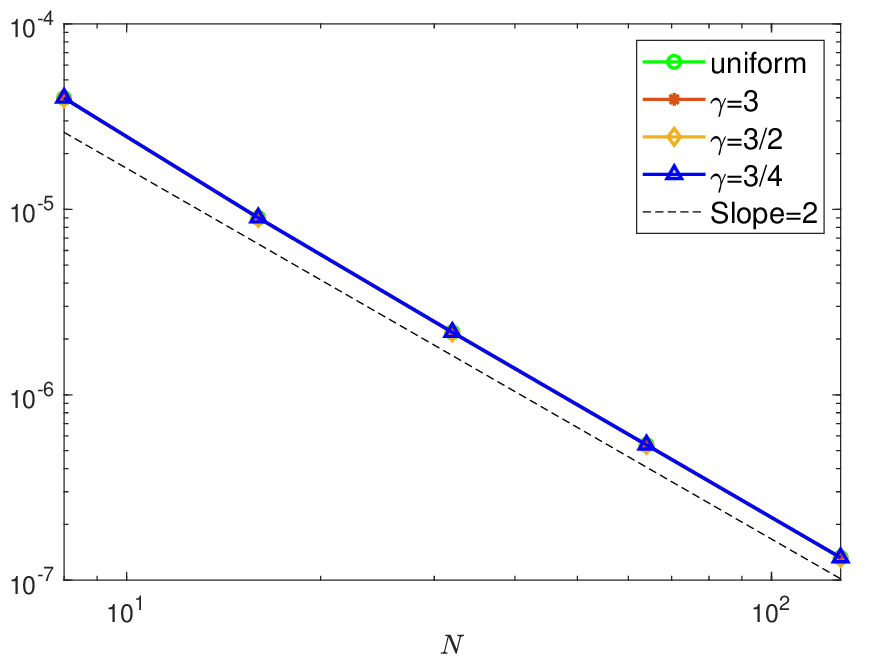} \\
  	  	\end{tabular}
  	  	\caption{Observed space errors for Example 4 for different values of $\beta$ and $N$, where $M=1024$ and $q=3$.}
  	  	\label{fig10}
  	  \end{figure}

  \begin{figure}[htbp]
  	\centering
  	\setlength{\tabcolsep}{4pt}
  	\begin{tabular}{m{0.08\textwidth} >{\centering\arraybackslash}m{0.35\textwidth} >{\centering\arraybackslash}m{0.35\textwidth}}
  		& $(p,q)=(2,2)$ & $(p,q)=(2,3)$ \\
  		\rotatebox{90}{\makebox[\linewidth]{Uniform}} &
  		\includegraphics[width=\linewidth]{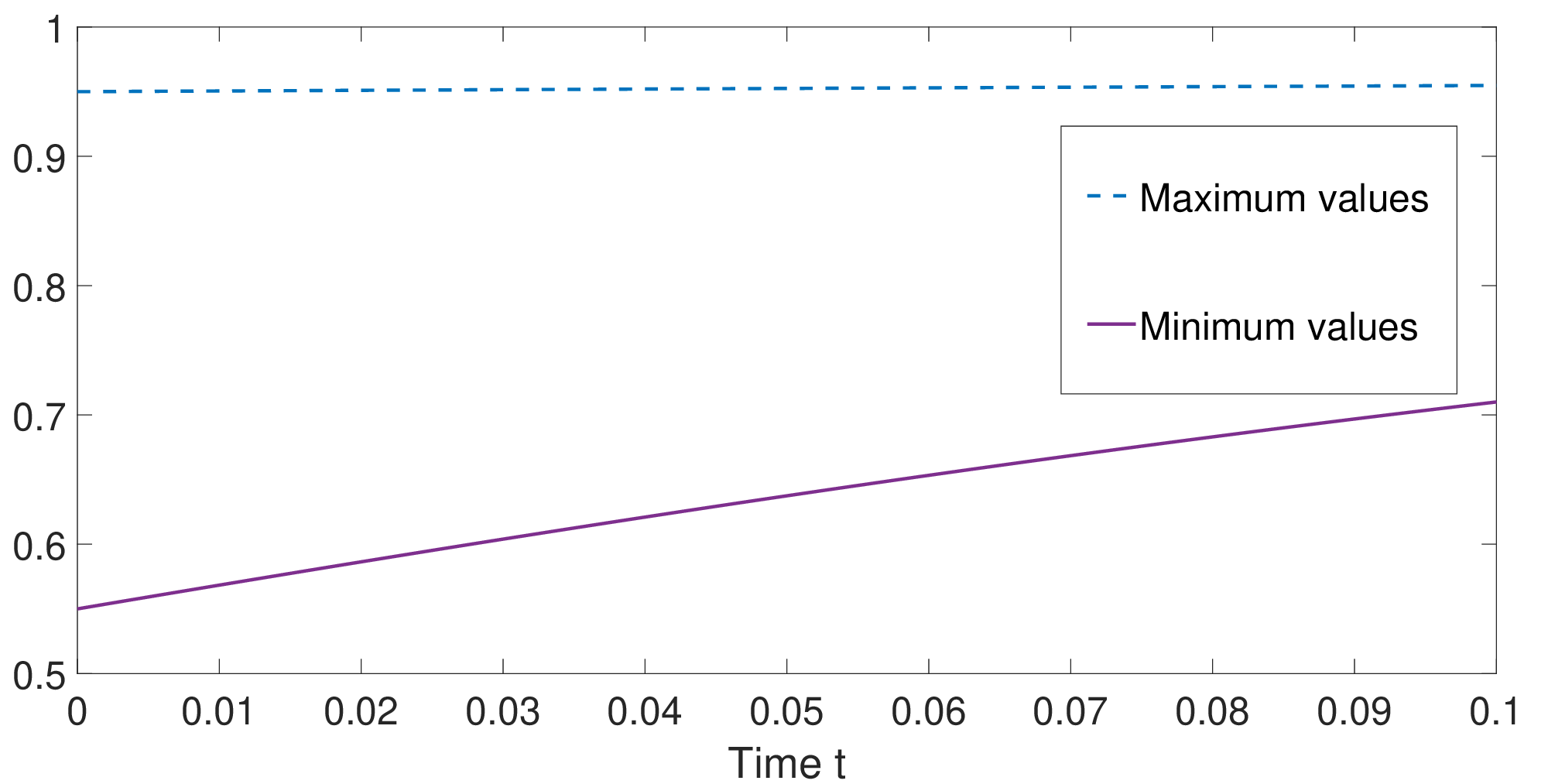} &
  		\includegraphics[width=\linewidth]{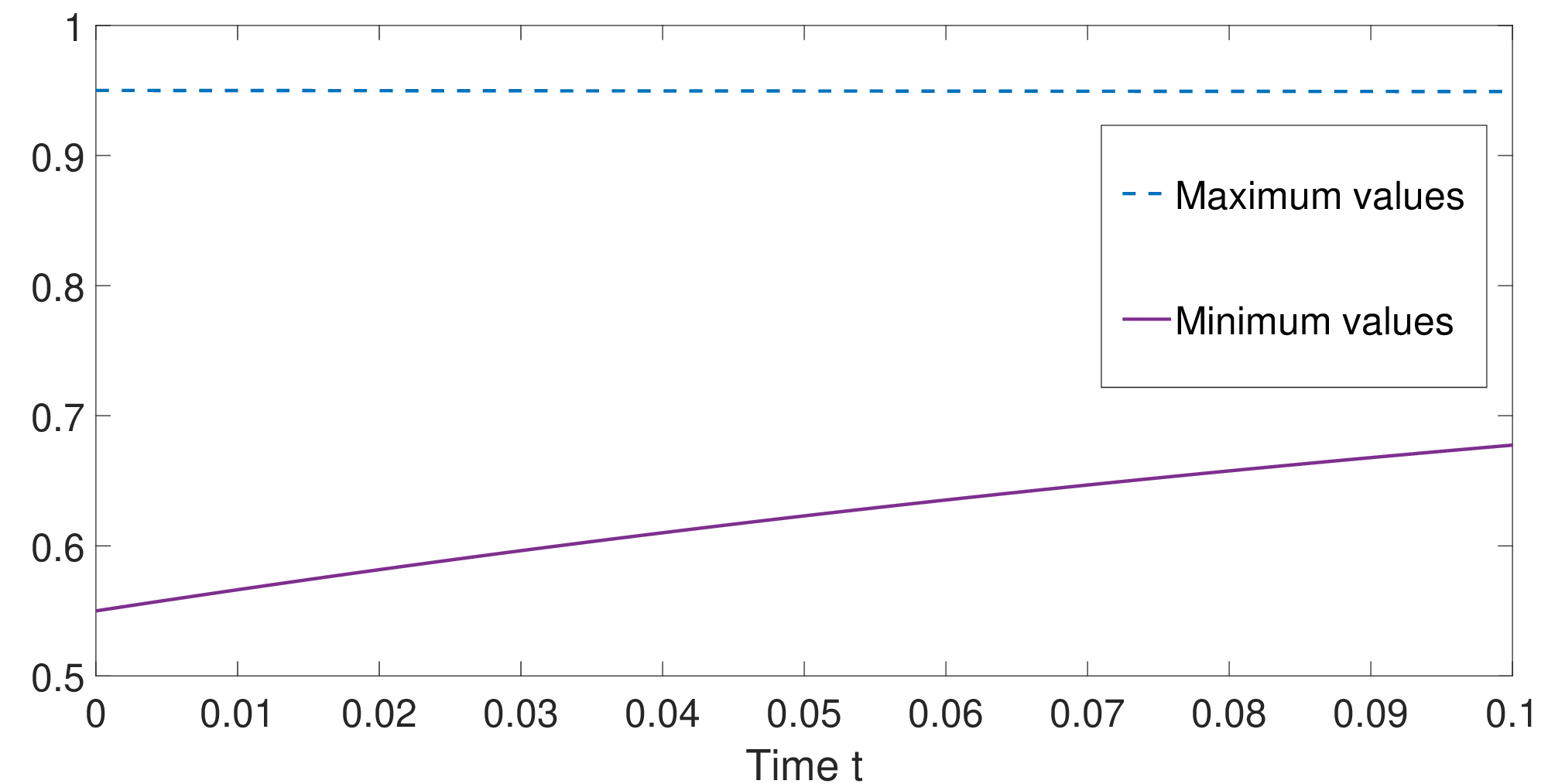} \\
  		\rotatebox{90}{\makebox[\linewidth]{Nonuniform}} &
  		\includegraphics[width=\linewidth]{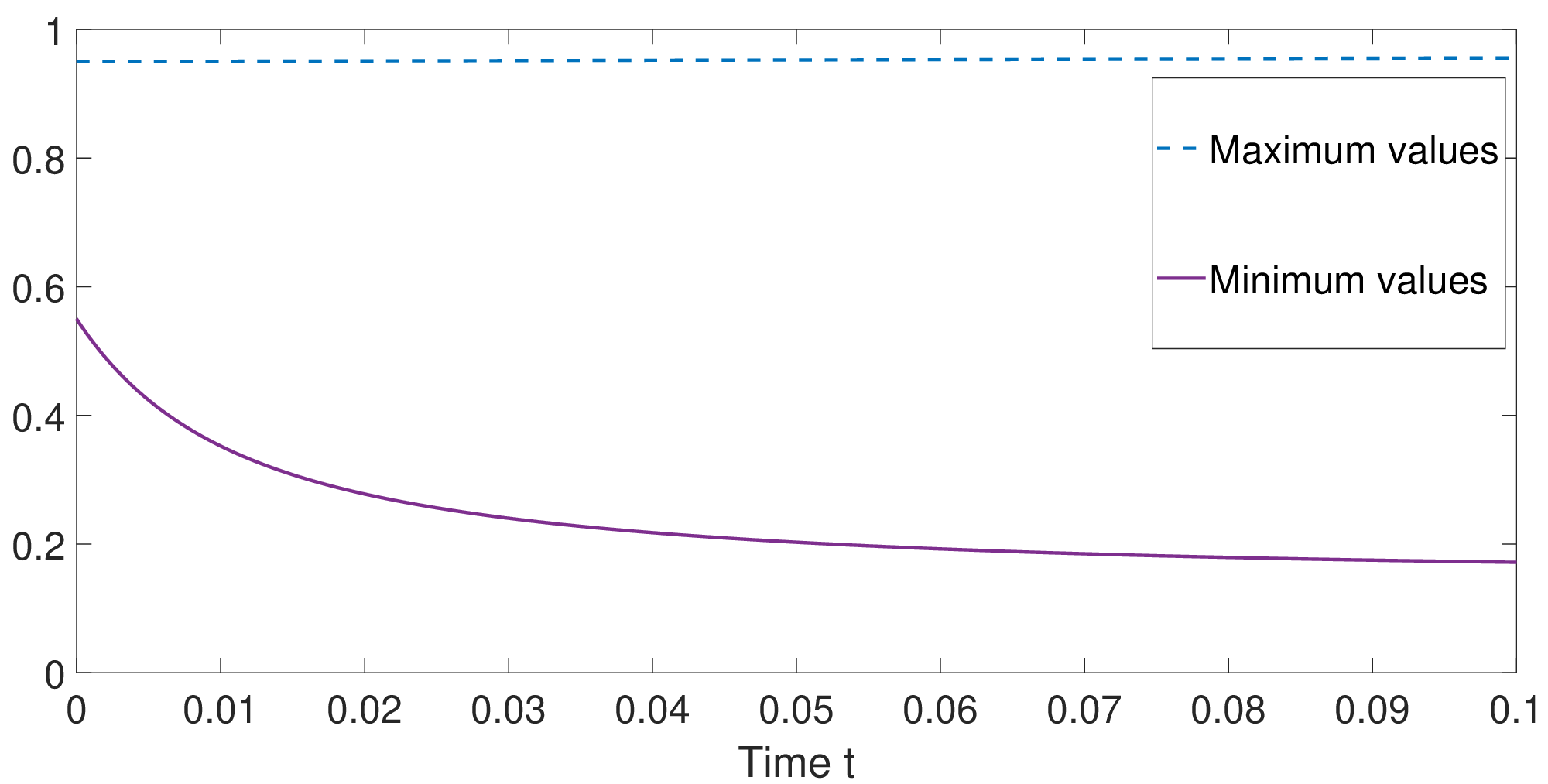} &
  		\includegraphics[width=\linewidth]{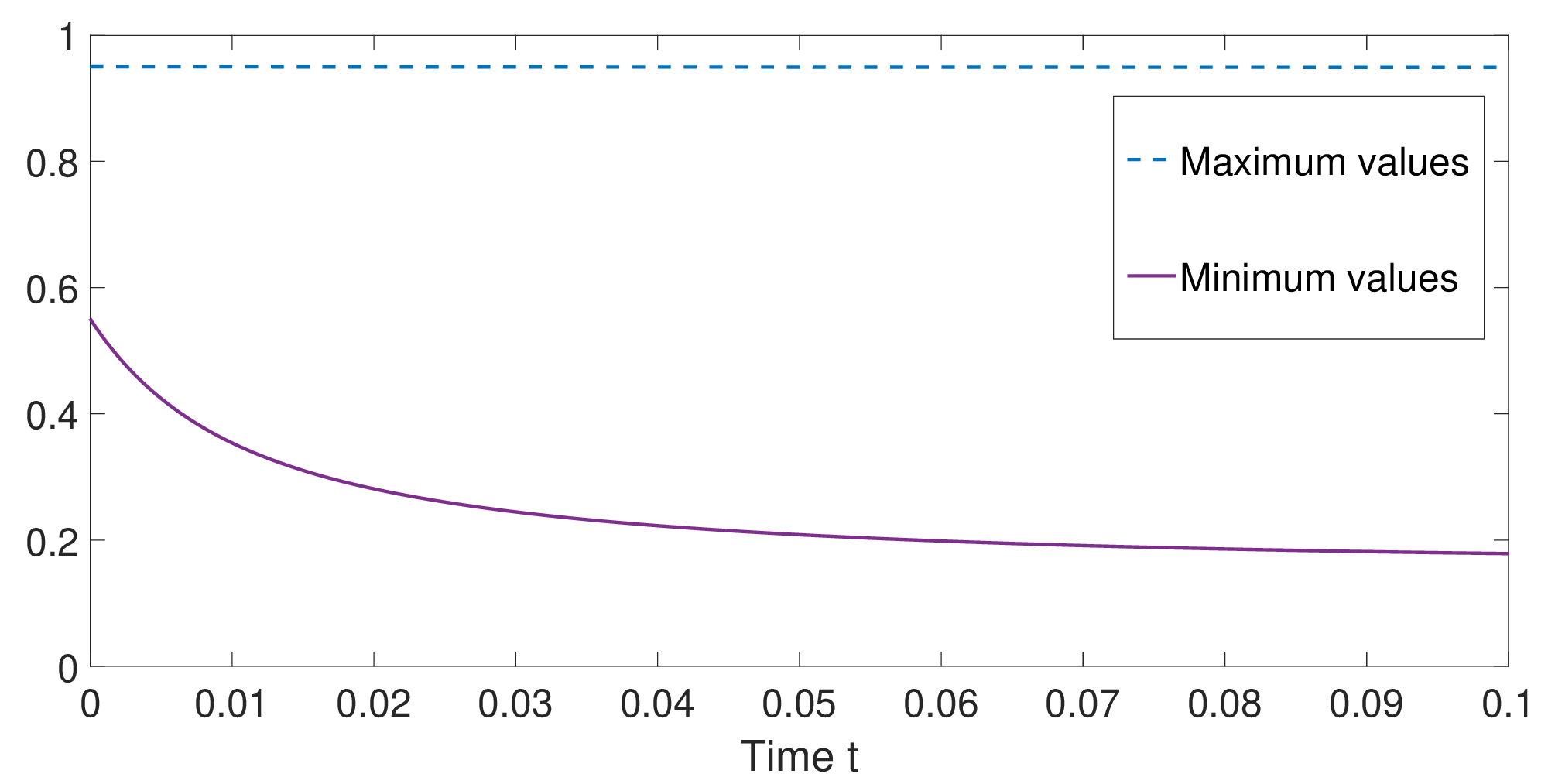} \\
  	\end{tabular}
  	\caption{Observed numerical solution maximum and minimum value for Example 4 with different parameter pairs $(p,q)$.}
  	\label{fig11}
  \end{figure}
  
  	  	\cref{fig9} verifies the second-order temporal convergence for \hyperlink{example4 }{Example 4} with the modified nonlinear terms  $ (q=3)$. 
  	  	The error curves perfectly follow the second-order slope, proving the scheme's adaptability to variations in equation parameters.
  	  
  	  	\cref{fig10} confirms second-order spatial convergence for \hyperlink{example4 }{Example 4}  $ (q=3)$. The overlapping curves for different $\beta$ values further validate the stability and reliability of the spatial discretization, whose accuracy is independent of both temporal strategy and nonlinear parameters.
  	  	
  	  	In summary, the results from  \cref{fig7}-\cref{fig10} indicate that the GBDF2-IMEX scheme maintains second-order convergence accuracy in both time and space when solving the  modified Eq.~\cref{eq1.1} with different parameters in the nonlinear term, demonstrating excellent robustness.
  	  
  	  In \cref{fig11}, $\beta=2$ for all cases and $\gamma=2/3$ for nonuniform time steps.
  	  The plots generated with other values of $\beta$ and $\gamma$ are similar and thus are omitted here for brevity.
  	  
  	\cref{fig11} clearly demonstrates that the GBDF2-IMEX scheme  can stably simulate the dynamic behavior of the  two-dimensional modified Fisher-KPP equation, regardless of whether uniform or nonuniform time steps are employed. 
  	The maximum and minimum values of the numerical solution maintain smooth variation throughout the temporal evolution, with no unphysical oscillations or divergence phenomena observed, which confirms the scheme's excellent numerical stability.
  	  
\section{Concluding remarks}
\label{sec5}
This paper focuses on a uniform/nonuniform GBDF2-IMEX scheme for  ~\cref{eq1.1} based on the Taylor expansion at time $t^{n+\beta}$, where $\beta >1$ is a tunable parameter.
A rigorous theoretical analysis of the uniform time-stepping scheme is provided, proving its stability and convergence.
 For  Eq.~\cref{eq1.1}, numerical experiments confirm its second-order convergence in both time and space.
 A key contribution of this work is the development of an alternative second-order scheme that remains effective on nonuniform time grids.

Future work will focus on  the rigorous  error analysis of  the nonuniform GBDF2-IMEX scheme for ~\cref{eq3.1}.
In the future, the effectiveness of other nonuniform time-step formats may surpass that of the nonuniform time-step format proposed in this paper, which remains to be explored and developed.







\bibliographystyle{siamplain}
\bibliography{Ref}
\end{document}